\documentclass[referee,sn-basic,Numbered]{sn-jnl}

\usepackage{graphicx}%
\usepackage{multirow}%
\usepackage{amsmath,amssymb,amsfonts}%
\usepackage{amsthm}%
\usepackage{mathrsfs}%
\usepackage[title]{appendix}%
\usepackage{xcolor}%
\usepackage{textcomp}%
\usepackage{manyfoot}%
\usepackage{booktabs}%
\usepackage{algorithm}%
\usepackage{algorithmicx}%
\usepackage{algpseudocode}%
\usepackage{listings}%

\usepackage{bbold}
\usepackage{dsfont}
\usepackage{mathtools}
\usepackage{stmaryrd}
\usepackage{verbatim}
\usepackage{titlesec}
\usepackage{mathrsfs}

\DeclareFontFamily{U}{mathx}{}
\DeclareFontShape{U}{mathx}{m}{n}{<-> mathx10}{}
\DeclareSymbolFont{mathx}{U}{mathx}{m}{n}
\DeclareMathAccent{\widehat}{0}{mathx}{"70}
\DeclareMathAccent{\widecheck}{0}{mathx}{"71}

\usepackage{float}
\usepackage{comment}

\usepackage{mhchem}
\graphicspath{ {images/} }
\usepackage{array}
\newcolumntype{P}[1]{>{\centering\arraybackslash}p{#1}}
\usepackage{longtable}
\usepackage{mathrsfs}
\usepackage{caption}
\DeclareCaptionLabelFormat{adja-page}{\hrulefill\\#1 #2 \emph{(previous page)}}
\usepackage{subcaption}

\allowdisplaybreaks

\usepackage{tikz}

\usepackage{tabulary}
\usepackage{booktabs}
\usepackage{array} % For \extracolsep
\usepackage{blkarray}

\usepackage{hyperref}
\hypersetup{
    colorlinks=true,
    linkcolor=blue,
    filecolor=blue,      
    urlcolor=blue,
    citecolor=blue,
    pdftitle={Overleaf Example},
    pdfpagemode=FullScreen,
    }

\usepackage[capitalise,nameinlink]{cleveref}

\usepackage[T3,T1]{fontenc}
\DeclareSymbolFont{tipa}{T3}{cmr}{m}{n}
\DeclareMathAccent{\invbreve}{\mathalpha}{tipa}{16}

\usepackage{soul,xcolor}
\setstcolor{red}

\crefname{supp}{Supplement}{Supplements}

\usepackage{amsthm}
\newtheorem{theorem}{Theorem}[section]
\newtheorem{corollary}{Corollary}[section]
\newtheorem{lemma}{Lemma}[section]

\newtheorem{assumption}{Assumption}[section]
\numberwithin{equation}{section}

\theoremstyle{definition}

\theoremstyle{definition}
\newtheorem{remark}{Remark}[section]

\theoremstyle{definition}
\newtheorem{example}{Example}[section]

\newcommand{\initialx}{x^\circ}
\newcommand{\initialxbreve}{\Breve{x}^\circ}

\renewcommand{\l}{\ell}

\newcommand{\eps}{\varepsilon}

\newcommand{\Ss}{\mathscr{S}}
\newcommand{\Rs}{\mathscr{R}}

\newcommand{\PP}{\mathds{P}}
\newcommand{\R}{\mathds{R}}

\newcommand{\Z}{\mathds{Z}}
\newcommand{\E}{\mathds{E}} %Make comments

\newcommand{\G}{\mathcal{G}}

\newcommand{\F}{\mathcal{F}}
\newcommand{\B}{\mathcal{B}}

\renewcommand{\L}{\mathcal{L}}

\newcommand{\X}{\mathcal{X}}

\newcommand{\V}{\mathcal{V}}

\newcommand{\one}{\mathds{1}}
\newcommand{\zero}{\mathbb{0}}

\newcommand{\Stot}{\text{S}_{\text{tot}}}
\newcommand{\Dtot}{\text{D}_{\text{tot}}}

\newcommand{\Ntot}{\text{N}_{\text{tot}}}

\newcommand{\prp}{\Lambda}  %Propensity function
\newcommand{\rate}{\Upsilon} %Rate of transition (sum of propensity)
\newcommand{\upb}{K} %Upper bound for generic birth-and-death process

\newcommand{\edge}{e}

\newcommand{\Eu}{\mathcal{E}_+}
\newcommand{\Ed}{\mathcal{E}_-}

\newcommand{\wccNum}{p}
\newcommand{\wccIter}{q}

\DeclareMathOperator{\rank}{rank}

\DeclareMathOperator{\diag}{diag}

\newcommand{\interior}[1]{%
  {\kern0pt#1}^{\mathrm{o}}%
}

\DeclarePairedDelimiterX{\abs}[1]{\lvert}{\rvert}{#1}
\DeclarePairedDelimiterX{\norm}[1]{\lVert}{\rVert}{#1}
\DeclarePairedDelimiterX{\floor}[1]{\lfloor}{\rfloor}{#1}
\DeclarePairedDelimiterX{\qvar}[1]{\langle}{\rangle}{#1}
\DeclarePairedDelimiterX{\inn}[1]{\langle}{\rangle}{#1}

\DeclareCaptionLabelFormat{cont}{#1~#2\alph{ContinuedFloat}}
\begin{document}

\title[\textcolor{black}{Coclique} Level Structure for Stochastic Chemical Reaction Networks]{\textcolor{black}{Coclique} Level Structure for Stochastic Chemical Reaction Networks}


\author[1,2]{\fnm{Simone} \sur{Bruno}}\email{sbruno@ds.dfci.harvard.edu}
\equalcont{These authors contributed equally to this work.}

\author[3]{\fnm{Yi} \sur{Fu}}\email{yif064@ucsd.edu}
\equalcont{These authors contributed equally to this work.}

\author[4]{\fnm{Felipe A.} \sur{Campos}}\email{fcamposv@ucsd.edu}

\author[2]{\fnm{Domitilla} \sur{Del Vecchio}}\email{ddv@mit.edu}

\author*[4]{\fnm{Ruth J.} \sur{Williams}}\email{rjwilliams@ucsd.edu}

\affil[1]{\orgdiv{Department of Data Science}, \orgname{Dana-Farber Cancer Institute}, \orgaddress{\street{450 Brookline Avenue}, \city{Boston}, \postcode{02115}, \state{MA}, \country{USA}}}

\affil[2]{\orgdiv{Department of Mechanical Engineering}, \orgname{Massachusetts Institute of Technology}, \orgaddress{\street{77 Massachusetts Avenue}, \city{Cambridge}, \postcode{02139}, \state{MA}, \country{USA}}}

\affil[3]{\orgdiv{Bioinformatics and Systems Biology Program}, \orgname{University of California, San Diego}, \orgaddress{\street{9500 Gilman Drive}, \city{La Jolla}, \postcode{92093-0112}, \state{CA}, \country{USA}}}

\affil*[4]{\orgdiv{Department of Mathematics}, \orgname{University of California, San Diego}, \orgaddress{\street{9500 Gilman Drive}, \city{La Jolla}, \postcode{92093-0112}, \state{CA}, \country{USA}}}

%\affil[*]{These authors contributed equally: S. Bruno and Y. Fu}
%\date{}                     %% if you don't need date to appear

%\doublespacing


\abstract{Continuous time Markov chains are commonly used as models for the stochastic behavior of chemical reaction networks. More precisely, these Stochastic Chemical Reaction Networks (SCRNs) are frequently used to gain a mechanistic understanding of how chemical reaction rate parameters impact the stochastic behavior of these systems. One property of interest is mean first passage times (MFPTs) between states. 
%However, except for one-dimensional models, deriving an explicit formula for MFPTs is highly complex for higher-dimensional SCRNs. 
However, deriving explicit formulas for MFPTs can be highly complex.
 In order to address this problem, we first introduce the concept of $\textcolor{black}{coclique}\; level\; structure$ and develop theorems to determine whether certain SCRNs have this feature by studying associated graphs. Additionally, we develop an algorithm to identify, under specific assumptions, all possible \textcolor{black}{coclique} level structures associated with a given SCRN. Finally, we demonstrate how the presence of such a structure in a SCRN allows us to derive closed form formulas for both upper and lower bounds for the MFPTs.
 Our methods can be applied to SCRNs taking values in a generic finite state space and can also be applied to models with non-mass-action kinetics. 
We illustrate our results with examples from the biological areas of epigenetics, neurobiology and ecology.

}




\maketitle

\section{Introduction}
\label{sec:introduction}

\subsection{Overview}

%\textbf{Overview of SCRNS}

    Stochastic Chemical Reaction Networks (SCRNs) are a class of continuous time Markov chain models used to describe the stochastic dynamics of a chemical system undergoing a series of reactions that change the numbers of molecules of a finite set of species over time. 
    These models can be used to conduct theoretical studies in different biological areas such as chromatin regulation (see, for instance, 
    %Bruno et al. 
    \cite{bib:BWD2022}), enzymatic kinetics (see, for instance, %Kang et al.
    \cite{KangKhudaBukhshKoepplRempala}), and intracellular viral kinetics (see, for instance, 
    %Srivastava et al. 
    \cite{Srivastava} and 
    %Haseltine \& Rawlings 
    \cite{HaseltineRawlings}).

%\textbf{what we are trying to solve}

Formulas for mean first passage times (MFPTs) between states of the Markov chains are helpful in studying the stochastic behavior. However, while calculating an explicit formula is relatively straightforward for certain one-dimensional models, such as birth-death processes (see SI - Section \ref{sec:Appendix1Dmodel}), this typically becomes significantly more difficult for higher-dimensional SCRNs.
\textcolor{black}{A potential approach to compute MFPTs involves matrix inversion techniques \cite{epifinite}. While this method can, in principle, lead to closed-form expressions, it becomes algebraically complex when symbolic parameter dependencies are retained, even for moderate system sizes. As the system size grows, the dimension of the matrix to be inverted increases, leading to higher computational costs and even more complex expressions. In this case, fixing parameter values can help reduce the computational burden, but this also reduces information about how parameters influence the MFPT.} 

The inability to obtain an explicit analytical expression for the MFPT poses a challenge because having one allows for a clear understanding of how reaction rate parameters affect certain stochastic behavior of the system.
One possible approach to overcome and study the effect of parameter variations on system dynamics, without calculating explicit formulae, is to exploit comparison theorems for stochastic processes (see, for instance, 
%Campos et al. 
\cite{Monotonicitypaper}). However, this would require applying the theorem, if applicable, to a fixed parameter or combination of parameters of interest, thereby increasing the complexity of the study compared to having an analytical formula that can be directly analyzed.

   %\textbf{how do we solve it? Ind level structure}
 
    In this paper, we first introduce the concept of $\textcolor{black}{coclique}\; level\; structure$ and develop theorems to determine whether certain SCRNs have this feature by studying associated graphs. We also develop an algorithm to identify, under specific assumptions, all possible \textcolor{black}{coclique} level structures associated with such SCRNs. Then, we demonstrate how the presence of such a structure in a SCRN allows us to derive closed form formulas for both upper and lower bounds for MFPTs.
    We apply the theoretical tools developed in this paper to multiple examples to illustrate how they can be used to determine \textcolor{black}{coclique} level structures, to derive formulas for upper and lower bounds for the MFPTs, and then to understand how key biological parameters affect the stochastic behavior of the system. While our focus is  on SCRNs,  our definition of \textcolor{black}{coclique} level structure is for associated  continuous time Markov chains with  state spaces that are  finite subsets of the non-negative integer orthant and in which the set of all possible transition vectors is finite as well. Consequently, our approach potentially has broader applications to other models that have similar characteristics to SCRNs.

    % \textbf{paper structure}
    
    The structure of the paper is the following: we first introduce some background on SCRNs and mean first passage times (Section \ref{sec:BasicDefinitions}). We then introduce the notion of \textcolor{black}{coclique} level structure (Section \ref{sec:IndLevelStructureApp}) and describe the main theoretical tools developed in this paper (Section \ref{sec:MainResults}). Finally, we apply our results to multiple examples (Section \ref{sec:Applications}) and present some concluding remarks (Section \ref{sec:conclusion}).
    
    %, with proofs provided in Section \ref{sec:proofmainresults}. In Section \ref{sec:examples} we apply our theoretical tools to several examples, such as epigenetic regulation by chromatin modifications and enzymatic kinetics. Concluding remarks are presented in Section \ref{sec:conclusion}.
    %The Supplementary information (SI) file contains some further details and extensions of the main results and examples in the paper.
    
%Then, to make these theoretical tools more usable, we tailor the theorems to the structure of SCRNs.    
%%together with the special structure of CRNs
\subsection{Related work}

In this paper, we introduce the concept of \textcolor{black}{coclique} level structure for suitable SCRNs, which will be a partition of the finite state space of an associated continuous time Markov chain,
%. This partition takes inspiration from the concept of an independent set in graph theory \cite{EBertram1983}, that is a set including vertices not connected by an edge in the graph $\G$.
in which the sets in the partition are level sets of a linear function $L$ on the state space and they are \textcolor{black}{coclique}s in the sense of graph theory \textcolor{black}{\cite{vanMieghem2010,KelseyRoneyDougal2022}}, i.e., there are no direct transitions between states within a set in the partition.
To the best of our knowledge, this concept of \textcolor{black}{coclique} level structure has not been previously used to derive a reduced stochastic process. Here, we apply the function $L$ to a continuous time Markov chain, and analyze the result to estimate MFPTs for SCRNs.
Previous theoretical tools developed to evaluate upper and lower bounds for MFPTs are mostly suitable for computational studies \cite{Latouche} and often focus on specific models, such as imprecise birth-death chains \cite{SLopatatzidis2017} or population continuous time Markov chains \cite{MBackenkohler2020}. In contrast to these existing works, the theoretical tools that we develop enable the derivation of closed form formulas for MFPT bounds, making them suitable for analytical analysis. 

\subsection{Terminology and Notation}
\label{sec:PreliminariesAndNotation}

    Denote the set of integers by $\Z$. For an integer $d \geq 2$, we denote by $\Z^d$ the set of $d$-dimensional vectors with entries in $\Z$. Denote by $\Z_+ = \{0,1,2, \ldots \}$, the set of non-negative integers. For an integer $d \geq 2$, we denote by $\Z_+^d$ the set of $d$-dimensional vectors with entries in $\Z_+$. We denote by $\zero$, respectively, $\one$, a vector of any dimension where all entries are $0$'s, respectively, $1$'s. The size of $\zero$ or $\one$ will be understood from the context. The set of real numbers will be denoted by $\R$, $\R_+=[0,\infty)$, and $d$-dimensional Euclidean space will be denoted by $\R^d$ for $d \geq 2$. For integers $m,n \geq 1$, the set of $m \times n$ matrices with real-valued entries will be denoted by $\R^{m \times n}$. 
    %The set of complex numbers will be denoted by $\C$.
    For a matrix $A \in \R^{m \times n}$, we denote the kernel of $A$ by $\ker(A) := \{ u \in \R^{n}:\: Au=\zero \}$. Vectors are column vectors unless indicated otherwise and a superscript of $T$ will denote the transpose of a vector or matrix.

\section{Stochastic Chemical Reaction Networks (SCRNs)}
\label{sec:BasicDefinitions}

In this section, we provide basic definitions for a class of continuous time Markov chains called Stochastic Chemical Reaction Networks (SCRNs). The reader is referred to Anderson \& Kurtz \cite{AndersonKurtzBook} for a more in depth introduction to this subject. The models considered in our examples will be SCRNs and the state space of all of our models will be finite. 
%Then, we will introduce the definition of mean first passage time. 

    We assume there is a finite non-empty set $\Ss = \{\mathrm{S}_1,\ldots,\mathrm{S}_d\}$ of $d$ \textbf{species}, and a finite non-empty set $\Rs \subseteq \Z_+^d \times \Z_+^d$ that represents chemical \textbf{reactions}. We assume that $(w,w) \notin \Rs$ for every $w \in \Z^d_+$. The set $\Ss$ represents $d$ different molecular species in a system subject to reactions $\Rs$ which change the number of molecules of some species. For each $(v^{-},v^+) \in \Rs$, the $d$-dimensional vector $v^{-}$ (the \textbf{reactant vector}) counts how many molecules of each species are consumed in the reaction, while $v^{+}$ (the \textbf{product vector}) counts how many molecules of each species are produced. The associated reaction is usually written as
    \begin{equation}
	   \label{eq:ReactionNotation}
	   \sum_{i=1}^d (v^{-})_{i}\mathrm{S}_i \longrightarrow  \sum_{i=1}^d (v^{+})_{i}\mathrm{S}_i.
    \end{equation}
    To avoid the use of unnecessary species, we will assume that for each $1 \leq i \leq d$, there exists a vector $w=(w_1, \ldots,w_d)^T \in \Z_+^d$ with $w_i >0$ such that $(w,v)$ or $(v,w)$ is in $\Rs$ for some $v \in \Z^d_+$, i.e., each species is either a reactant or a product in some reaction.
    The net change in the quantity of molecules of each species due to a reaction $(v^{-},v^{+}) \in \Rs$ is described by $v^{+}-v^{-}$ and it is called the associated \textbf{reaction vector}. We denote the set of reaction vectors by $\V := \{ v \in \Z^d : v = v^{+}- v^{-} \text{ for some } (v^{-},v^{+}) \in \Rs\}$; we let $n := |\V|$, the size of $\V$; and we enumerate the members of $\V$ as $\{v_1,\ldots,v_n\}$. Note that $\V$ does not contain the zero vector because $\Rs$ has no elements of the form $(w,w)$. Different reactions might have the same reaction vector. For each $v_k \in \V$ we consider the set $\Rs_{v_k} := \{(v^{-},v^{+}) \in \Rs : v_k =v^{+}-v^{-} \}$. {\color{black} The matrix $S \in \R^{d \times n}$ whose columns are the elements in $\V$ will be called the \textbf{stoichiometric matrix}\footnote{\color{black} In this stoichiometric matrix, there are no duplicate reaction vectors. This parallels combining reaction intensity functions associated with the same reaction vector as in \eqref{defrate}.}.} In addition, we define a \textbf{conservation vector} $m$ (if there is one) as a $d$-dimensional non-zero vector such that $m^TS=0$ and we say that the conservation vector is $unique$ if $m$ is unique, up to  multiplication by a scalar.

%    Given $(\Ss,\Rs)$ we will consider an associated continuous time Markov chain $X=(X_1,$ $\ldots,X_d)^T$, with a state space $\X$ contained in $\Z^d_+$, which tracks the number of molecules of each species over time. Roughly speaking, the dynamics of $X$ will be given by the following: given a current state  $x=(x_1,\ldots,x_d)^T \in \X \subseteq \Z_+^{d}$, for each reaction $(v^{-},v^{+}) \in \Rs$, there is a clock which will ring at an exponentially distributed time (with rate $\Lambda_{(v^{-},v^{+})}(x)$). The clocks for distinct reactions are independent of one another. If the clock corresponding to $(v^{-},v^{+})\in \Rs$ rings first, the system moves from $x$ to $x+v^{+}- v^{-}$ at that time, and then the process repeats. We now define the continuous time Markov chain in more detail.

    Consider sets of species $\Ss$ and reactions $\Rs$, a non-empty set $\X \subseteq \Z^d_+$, and a collection of functions $\prp= \{\prp_{(v^{-},v^{+})}:\X \longrightarrow \R_+\}_{(v^{-},v^{+}) \in \Rs}$ such that for each $x \in \X$ and $(v^{-},v^{+}) \in \Rs$, if $x+v^{+}-v^{-} \notin \X$, then $\Lambda_{(v^{-},v^{+})}(x)=0$. Now, for $1 \leq k \leq n$, define
    \begin{equation}\label{defrate}
	   \rate_k(x) := \sum_{(v^{-},v^{+}) \in \Rs_{v_k}} \prp_{(v^{-},v^{+})}(x).
    \end{equation}
    Note that for each $x \in \X$ and $1 \leq k \leq n$, if $x +v_k \notin \X$, then $\rate_k(x) = 0$. The functions $\{\prp_{(v^{-},v^{+})}:\X \longrightarrow \R_+\}_{(v^{-},v^{+}) \in \Rs}$ are called \textbf{propensity} or \textbf{intensity} functions. A common form for the propensity functions is the following, which is associated with \textbf{mass action kinetics}:
    \begin{equation}
	   \label{eq:PropensityFunctions}
	   \prp_{(v^{-},v^{+})}(x) = \kappa_{(v^{-},v^{+})}\prod_{i=1}^{d}(x_i)_{(v^{-})_i},
        \end{equation}
    where $\{\kappa_{(v^{-},v^{+})}\}_{(v^{-},v^{+}) \in \Rs}$ are non-negative constants and for $m,\l \in \Z_+$, the quantity $(m)_\l$ is the falling factorial, i.e., $(m)_0 := 1$ and $(m)_\l := m(m-1)\ldots(m-\l+1)$.
    
    A \textbf{stochastic chemical reaction network (SCRN)} (associated with $(\Ss,\Rs,\X,\prp)$) is a continuous time Markov chain $X$ with state space $\X$ and infinitesimal generator $Q$ given for $x,y \in \X$ by
    \begin{equation}
	   \label{eq:TransitionMatrixQ}
	   Q_{x,y} = \begin{cases}
		\rate_k(x) & \text{ if } y-x = v_k \text{ for some } 1 \leq k \leq n, \\
		- \sum_{k=1}^n\rate_k(x) & \text{ if } y = x, \\
		0 & \text{ otherwise.}
	\end{cases}
    \end{equation}
    
    If a SCRN associated with $(\Ss,\Rs,\X,\prp)$ has a conservation vector $m \ne 0$ and $m^TX(0)=x_{\mathrm{tot}}$ for some integer $x_{\mathrm{tot}}\ge0$, then $m^TX(t)=x_{\mathrm{tot}}$ for every $t \geq 0$. Consequently, we can reduce the dimension of the continuous time Markov chain describing the system by one. 
    In this paper, %we will be assuming that 
    we will initially be considering SCRNs for which $m=(1,\ldots,1)^T$ is a conservation vector. Then, the projected process $\widecheck{X}=(X_1,\ldots,X_{d-1})^T$ is again a continuous time Markov chain with finite state space $\widecheck{\X}=\{(x_1,\ldots,x_{d-1})^T \in \Z^{d-1}_+ : (x_1,\ldots,x_{d-1},x_{\mathrm{tot}} - \sum_{i=1}^{d-1} x_i)^T \in \X \}$.
    We denote its infinitesimal generator by $\widecheck{Q}$. 
    We will assume that $|\widecheck{\X}| > 1$. Let $\B$ be a nonempty subset of $\widecheck{\X}$ such that $\B \neq \widecheck{\X}$, and let 
    \begin{equation}\label{MFPTNorrisformula}
    \tau_\B := \inf \{t \geq 0: \widecheck{X}(t) \in \B\}.
    \end{equation}
    Then, the \textbf{mean first passage time (MFPT)} (for $\widecheck{X}$) from $x \in \widecheck{\X}$ to $\B$ can be defined as 
    \begin{equation}\label{MFPTxB}
    h_{x,\B}=\E[\tau_\B \;|\; \widecheck{X}(0) = x].
    \end{equation}
    If $\B = \{y\}$ for some $y \in \widecheck{\X}$, we use the notation $h_{x,y} := h_{x,\{y\}}$. 
    {\color{black} After considering this simple projection first, later on in Section \ref{generalization}, we will consider the situation where there are $\wccNum>1$ linearly independent conservation vectors for a SCRN. Then, we can reduce the dimension of the Markov chain $X$ by $\wccNum$, in which case we also denote the projected process by $\widecheck{X}$, its state space by $\widecheck{\X} \subset \Z^{d-\wccNum}_+$ and its infinitesimal generator by $\widecheck{Q}$. The mean first passage time for $\widecheck{X}$ is then defined by \eqref{MFPTNorrisformula}--\eqref{MFPTxB}.}

%%%%%%%%%%%%%%%%%%%%%%%%%%%%%%%%%%%%%%%%%%%%%%%%%%%%%

\section{\textcolor{black}{Coclique} level structure for SCRNs}
\label{sec:IndLevelStructureApp}

In this section, we define the notion of $\textcolor{black}{coclique}\; level\; structure$ for SCRNs %continuous time Markov chains, 
and introduce some useful assumptions, definitions, and lemmas that we use in this paper.

First, consider a SCRN, as defined in Section \ref{sec:BasicDefinitions}, and  assume the following:
     \begin{assumption}
    \label{assumption:Unimolecular_change}
    %For each $(v^{-},v^+) \in \Rs$, the entries of the reaction vector $v= v^+-v^-$ consist of $d-2$ zeros, a single one and a single minus one, in any order.
    Assume that $d\ge 2$. For each $v \in \V = \{v_1,\dots,v_n\}$, the entries of the reaction vector $v$ consist of $d-2$ zeros, a single one and a single minus one, in any order.
    \end{assumption}
    When satisfied, this assumption implies that the net change in the number of species stemming from a reaction consists of consuming one molecule of a given species to produce a molecule of a different species. For example, when $d=3$,  reactions satisfying this assumption include  $\mathrm{S}_1 \longrightarrow \mathrm{S}_2$ and  $\mathrm{S}_1 + \mathrm{S}_3 \longrightarrow \mathrm{S}_2 + \mathrm{S}_3$, both of which have associated reaction vector $(-1,1,0)^T$.

    %In terms of \textbf{SCRN graph structure}, 
    We consider a directed graph $\G$ associated with a SCRN satisfying Assumption \ref{assumption:Unimolecular_change}, in which the vertices represent the species of the system \textcolor{black}{(we will label these vertices $1,\dots,d$ or $\mathrm{S}_1,\dots,\mathrm{S}_d$)} and the directed edges between vertices represent the reaction vectors $\{v_1,\dots,v_n\}$. In this graph, an edge $\edge_k=(i,j)$, with $k\in \{1,2,\dots,n\}$ and $i \ne j$, exists if and only if there is a reaction vector $v_k\in \V$ where $v_k(i) =-1$, $v_k(j)=1$, $v_k(\ell) =0$ for $\ell \notin \{i  ,j\}$. Note that, if we have two edges $\edge_{k}=(i,j)$ and $\edge_{k'}=(j,i)$, with $i \ne j$, then $v_{k}=-v_{k'}$.
%    {\color{black} it the following some kind of definition or a fact from the previous??} ADD EXAMPLE OF ALLOWED REACTIONS AND REMOVE DEFINITION OF REVERSIBLE Furthermore, if a reaction between two species $\mathrm{S}_i,\mathrm{S}_j\in \Ss$ is reversible and the edge $\edge_a$ going from the vertex $i$ to the vertex $j$ (with the associated reaction vector defined as $v_a$) corresponds to the forward reaction, then the backward reaction is associated with a directed edge $\edge_b$ going from the vertex $j$ to the vertex $i$ and the reaction vector $v_b$ associated with the edge $\edge_b$ is $v_b=-v_a$.

Before introducing the definition of \textcolor{black}{coclique} level structure, we provide some useful definitions and lemmas that we will use in this section. Given the directed graph $\G$, a \textbf{\textcolor{black}{coclique}} is a non-empty set of vertices such that any two vertices in the set do not have a direct edge between them. The graph $\G$ is \textbf{bipartite} if the vertex set can be partitioned into two non-empty \textcolor{black}{coclique}s. 
%{\color{black} the following concepts are only used in Lemma 3.1?? move to SI} The \textbf{degree} of a vertex is the number of edges that are incident to the vertex, where we count both incoming and outgoing edges. A \textbf{weakly directed tree} is a directed graph whose underlying undirected graph is a tree. Given a directed graph with $d$ vertices, a \textbf{weakly directed spanning tree} is a subgraph of the graph with all $d$ vertices and such that it is a weakly directed tree, and it will have $d-1$ edges. We abbreviate weakly directed spanning tree as wd-spanning tree. {\color{black} end here.} 
The directed graph $\G$ is \textbf{weakly connected} if the underlying undirected graph is connected. Note that this implies $n\ge d-1$.
%{\color{black} start here again} A weakly connected graph $\G$ has a wd-spanning tree. {\color{black} end here again} 
%
%
%For the rest of this paper we will focus on SCRNs whose associated graphs $\G$ are weakly connected. This is without loss of generality, since if $\G$ has more than one weakly connected component, the following lemma and theorems can be applied to each of the weakly connected components separately \textcolor{black}{(see SI - Section \ref{SI:disconnected-comp-example} for an illustrative example showing how weakly connected components can be analyzed separately, as they represent non-interacting subsystems)}.
%
{\color{black} In this section,} we will focus on SCRNs whose associated graphs $\G$ are weakly connected{\color{black}, and we will prove results for such SCRNs in Sections \ref{suffcond}--\ref{indeplevelstructure}. }

{\color{black} 
For a SCRN whose associated graph $\G$ is not weakly connected, we can decompose $\G$ into finitely many disjoint \textbf{weakly connected components}. Each of the weakly connected components of $\G$ is a weakly connected graph and there are no edges between different weakly connected components. 
In Section \ref{generalization}, we will show how to leverage the results in Sections \ref{suffcond}--\ref{indeplevelstructure} for the single weakly connected component case to treat the general case where $\G$ is a finite disjoint union of weakly connected components.
}

%\begin{definition}\label{conslawdefSI}
%Given a graph with $d$ species and $n$ reaction vectors, we define the \textbf{stoichiometric matrix} $S$ as the matrix $S=[v_1,v_2,\dots,v_{p}]\in \Z^{d \times n}$. In addition, we define a \textbf{conservation vector} $m$ associated with $S$ as a $d$-dimensional non-zero vector such that $m^TS=0$. This means that $m^Tx$ is constant for all the possible choices of $x=(x_1,x_2,\dots,x_d)^T$. Furthermore, we say that the conservation vector is $unique$ if $m$ is unique up to a scalar.
%\end{definition}

%\subsubsection{Proofs of Lemma \ref{lemma2}, Theorems \ref{thm2a} and \ref{thm2b}}
%\label{proofslemmastheorems}
%Let us start providing the proof of Lemma \ref{lemma2}. To this end, let us first %introduce a Lemma needed in the proof of Lemma \ref{lemma2}.

\begin{lemma}\label{lemmabis}
	Consider a SCRN satisfying Assumption \ref{assumption:Unimolecular_change} with $d \ge 2$ species, $\Ss = \{\mathrm{S}_1,$ $\ldots,\mathrm{S}_d\}$, and $n$ reaction vectors, $\V =\{v_1,\dots,v_n\}$. Assume that the associated graph $\G$ is weakly connected. Then, the rank of the associated stoichiometric matrix $S$, whose columns are given by the elements of $\V$, is equal to $d-1$.
Furthermore, $(1,\ldots,1)^T$ is the unique (up to scalar multiplication) conservation vector associated with the stoichiometric matrix $S$.
\end{lemma}
%\begin{proof}
%From the proof of Lemma \ref{lemmabis}, we have that $\one^T S=0$ and $S^T \one=0$. Thus, $\one$ is a conservation vector and, since rank[$S^T$] = $d-1$ (see Lemma \ref{lemmabis}), the nullspace of $S^T$ is one-dimensional by the rank-nullity theorem \cite{LinearAlgebra}. Hence, $\one$ is the only vector such that $S^T\one=0$ up to scalar multiplication. We conclude that $\one$ is the unique conservation vector of the SCRN considered (up to scalar multiplication).
%\end{proof}

The proof of Lemma \ref{lemmabis} is given in SI - Section \ref{SI:proofs}.

{\color{black}
\begin{remark}
\label{rmk:deficiencyzero}
    If each vertex in $\G$ is regarded as a species and each edge in $\G$ is regarded as a reaction, then $\G$ can be interpreted as a chemical reaction network in the sense of Feinberg \cite{Feinberg} (where the complexes in this network are simply the species). For this, one can determine its deficiency $\delta=d-\ell-s$, where $d$ is the number of complexes (which is the number of species in this case), $\ell$ is the number of weakly connected components in $\G$ and $s$ is the rank of the stoichiometric matrix $S$ associated with $\G$. One can see by Lemma \ref{lemmabis} that, when $\G$ is weakly connected, the chemical reaction network associated with $\G$ has deficiency zero\footnote{{\color{black} If one starts with the chemical reactions associated with a SCRN, for example as in \eqref{reacs2D}, one can form the associated complex-reaction graph, also known as the Feinberg-Horn-Jackson graph, and compute the deficiency associated with that graph. This graph need not be the same as our graph $\G$, and its deficiency need not necessarily be zero, for example as in \eqref{reacs2D}.}}. 
\end{remark}
}

%
%%%Let us now provide the proof of Theorems \ref{thm2a} and \ref{thm2b}. To this end, let us first introduce Lemma \ref{lemmaweaklyconnect}, that will be needed to prove Theorem \ref{thm2b}.
% 
%\begin{remark}
%\label{remark1}
	%Lemma \ref{lemmabis} implies that SCRNs, satisfying Assumption \ref{assumption:Unimolecular_change} and for which the graph $\G$ is weakly connected, have a unique conservation vector $m = (1,\ldots,1)^T$ (up to a scalar multiplication) and then, as pointed out at the end of Section \ref{sec:BasicDefinitions},
 Under the assumptions of Lemma \ref{lemmabis}, 
 %if $m^TX(0) = x_{tot}
 %\ge 0$, 
 we can introduce a \textbf{projected continuous time Markov chain} $\widecheck{X}=\{\widecheck{X}(t):\: t \ge 0\}$ in which the state $\widecheck x$ represents the number of molecules of the first $d-1$ species, that is, $\widecheck x=(x_1,\dots,x_{d-1})^T$. Note that the choice to express $x_d$ as a function of the other $x_i$ with $i\in\{1,\dots,d-1\}$ is without loss of generality, since the labeling of the species can always be reordered so that the species chosen to be expressed as a function of the others is the last one. 
%\end{remark}
%
%Assuming that we have the assumptions of Lemma \ref{lemmabis}
% %with $d$ species, $\Ss = \{\mathrm{S}_1,\ldots,\mathrm{S}_d\}$ and $p$ reactions,
 %we introduce the projected Markov chain $\widecheck{X}=(X_1,\dots,X_{d-1})$.
 The process $\widecheck{X}$ is a continuous time Markov chain defined on the finite state space 
 \begin{equation}
     \begin{aligned}
      \widecheck{\mathcal{X}}&\vcentcolon=\left\{\widecheck x=(x_1,\dots,x_{d-1})^T\in \Z_{+}^{d-1}:\left(x_1,\dots,x_{d-1},x_{\mathrm{tot}}-\sum_{i=1}^{d-1}x_i\right)^T \in \X\right\}\\
      &\;\subset\left\{\widecheck x=(x_1,\dots,x_{d-1})^T\in \Z_{+}^{d-1}:x_1+\dots+x_{d-1}\le x_{\mathrm{tot}}\right\},   
     \end{aligned}
 \end{equation}
where $x_{\mathrm{tot}}=\sum_{i=1}^{d}X_i(0)$. We will assume that $|\widecheck{\X}| > 1$, and the infinitesimal generator of $\widecheck{X}$ will be denoted by $\widecheck Q$.

We now introduce the definitions of \textcolor{black}{coclique} level function and \textcolor{black}{coclique} level structure. A \textbf{\textcolor{black}{coclique} level function} for $\widecheck{X}$ is a linear function $L: \Z^{d-1} \to \Z$ such that for each $k= 1,\dots,n$,
 %
% one of the following is true:
 %
%		\begin{itemize}
%			\item for every $x,y\in \widecheck{\mathcal{X}}$ such that $\widecheck Q_{x,y}>0$ and $y-x=\widecheck{v}_k$, we have $L(y)=L(x)+1$;
%			\item for every $x,y\in \widecheck{\mathcal{X}}$ such that $\widecheck Q_{x,y}>0$ and $y-x=\widecheck{v}_k$, we have $L(y)=L(x)-1$.
%		\end{itemize}
	%
 \begin{equation}\label{indlevfunctiondef}
L(\widecheck{v}_k) \in \{-1,+1\},
 \end{equation}
where $\widecheck{v}_k\in \Z^{d-1}$ is the vector obtained from $v_k$ by removing the last element.
 If such an $L$ exists, it can be written as 
 \begin{equation}\label{ILFformula}
  L(x) = b^T x\;\; \mathrm{for}\; x \in \Z^{d-1}\; \mathrm{and\; some}\; b \in \Z^{d-1},   
 \end{equation}
where, upon partitioning the set of edges of the associated graph $\G$ into two disjoint subsets $\Eu=\{\edge_k:L(\widecheck{v}_k)=1\}$ and $\Ed=\{\edge_k:L(\widecheck{v}_k)=-1\}$ (where one of these may be empty), the vector $b=(b_1,\dots,b_{d-1})^T$ solves the system of equations
  %$\{b_i,i=1,\dots,d-1\}\subset \R$:
	\begin{equation}\label{syst}
	\sum_{i=1}^{d-1}b_{i} \widecheck{v}_k(i)=\begin{cases}
	+1\;\;\;\text{if }\; \edge_k \in \Eu\\
	-1\;\;\;\text{if }\; \edge_k \in \Ed\\
	\end{cases}\;\;\;\;\; \text{for } k=1,\dots,n.
	\end{equation}
 Finally, for \textcolor{black}{a coclique} level function $L$, the (ordered) partition $\{\L_{\ell}, \dots, \L_u \}$, with \begin{equation}\label{indlevstructure}
	    \L_z:=\{x \in \widecheck{\mathcal{X}}:\: L(x)=z\}\; \mathrm{for}\; z=\ell,\ell+1,\dots,u-1,u,
	\end{equation}
 is called a \textbf{\textcolor{black}{coclique} level structure} for $\widecheck{X}$, with 
 \begin{equation}\label{lowhighletters}
 \ell=\min \{L(x):x \in \widecheck{\X}\},\quad \quad u=\max \{L(x):x \in \widecheck{\X}\}.
 \end{equation}
The sets $\L_{\ell}, \dots, \L_u$ are {\color{black} cocliques} in the {\color{black} Markov chain graph for $\widecheck{X}$ which consists of states of $\widecheck{\X}$ with edges given by $\{(\widecheck{x}, \widecheck{x}+\widecheck{v}_k):k=1,\dots,n,\;\; \widecheck{x}\in \widecheck{\X},\;\; \widecheck{x}+\widecheck{v}_k\in \widecheck{\X}\}$}.
%\label{indleveldef}
%

%\noindent \textbf{Proof of Theorem \ref{thm2a}}
%
\begin{lemma}\label{thm2a}
%\textcolor{black}{remove first part and put this above (after equation 3.3 (To add: note that \{Eu,El\}, in each one can be empty, forms a partition of graph  G). Leave the second part in lemma 3.2}
%	Let us consider a SCRN satisfying Assumption \ref{assumption:Unimolecular_change} and having a single linkage class. Partition the set of edges of the associated $\G$ graph into two subsets $\Eu$ and $\Ed$. Furthermore, let us define $\widehat{\Ss}$ as the set of d-1 species obtained from $\Ss$ by removing the last species $\mathrm{S}_d$ ($\widehat{\Ss}={\Ss} \backslash \mathrm{S}_d$). Now, introducing the variables $b_{i}$ with $\mathrm{S}_i\in \widehat{\Ss}$, let us define the following system of equations:	 
%	\begin{equation}\label{syst}
%	\sum_{\mathrm{S}_i\in \widehat{\Ss}}b_{i}v_k(i)=\begin{cases}
%	+1\;\;\;if\; e_k \in \Eu\\
%	-1\;\;\;if\; e_k \in \Ed\\
%	\end{cases}\;\;\;\;\;
%	\end{equation}
%	Then, the system (\ref{syst}) admits either one or zero solutions.
Consider a SCRN satisfying the assumptions in Lemma \ref{lemmabis}. 
%If the projected Markov chain $\widecheck{X}$ admits an \textcolor{black}{coclique} level function $L$ as in \eqref{ILFformula},
%%, that can be written as $L(x_1,\dots,x_{d-1})=\sum_{i=1}^{d-1}b_{i}x_i$ for $(x_1,\dots,x_{d-1}) \in \Z^{d-1}$ and with $b\in \Z^{d-1}$,
%then letting $\Eu=\{\edge_k:L(\widecheck{v}_k)=1\}$ and $\Ed=\{\edge_k:L(\widecheck{v}_k)=-1\}$, the vector $b=(b_1,\dots,b_{d-1})^T$ solves the system of equations
  %%$\{b_i,i=1,\dots,d-1\}\subset \R$:
%	\begin{equation}\label{syst}
%	\sum_{i=1}^{d-1}b_{i}v_k(i)=\begin{cases}
%	+1\;\;\;if\; e_k \in \Eu\\
%	-1\;\;\;if\; e_k \in \Ed\\
%	\end{cases}\;\;\;\;\; \text{for } k=1,\dots,n.
%	\end{equation}
%%
%%admits a solution equal to $b\in \Z^{d-1}$.
If we partition the set of the edges of the associated graph $\G$ into two disjoint subsets $\Eu$ and $\Ed$, where one of these may be empty, then, for this partition $\{\Eu,\Ed\}$, the system \eqref{syst} %with $b=(b_1,\dots,b_{d-1})^T \in \R^{d-1}$ 
admits either zero or one solution, and the solution, if one exists, has integer entries.
\end{lemma}

\begin{proof}
%
%{\color{black} START VERSION ONE} An \textcolor{black}{coclique} level function for $\widecheck{X}$ is a linear function $L: \Z^{d-1} \to \Z$ such that \eqref{indlevfunctiondef} holds. Such an $L$ is of the form $L(x)=b\cdot x$, $x\in \Z^{d-1}$ for some $b \in \Z^{d-1}$ and satisfies \eqref{syst}  where 
%$\Eu=\{\edge_k:L(\widecheck{v}_k)=b^T \widecheck{v}_k=+1\}$ and $\Ed=\{\edge_k:L(\widecheck{v}_k)=b^T \widecheck{v}_k=-1\}$. Thus, the first part of the theorem follows. {\color{black} END VERSION ONE}
%{\color{black} START VERSION TWO}
%The first part of the theorem follows from \eqref{indlevfunctiondef}. 
%{\color{black} END VERSION TWO}

%Define $\widehat{\Ss}$ to be the set containing $d-1$ species obtained from $\Ss$ by removing the last species $\mathrm{S}_d$ ($\widehat{\Ss}=\Ss \setminus \mathrm{S}_d$). Then, 
%We rewrite the system (\ref{syst}) in matrix-vector form. In particular, 
%For the second part, 
Let $\widecheck{S}\in \Z^{(d-1) \times n}$ be the first $(d-1)$ rows of the stoichiometric matrix $S$ associated with the SCRN. %with the last row removed.
%(that is, the row corresponding to the species removed from $\Ss$ to get the set $\widehat{\Ss}$). 
Then, the system (\ref{syst}) can be re-written in matrix-vector form as
\begin{equation}\label{systMATRIX}
\widecheck{S}^T b = w,
\end{equation}
where $w_k$ is $+1$ or $-1$ depending on whether $\edge_k \in \Eu$ or $\edge_k \in \Ed$, respectively, where $k=1,\dots,n$. %{\color{black} BEGIN OLD} Now, by Lemma \ref{lemmabis}, we know that rank[$S^T$] $=d-1$. Furthermore, since $Null(S^T)=span\{\one\}$ (by Assumption \ref{assumption:Unimolecular_change}) {\color{black} END OLD} {\color{black} BEGIN NEW} From Assumption \ref{assumption:Unimolecular_change}, we have that $\ker (S^T) \supseteq \text{span} \{\one\}$. Since  $\rank(S^T)$ $=d-1$ (see Lemma \ref{lemmabis}), we have that $\dim \ker (S^T) = 1$ and thus $\ker (S^T) = \text{span} \{\one\}$. {\color{black} this is actually the proof for Lemma \ref{lemmabis}} {\color{black} END NEW}, 
By Lemma \ref{lemmabis}, we know that $\ker (S^T) = \text{span} \{\one\}$, and thus the last row of $S$ is a linear combination of the first $(d-1)$ rows of $S$. This means that removing the last row of $S$ does not affect its rank, and thus $\rank\left(\widecheck{S}^T\right) = d-1$. 
%$Null(\widecheck{S}^T)$ is the zero vector in $\R^{d-1}$. It follows that $\rank\left(\widecheck{S}^T\right) = d-1$. 
Hence, if $w$ is in the \text{range} of $\widecheck{S}^T$ the system \eqref{syst} admits a unique solution, while if it is not in the \text{range}, the system does not admit a solution. 

Suppose that \eqref{syst} admits a solution $b \in \R^{d-1}$. Since $\G$ is weakly connected, there is an edge connecting the vertex $d$ with another vertex in $\G$, say vertex $i$. This means that there exists $v_k \in \V$ such that %$|v_k(i)|=|v_k(d)|=1$
$v_k(i)=-v_k(d) \in \{-1,1\}$ and $v_k(\ell)=0$ for $\ell \in \{1,2,\dots,d-1\} \setminus \{i\}$. Then, since $b \in \R^{d-1}$ solves the $k^{th}$ equation of \eqref{syst}, we have that $\abs{b_i}=1$, which means $b_i \in \Z$. Using similar logic, for a vertex %$j \neq d$
$j \notin \{i,d\}$ in $\G$ such that there is an edge between $j$ and $\{i,d\}$, we have that either $\abs{b_i-b_j}=1$ or $\abs{b_j}=1$, and thus $b_j \in \Z$. Since $\G$ is weakly-connected and has finitely many vertices, we can iteratively show that all of the entries of $b$ have integer values.
\end{proof}

{\color{black}
We use the following example to illustrate our theory in a simple context before giving more complex examples in Section \ref{sec:Applications}.
\begin{figure}[t]
            \centering
            \includegraphics[scale=0.39]{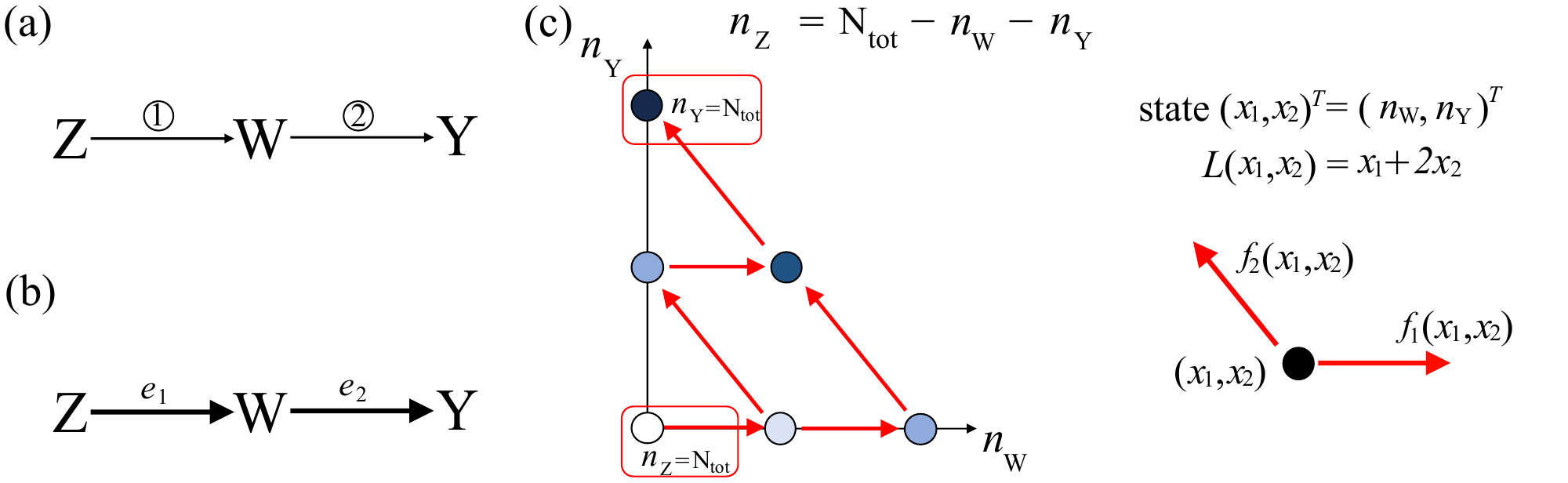}
            \caption{\small {\color{black} { \bf Three-species cascade motif: reaction diagram, graph $\G$ and associated Markov chain (Example \ref{ex:illustrativeexample}).}
(a) Chemical reaction system diagram.  The numbers on the arrows correspond to the reactions associated with the arrows as described in (\ref{reacs-3speciesexample}) in the main text. 
(b) Graph $\G$ associated with the chemical reaction system in panel (a).
(c) State space and transitions for the projected continuous time Markov chain $\widecheck{X}=\{
(X_1(t),X_2(t))^T:\: t \ge 0\}$, which keeps track of $(n_{\mathrm{W}},n_{\mathrm{Y}})$ through time. Here, we consider $\mathrm{N_{tot}}=2$ and we use dots to represent the states, and red double-ended arrows to represent transitions in both directions. 
Additionally, we use shades of blue to distinguish the level to which each state belongs. The function $L(x_1,x_2)$ associated with the coclique level structure is $L(x_1,x_2)=x_1+2x_2$.
%
%(c)
%
The rates associated with the one-step transitions for the projected Markov chain $\widecheck{X}$ are given in (\ref{rates3species}).}
            }
   \label{fig:illexample}
        \end{figure}
    \begin{example}
    \label{ex:illustrativeexample}
    Consider the following SCRN with mass action kinetics, involving three species and two irreversible reactions:
\begin{equation}\label{reacs-3speciesexample}
\begin{aligned}
&{\large \textcircled{\small 1}}\;\ce{Z ->[$\alpha$] W}, \quad
{\large \textcircled{\small 2}}\;\ce{W ->[$\beta$] Y},
\end{aligned}
\end{equation}
in which $\alpha, \beta > 0$. The diagram of this SCRN is shown in Figure \ref{fig:illexample}(a).  Let the species vector be $ x = (n_{\mathrm{W}}, n_{\mathrm{Y}}, n_{\mathrm{Z}})$, where 
$n_{\mathrm{W}}$, $n_{\mathrm{Y}}$, $n_{\mathrm{Z}}$ denote the number of molecules of species W, Y, Z, respectively.
The total number of molecules in this system is conserved, that is, $n_{\mathrm{Z}} + n_{\mathrm{W}} + n_{\mathrm{Y}} = \Ntot$. The reaction vectors associated with this SCRN are $v_1 = (1, 0, -1)^T$ and  $v_2 = (-1, 1, 0)^T$. Then, Assumption \ref{assumption:Unimolecular_change} is satisfied, and the associated graph $\G$ is shown in Figure \ref{fig:illexample}(b),
% %
% \begin{equation}
% \label{3species_graph}
% \mathrm{Z} \xrightarrow{e_1} \mathrm{W} \xrightarrow{e_2} \mathrm{Y},
% \end{equation}
% %
which is weakly connected.
%and that $\G$ is bipartite. 
%
By Lemma \ref{lemmabis}, our SCRN has a unique conservation vector $m = (1,1,1)^T$, and then we can introduce a projected continuous time Markov chain $\widecheck{X}=\{
(X_1(t),X_2(t))^T:\: t \ge 0\}$, which keeps track of $(n_{\mathrm{W}},n_{\mathrm{Y}})$ through time.
Since the total number of molecules $\Ntot$ is conserved, the state space is $\widecheck{\X}= \{\widecheck{x}=(x_1,x_2)^T \in \Z_+^2 :\: x_1 + x_2 \leq \Ntot \}$. 
The potential one-step transitions for $\widecheck X$ from $ x \in \widecheck{\X}$ are illustrated in Figure \ref{fig:illexample}(c) (for $\Ntot=2$), where the associated transition vectors are 
\begin{equation}\label{vectors3species}
\widecheck{v}_1=(1,0)^T,\;\;\; \widecheck{v}_2=(-1,1)^T,
\end{equation}
and the infinitesimal transition rates are
   \begin{equation}\label{rates3species}
    \begin{aligned}
&\widecheck{Q}_{\widecheck{x},\widecheck{x}+\widecheck{v}_1}=f_1(\widecheck{x}) = \alpha(\Ntot -(x_1+x_2)),\;\;\; \widecheck{Q}_{\widecheck{x},\widecheck{x}+\widecheck{v}_2}=f_2(\widecheck{x}) = \beta x_1.
    \end{aligned}
   \end{equation}

We now determine a \textcolor{black}{coclique} level function and the associated \textcolor{black}{coclique} level structure for the projected Markov chain $\widecheck{X}$. To this end, consider the partition of the edge set given by $\Eu=\{\edge_1,\edge_2\},\; \Ed=\emptyset$.
According to Lemma \ref{thm2a}, since $\G$ is weakly connected, and Assumption~\ref{assumption:Unimolecular_change} is satisfied, the system of equations \eqref{syst} admits either zero or one solution. In this case, solving the system of equations
$$ b_1 =1,\; -b_1 + b_2 = 1$$
yields the unique solution $b =(b_1,b_2)^T= (1,2)^T$. Therefore, the function $L(x_1,x_2) = x_1 + 2x_2$ is a \textcolor{black}{coclique} level function for $\widecheck{X}$.
The corresponding \textcolor{black}{coclique} level structure is the (ordered) partition $\{\L_{\ell},\dots,\L_u\}$ of the state space $\widecheck{\X}$, where each level set is defined as
$\L_z := \left\{ x \in \widecheck{\X} : L(x) = z \right\}, z = \ell, \dots, u$, with
$\ell = \min\{L(x) : x \in \widecheck{\X}\}=0$ and $u = \max\{L(x) : x \in \widecheck{\X}\}=2\Ntot$.

It is important to note that this is only one possible \textcolor{black}{coclique} level structure for $\widecheck{X}$. In the next section (Section \ref{sec:MainResults}), we develop theoretical tools to identify all possible \textcolor{black}{coclique} level structures for SCRNs satisfying Assumption \ref{assumption:Unimolecular_change}.
\end{example}
}

\section{Main Results}
\label{sec:MainResults}

    %{\color{teal} BEGIN REMOVE} In this section, we first describe the main theoretical tools developed in the paper (see Section \ref{suffcond}). {\color{teal} END REMOVE}
    
    {\color{black} In Sections \ref{suffcond}--\ref{indeplevelstructure}, we develop theoretical tools for SCRNs satisfying Assumption \ref{assumption:Unimolecular_change} and whose associated graph $\G$ is weakly connected, and in Section \ref{generalization}, we relax the latter condition and consider the case where $\G$ does not need to be weakly connected. }
    %\textcolor{black}{SB COMMENT: I think this new version is a little repetitive. You say upfront that section 4.1-4.3 are for weakly connected case, and 4.4 is the extension. Then, there is "More precisely", and everything said before is repeated here. In my opinion, I would live the old version removing "(see Section 4.1)".}
    More precisely, {\color{black} in Section \ref{suffcond},} we introduce a theorem to determine all of the possible \textcolor{black}{coclique} level functions associated with SCRNs satisfying Assumption \ref{assumption:Unimolecular_change} and with $\G$ weakly connected. We then derive theoretical tools to determine when a projected continuous time Markov chain associated with a SCRN admits a \textcolor{black}{coclique} level function by studying {\color{black} the structure of the graph $\G$ associated with the SCRN}.
    Additionally, we develop an algorithm to find, under certain assumptions, all of the possible \textcolor{black}{coclique} level functions associated with these SCRNs (see Section \ref{algIndLevFac}).
    We derive analytical formulas for upper and lower bounds on MFPTs for SCRNs having this type of structure (see Section \ref{indeplevelstructure}). 
    {\color{black} In Section \ref{generalization}, we generalize these results to study SCRNs that satisfy Assumption \ref{assumption:Unimolecular_change} and whose associated graph $\G$ does not need to be weakly connected.}
    We apply our results to several examples in Section \ref{sec:Applications}.

\subsection{Existence and characterization of \textcolor{black}{coclique} level structures}
\label{suffcond}

In the first theorem in this section, we show that if \eqref{syst} admits a solution, % and its coordinates are all integers {\color{black} check and add to theorem 3.1: if there is a solution then it has integer entries}, 
then there is \textcolor{black}{a coclique} level function for $\widecheck{X}$. Moreover, all \textcolor{black}{coclique} level functions for  $\widecheck{X}$ can be obtained in this way.

\begin{theorem}\label{thm2b}
    %Consider a SCRN satisfying the assumptions in Lemma \ref{thm2a}. 
    Consider a SCRN satisfying Assumption \ref{assumption:Unimolecular_change}.
    %with $d \ge 2$ species. %, $\Ss = \{\mathrm{S}_1,$ $\ldots,\mathrm{S}_d\}$, and $n$ reaction vectors, $\V =\{v_1,\dots,v_n\}$. 
    Assume that the associated graph $\G$ is weakly connected. 
 %    For each partition $\{ \Eu, \Ed \}$ of the edge set of $\G$ (where one of $\Eu$ and $\Ed$ may be empty), the system \eqref{syst} either has no solution or admits a unique solution $b$ in $ \Z^{d-1}$. In the latter case, the vector $b$ gives an \textcolor{black}{coclique} level function 
 %    \begin{equation*}
 %        \quad \quad \quad \quad \quad \quad L(x)=\sum_{i=1}^{d-1}b_{i}x_i,  \quad \quad \quad \text{ for each } x=(x_1,\dots,x_{d-1}) \in \Z^{d-1}.
 % \end{equation*}
% \textcolor{black}{If the system \eqref{syst} admits a solution $b\in \Z^{d-1}$, then the projected Markov chain $\widecheck{X}$ has an \textcolor{black}{coclique} level structure associated with the \textcolor{black}{coclique} level function $L$ as in \eqref{ILFformula}.}
%
The set of all \textcolor{black}{coclique} level functions for $\widecheck{X}$ is the set of all functions of the form \eqref{ILFformula},
% all of the \textcolor{black}{coclique} level functions of the projected Markov chain $\widecheck{X}$ have the form of
    % \begin{equation*}
    %     \quad \quad \quad \quad \quad \quad L(x)=\sum_{i=1}^{d-1}b_{i}x_i,  \quad \quad \quad \text{ for each } x \in \Z^{d-1},
    % \end{equation*}
    where %\textcolor{black}{ 
    $b = (b_1,\dots,b_{d-1})^T \in \Z^{d-1}$ is a solution of the system \eqref{syst} %\textcolor{black}{\st{
    for some partition $\{ \Eu, \Ed \}$ of the edge set of $\G$ (and one of $\Eu$ and $\Ed$ may be empty).
    %} vectors $b = (b_1, \dots, b_{d-1})^T \in \Z^{d-1}$ are obtained as solutions of the system \eqref{syst} for all partitions ${ \Eu, \Ed }$ of the edge set of $\G$ (where either $\Eu$ or $\Ed$ may be empty) for which the system \eqref{syst} admits a solution. The projected Markov chain $\widecheck{X}$ thus possesses an \textcolor{black}{coclique} level structure for any \textcolor{black}{coclique} level function derived as described above.}
\end{theorem}

\begin{proof}
    Let $L$ be \textcolor{black}{a coclique} level function for $\widecheck{X}$. Then, by Section \ref{sec:IndLevelStructureApp}, $L$ has the form \eqref{ILFformula}, where $b \in \Z^{d-1}$ satisfies \eqref{syst} with $\Eu=\{\edge_k:L(\widecheck{v}_k)=1\}$ and $\Ed=\{\edge_k:L(\widecheck{v}_k)=-1\}$.
    Conversely, suppose that $b$ is a solution of the system \eqref{syst} for some partition $\{ \Eu, \Ed \}$ of the edge set of $\G$ (where $\Eu$ or $\Ed$ may be empty). Then, $b \in \Z^{d-1}$ by Lemma \ref{thm2a}, and  with $L(x)=b^T x$ for all $x \in \Z^{d-1}$, $L$ satisfies
    \eqref{indlevfunctiondef}.  Thus, $L$ is \textcolor{black}{a coclique} level function for $\widecheck{X}$. %\textcolor{black}{ Define $\ell$, $u$ as in \eqref{lowhighletters}, and $\L_z$ for $z=\ell,\dots,u$, as in \eqref{indlevstructure}. Then, the (ordered) partition $\L_{\ell}, \dots, \L_u$ is an \textcolor{black}{coclique} level structure for $\widecheck{X}$.}
\end{proof}

\begin{remark}\label{correlationpartitionstructure}
    For a given partition $\{ \Eu, \Ed \}$ of the edge set of $\G$, if the system \eqref{syst} has a solution, then by Lemma \ref{thm2a}, it has a unique solution. Furthermore, if we switch the labels of the two subsets, i.e., $\Eu$ becomes $\Ed$ and vice versa, the function  $L$ is replaced by $-L$, and the sets in the partition $\{\L_{\ell}, \dots, \L_u \}$ remain the same but the order is reversed. Therefore, we may consider the \textcolor{black}{coclique} level structures associated with $L$ and $-L$ to be the same.\\
   % 
   % the edge set partition remains the same. The associated \textcolor{black}{coclique} level structure will have the same levels, but ordered in the opposite way, i.e., switching from $\L_{\ell}, \ldots, \L_u$ to $\L_{u}, \ldots, \L_{\ell}$.
    %
   % Thus, a partition of the edge set of $\G$ uniquely determines a partition $\L_{\ell}, \ldots, \L_u$ of $\widecheck{X}$ (associated with $L$ or $-L$)
%   in a specific \textcolor{black}{coclique} level structure, whose order depends on the labels of the two subsets partitioning the edge set.
\end{remark}

{\color{black}
\noindent \textbf{Application to Example \ref{ex:illustrativeexample}:} Consider the SCRN introduced in \eqref{reacs-3speciesexample}, with associated graph $\G$ shown in Figure \ref{fig:illexample}(b). As previously mentioned, by applying Lemma \ref{lemmabis}, we know that this SCRN has a unique conservation vector $m = (1,1,1)^T$. This allows us to introduce a projected continuous time Markov chain $\widecheck{X}=\{
(X_1(t),X_2(t))^T:\: t \ge 0\}$, which keeps track of $(n_{\mathrm{W}},n_{\mathrm{Y}})$ through time, with transition vectors and associated rates given in \eqref{vectors3species} and \eqref{rates3species}, respectively.

Given that Assumption \ref{assumption:Unimolecular_change} is satisfied and the associated graph $\G$ is weakly connected, we can use Theorem \ref{thm2b} to determine all of the coclique level structures for $\widecheck{X}$. To this end, consider all of the possible partitions $\{ \Eu, \Ed \}$ of edges of $\G$ that could allow us to determine a coclique level structure. These partitions are the following:
\begin{align}
    &\Eu=\{\edge_1,\edge_2\},\Ed=\emptyset\;\;\;\;\mathrm{and}\;\;\;\;\Eu=\{\edge_1\},\Ed=\{\edge_2\}.
\end{align}
We did not consider the partition $\Eu=\emptyset, \Ed=\{\edge_1,\edge_2\}$ or the partition $\Eu=\{\edge_2\}, \Ed=\{\edge_1\}$ because, as explained in Remark \ref{correlationpartitionstructure}, the associated functions $L$ would be the opposite of the ones obtained for the partitions considered above and the resulting coclique level structures are considered to be the same.
 For each partition, the system of equations in \eqref{syst} has a unique solution, these being $(b_1,b_2)^T=(1,2)^T$ and $(b_1,b_2)^T=(1,0)^T$, respectively. Then, by applying Theorem \ref{thm2b}, we can conclude that the projected Markov chain $\widecheck{X}$ has two coclique level structures, with associated coclique level functions $L(x_1,x_2)=x_1+2x_2$ and $L(x_1,x_2)=x_1$, respectively.\\}

Now, we characterize when there exists \textcolor{black}{a coclique} level function for $\widecheck{X}$, assuming Assumption \ref{assumption:Unimolecular_change} and the associated graph $\G$ is weakly connected (see Theorem \ref{thm}).
%In this section, we provide {\color{black} necessary and sufficient??} conditions for the existence of \textcolor{black}{coclique} level structures for $\widecheck{X}$.
%
\textcolor{black}{In the course of this, we} show how to derive a simple \textcolor{black}{coclique} level function, when one exists, \textcolor{black}{which exploits bipartite structure of $\G$}. \textcolor{black}{There can be other coclique level functions and} we can use 
 %we can use Theorem \ref{suffcondthm} %to rule out certain choices of %partition $\Eu$, %$\Ed$?? the partition %should write as $\{ \Eu, %\Ed \}$??} Then, we can use 
 Theorem \ref{thm2b} above 
 %introduced in the previous section
 to identify all possible \textcolor{black}{coclique} level functions. 
%
%Conversely, when these conditions {\color{black} what conditions??} are not verified, no \textcolor{black}{coclique} level structures exist.

%
%\begin{definition}%\label{def3}
%{\color{black} Before introducing these conditions, }
   Assume the SCRN satisfies Assumption \ref{assumption:Unimolecular_change} and its associated graph $\G$ is weakly connected. Consider a weakly directed cycle $c$ in $\G$, i.e., a directed subgraph of $\G$ whose underlying undirected graph is a cycle. We abbreviate weakly directed cycle as wd-cycle\footnote{Note that we allow two vertices connected by edges $\edge$ and the reverse of $\edge$ to be a wd-cycle.}. Choose an orientation for the cycle. Then, define the vector $\vartheta_c\in \Z^n$ associated with the wd-cycle $c$ such that for $k=1,\dots,n$,
	\begin{equation}\vartheta_c(k)=\begin{cases} +1 & \mbox{if $\edge_k$ is part of the wd-cycle and $\edge_k$ is in}\\
 & \mbox{the direction of the orientation of the wd-cycle,}\\-1 & \mbox{if $\edge_k$ is part of the wd-cycle and $\edge_k$ is in }\\	 & \mbox{the opposite direction of the orientation of the wd-cycle,}\\	
	0 & \mbox{if $\edge_k$ is not part of the wd-cycle.} \end{cases}
 \end{equation}
\begin{theorem}\label{suffcondthm}
	Consider a SCRN satisfying Assumption \ref{assumption:Unimolecular_change} and assume its associated graph $\G$ is weakly connected. 
 %Then, given a basis for $Null(S)$, that is $\{\vartheta_1,\vartheta_2,\dots,\vartheta_c,\dots,\vartheta_{n-d+1}\}$, and
 Let $\mathcal{C}$ be the set of all wd-cycles in $\G$. For each $c \in \mathcal{C}$, we have $\vartheta_c\in \ker(S)$. Furthermore, given a partition $\{\Eu,\Ed\}$ of the edge set of $\G$ (where one of $\Eu,\Ed$ may be empty), we can define a vector $w$ such that $w_k=+1$ if $\edge_k \in \Eu$ and $w_k=-1$ if $\edge_k \in \Ed$.
 %defined as in (\ref{systMATRIX}), 
 Then, if $w^T\vartheta_c \ne 0$ for some $c \in \mathcal{C}$,  the system (\ref{syst}) does not admit a solution $b \in \Z^{d-1}$.	%{\color{black} for $\{\Eu,\Ed\}$??}\textcolor{black}{yes, I would add it.}. 
 %\textcolor{black}{also there is $\E$ and $\mathcal{E}$, two kinds of curly $\edge$??}\textcolor{black}{Do we use $\E$ somewhere? Because I could not find any $\E$ in the paper.} {\color{black} $\E$ is the expectation, in the mean first passage time formula}
\end{theorem}

\begin{proof}
The $i^{th}$ entry of the vector $S\vartheta_c$ represents the amount of species $\mathrm{S}_i$ consumed/produced after all the reactions associated with the wd-cycle $c$ are triggered, where the signs in  $\vartheta_c$ ensure that the edges are followed in the direction of the chosen orientation.
Since $c$ is a wd-cycle and Assumption \ref{assumption:Unimolecular_change} is satisfied, $S\vartheta_c=0$ and then $\vartheta_c\in \ker(S)$.

 %\textcolor{black}{The reason why $\vartheta_c\in \ker(S)$ is that the product between the $i^{th}$ row of $S$ and $\vartheta_c$ corresponds to the amount of the species $\mathrm{S}_i$ consumed/produced after all the reactions associated with the wd-cycle $c$ are triggered, where the signs in  $\vartheta_c$ ensure that the edges are followed in the direction of the chosen orientation. Since $c$ is a wd-cycle, this means that the overall change in the amount of $\mathrm{S}_i$ is zero. Since this is true for all $\mathrm{S}_i\in \Ss$, then $S\vartheta_c=0$. TRY TO MAKE IT SHORTER.ADD JUST WHAT THE PRODUCT REPRESENTS, AND SINCE IT IS A CYCLE, THEN THE PRODUCT WILL BE ZERO}

Since $\widecheck{S}^T$ is $S^T$ without the last column, we have that $\text{range}(\widecheck{S}^T) \subseteq \text{range}(S^T)   = \ker(S)^\perp$. 
%The system \eqref{syst} has a solution if and only if $w \in \text{range}(\widecheck{S}^T)$. If $w \in \text{range}(\widecheck{S}^T)$, then $w \in \text{range}(S^T)$. If $w \in \text{range}(S^T)$, then $x \cdot w= 0$ for all $x \in Null(S)$. 
Since $\vartheta_c\in \ker(S)$ for each $c \in \mathcal{C}$, we have that $\text{span}\{\vartheta_c: c \in \mathcal{C} \} \subset \ker (S)$. 
It follows that if $w^T\vartheta_c  \neq 0$ for some $c \in \mathcal{C}$, then $w \notin \ker(S)^{\perp} \supseteq \text{range}(\widecheck{S}^T)$, and thus (\ref{syst}) does not admit a solution.
% The system \eqref{syst} has a solution if and only if $w \in \text{range}(\widecheck{S}^T)$. 
\iffalse
	By definition we know that, given a matrix $S$, the orthogonal complement of $Null(S)$ is $\text{range}(S^T)$, that is, for the stochiometry matrix $S\in \Z^{d \times n}$, $Null(S) \bigoplus \text{range}(S^T)=\R^n$. This implies that, for any $w\in \text{range}(S^T)$ and any $x \in Null(S)$, we have that $x \cdot w= 0$. Now, since we know that $\{\vartheta_1,\vartheta_2,\dots,\vartheta_c,\dots,\vartheta_{n-d+1}\}$ is a basis for $Null(S)$, any $x \in Null(S)$ can be expressed as a linear combination of the vector of the basis and then $x \cdot w= 0$ if and only if $\vartheta_c \cdot w= 0$. Now, let us introduce $\widecheck{S}\in \Z^{(d-1) \times n}$ that corresponds to the stoichiometric matrix associated with the system without the last row and that is full rank (see Lemma \ref{thm2a}) and let us consider a possible partition of the edges, \{$\Eu$, $\Ed$\}, and the corresponding vector $w$. Then, we can conclude that for any $c\in\{1,\dots,n-d+1\}$, if $\vartheta_c \cdot w=0$, the system (\ref{syst}) admits a unique solution, and if $\vartheta_c \cdot w\ne 0$, the system (\ref{syst}) does not admit any solution.
 \fi
% 
\end{proof}

\begin{corollary}\label{thm2bis}
	Consider a SCRN satisfying Assumption \ref{assumption:Unimolecular_change} and assume its associated graph $\G$ is weakly connected. Suppose $\{\Eu,\Ed\}$ is a partition of the edge set of $G$ (where one of $\Eu$ and $\Ed$ can be empty). If %\st{there are two reaction vectors such that one is the opposite of the other one (i.e.,} 
 there are $v_k,v_{k'} \in \V$ such that $v_k=-v_{k'}$ and the edges associated with them, $\edge_k$ and $\edge_{k'}$, belong to the same subset $\Eu$ or $\Ed$, then (\ref{syst}) does not admit a solution $b \in \Z^{d-1}$, and so there is no \textcolor{black}{coclique} level structure for this partition.
%	Let us consider a SCRN as introduced in Lemma \ref{nullspacebasis}. If there are reversible reactions and if the edges associated with the forward and backward steps belong to the same subset $\Eu$ or $\Ed$ (with $\Eu$, $\Ed$ defined in Theorem \ref{thm2a}), then system (\ref{syst}) does not admit solution.
\end{corollary}
%
%\begin{proof}
%OLD
%	By renumbering if necessary, we may suppose there are two reaction vectors, $v_1$ and $v_2$, such that $v_1=-v_2$, with corresponding edges ($\edge_1$ and $\edge_2$, respectively). Assume that both of them belong to $\Eu$. Consider the wd-cycle $c$ given by $\edge_1$ and $\edge_2$, with corresponding vector $\vartheta_c$. Then, $\vartheta_c(k)=0$ for any $k=3,4,\dots,n$ and $\vartheta_c(1)=\vartheta_c(2)= +1$ or $\vartheta_c(1)=\vartheta_c(2)= -1$, depending on the orientation of the cycle. Then, since both the edges belong to $\Eu$, we have that $w(1)=w(2)=1$ where $w$ is defined as in \eqref{systMATRIX}, and this implies that $\abs{\vartheta_c \cdot w} = 2 \neq 0$. Then, the result follows from Theorem \ref{suffcondthm}. A similar reasoning can be used for the case where $\edge_1,\edge_2\in\Ed$.
 %, this proves that, if the edges associated with reaction vectors such that one is the opposite of the other one belong to the same subset $\Eu$ or $\Ed$, then system (\ref{syst}) does not admit any solution.
%\end{proof}

\begin{proof}
	Suppose that $\edge_k,\edge_{k'} \in \Eu$ have the properties described. Consider the wd-cycle $c$ given by $\edge_k$ and $\edge_{k'}$, with corresponding vector $\vartheta_c$. Then, $\vartheta_c(\ell)=0$ for $\ell \notin \{k,k'\}$ and, depending on the orientation of the cycle, $\vartheta_c(k)=\vartheta_c(k')= +1$ or $\vartheta_c(k)=\vartheta_c(k')= -1$. Then, since both the edges belong to $\Eu$, we have $w$, defined as in \eqref{systMATRIX}, to be such that $w(k)=w(k')=1$ and $w(\ell)=0$ for $\ell \notin \{k,k'\}$. This implies that $\abs{w^T\vartheta_c} = 2 \neq 0$. Then, the result follows from Theorem \ref{suffcondthm}. Similar reasoning can be used for the case where $\edge_k,\edge_{k'}\in\Ed$.
 %, this proves that, if the edges associated with reaction vectors such that one is the opposite of the other one belong to the same subset $\Eu$ or $\Ed$, then system (\ref{syst}) does not admit any solution.
\end{proof}

\begin{remark}
\label{rmk:reversible}
Based on the results of Corollary \ref{thm2bis}, given a partition $\{ \Eu, \Ed \}$, we can obtain \textcolor{black}{a coclique} level function only if edges $\edge_k$ and $\edge_{k'}$ associated with $v_k,v_{k'} \in \V$ such that $v_k=-v_{k'}$ belong to different members of the partition. If this property holds, then the equation \eqref{syst} associated with $v_k$ is the opposite of the one associated with $v_{k'}$.
%(that is, we can obtain one equation by multiplying the other one by $-1$). 
This implies that the two  equations are providing the same information. Thus, in order to solve system (\ref{syst}) and reduce the computational cost, we can just write the system for $\widehat{\V}$, in which $\widehat{\V} \subseteq \V$ is \textcolor{black}{a} maximal set of reaction vectors such that, for any $v_k,v_{k'} \in \widehat{\V}$, $v_k \ne - v_{k'}$.
%$\widehat{\mathcal{R}}$ is the set of the $\bar p$ reactions given by all the irreversible reactions and only the forward steps of the reversible ones.
\end{remark}

%\begin{corollary}\label{thm3}
%	Consider a SCRN satisfying Assumption \ref{assumption:Unimolecular_change} and assume its graph $\G$ is weakly connected. If the graph $\G$ is not bipartite, then the system (\ref{syst}) does not admit any solutions\textcolor{black}{, indicating the nonexistence of an \textcolor{black}{coclique} level structure.}
%\end{corollary}
%
%\begin{proof}
%	If the graph $\G$ is not bipartite, it is well known that it contains at least one odd wd-cycle, $c^{odd}$. Now, let us consider the vector $\vartheta_{c^{odd}}$ associated with this cycle. By definition, $\vartheta_{c^{odd}}$ has all of its the entries equal to zero except for an odd number of entries, associated with the edges in the wd-cycle $c^{odd}$. These latter entries are in $\{-1,+1\}$. Furthermore, let us define a vector $w$ as done in (\ref{systMATRIX}) that has entries that can only be $-1$ and $+1$. This means that the product $\vartheta_{c^{odd}} \cdot w$ will result in a sum of an odd number of terms that can be $+1$ or $-1$ and this can never be equal to zero. By Theorem \ref{suffcondthm}, this implies that system (\ref{syst}) does not admit a solution.	
%	% Given any two equations we define their sum to be the equation whose left hand side is the sum of the two left hand sides, and whose right hand side is the sum of the two right hand sides.
%\end{proof}

Finally, we introduce a theorem demonstrating that for a SCRN satisfying Assumption \ref{assumption:Unimolecular_change} and having a weakly connected graph $\G$,  there is \textcolor{black}{a coclique} level structure for $\widecheck X$ if and only if $\G$ is bipartite.

\begin{theorem}\label{thm}
	Consider a SCRN satisfying Assumption \ref{assumption:Unimolecular_change} and suppose the associated graph $\G$ is weakly connected. 
 If the graph $\G$ is bipartite,
 %Consider the associated projected continuous time Markov chain, $\widecheck{X}$. 
 then there exists \textcolor{black}{a coclique} level function for the projected continuous time Markov chain $\widecheck{X}$. In this case, by relabeling vertices if necessary, there are two disjoint, non-empty \textcolor{black}{coclique}s of vertices in $\G$, $B=\{1,\dots,\bar i\}$ and $C=\{\bar i+1,\dots,d\}$, and then 
 \begin{equation}
 \label{Lformula}
 L(\widecheck{x})=L(x_1,\dots,x_{\bar i},x_{\bar i+1},\dots,x_{d-1})=x_1+\dots+x_{\bar i}
 \end{equation} 
 defines \textcolor{black}{a coclique} level function for $\widecheck{X}$.
 If the graph $\G$ is not bipartite, 
 %{\color{black} version 1:} then the system (\ref{syst}) does not admit any solutions for any partition $\{ \Eu, \Ed \}$ of the edge set of $\G$ (where one of $\Eu$ and $\Ed$ may be empty), indicating the nonexistence of any \textcolor{black}{coclique} level structure. {\color{black} version 2:} 
 there is no \textcolor{black}{coclique} level structure for $\widecheck{X}$.
\end{theorem}

\begin{proof}
We start by assuming the graph $\G$ is bipartite, and show there exists \textcolor{black}{a coclique} level function for the projected continuous time Markov chain $\widecheck{X}$ as described in the theorem.
Suppose $\G$ is bipartite. By relabeling vertices if necessary, there are sets $B$ and $C$ as described in the theorem.
Let $\Eu = \{ \edge_k = (i,j) : v_k \in \V, i \in C \text{ and } j \in B\}$ and $\Ed = \{ \edge_k = (i,j) : v_k \in \V, i \in B \text{ and } j \in C\}$. Note that this $\{ \Eu, \Ed \}$ is a partition of the edge set of $\G$ where one of $\Eu$ or $\Ed$ could be empty.
Let $b \in \Z^{d-1}$ be the vector where the first $\bar i$ entries are 1 and the remaining $d-1-\bar i$ entries are zero and define $L(x) = b^Tx$ for $x \in \Z^{d-1}$.
Since $\G$ is bipartite and $B$ and $C$ are two \textcolor{black}{coclique}s, each reaction corresponds either to an edge from a vertex in $C$ to a vertex in $B$ or an edge from a vertex in $B$ to a vertex in $C$. In particular, when a reaction $v_k \in \V$ is triggered, if $\edge_k \in \Eu$ then $L(\widecheck{v}_k)=+1$ and if $\edge_k \in \Ed$ then $L(\widecheck{v}_k)=-1$. 
%
%the quantity $x_1+\dots+x_{\bar i}$ is increased by one if $\edge_k \in \Eu$ or decreased by one if $\edge_k \in \Ed$. 
%Thus, $b$ satisfies \eqref{syst}.
%%
%Then, by Theorem \ref{thm2b}, the function $L: \Z^{d-1} \to \Z$ given by \eqref{Lformula} is an \textcolor{black}{coclique} level function for $\widecheck{X}$.
Thus, $L$ satisfies \eqref{indlevfunctiondef} and is \textcolor{black}{a coclique} level function for $\widecheck{X}$.

Now suppose the graph $\G$ is not bipartite. 
%then the system (\ref{syst}) does not admit any solutions, indicating the nonexistence of an \textcolor{black}{coclique} level structure.
%%
%If the graph $\G$ is not bipartite, 
Then $\G$ contains at least one odd wd-cycle, $c^{odd}$ \textcolor{black}{(see Theorem 5.1 on page 27 of \cite{WilsonRJ})}. Consider the vector $\vartheta_{c^{odd}}$ associated with this cycle. By definition, $\vartheta_{c^{odd}}$ has all of its entries equal to zero except for an odd number of entries, associated with the edges in the wd-cycle $c^{odd}$, which take values in $\{-1,+1\}$. For any partition $\{\Eu,\Ed\}$ of the edge set of $\G$ (where one of $\Eu$ or $\Ed$ may be empty), and for $w$ as in \eqref{systMATRIX}, the product $w^T\vartheta_{c^{odd}} $ will result in a sum of an odd number of terms, each of which is $+1$ or $-1$, and thus the sum can never be equal to zero. By Theorem \ref{suffcondthm}, this implies that the system (\ref{syst}) does not admit a solution. By Theorem \ref{thm2b}, there cannot be \textcolor{black}{a coclique} level function for $\widecheck{X}$.
%for any partition $\{ \Eu, \Ed \}$ of the edge set of $\G$ (where one of $\Eu$ and $\Ed$ may be empty).	
\end{proof}

{\color{black}
\noindent \textbf{Application to Example \ref{ex:illustrativeexample}:} Consider the SCRN introduced in \eqref{reacs-3speciesexample}, with associated graph $\G$ shown in Figure \ref{fig:illexample}(b). Earlier, we showed how, by applying Theorem \ref{thm2b}, one can determine all of the coclique level functions for $\widecheck{X}$. On the other hand, by inspecting $\G$, it is also possible to note that the graph is bipartite, and thus we can apply Theorem \ref{thm} to directly identify one of its coclique level functions. Specifically, the two disjoint, non-empty cocliques of vertices in $\G$ are $B=\{1\}$ and $C=\{2, 3\}$, and then $L(x)=x_1$ defines a coclique level function for $\widecheck{X}$.}

\subsection{Algorithm to determine all \textcolor{black}{coclique} level functions 
associated with a SCRN satisfying Assumption \ref{assumption:Unimolecular_change} and whose associated graph $\G$ is weakly connected}\label{algIndLevFac}

\begin{figure}[H]%[h!]
            \centering
            \includegraphics[scale=0.43]{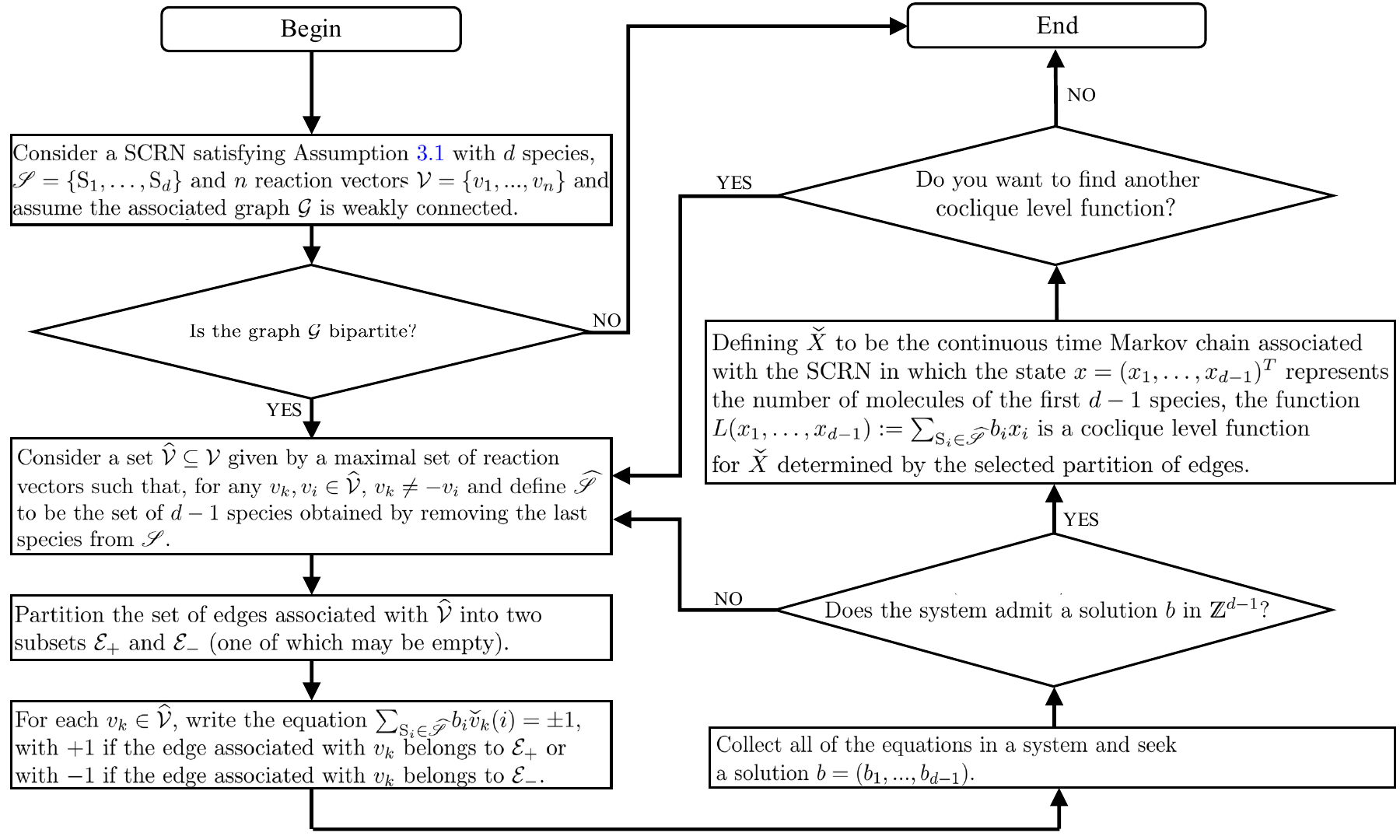}
            \caption{\small { \bf Key steps of the algorithm for identifying all \textcolor{black}{coclique} level functions for the projected continuous time Markov chain $\widecheck X$ 
associated with a SCRN under 
Assumption \ref{assumption:Unimolecular_change} and assuming that $\G$ is weakly connected.}
}
   \label{Alg}
        \end{figure}

Based on the results given in Theorems \ref{thm2b} and \ref{thm}, given a SCRN satisfying Assumption \ref{assumption:Unimolecular_change} and assuming that $\G$ is weakly connected, we can obtain an algorithm that allows us to find all of the \textcolor{black}{coclique} level functions for the projected continuous time Markov chain $\widecheck{X}$ associated with the SCRN. The steps of the algorithm are described in Figure \ref{Alg}. {\color{black} The first step of the algorithm involves checking whether a graph is bipartite. The reader may refer to Section 3.4 in Kleinberg and Tardos \cite{bib:JE2005} for a discussion about a Breadth-First Search implementation to check whether a graph is bipartite, whose run time is linear in terms of the number of vertices in the graph (or, in our case, the number of species in the SCRN). 
}

{\color{black}
As a second step, the algorithm requires selecting a subset $\widehat{\V} \subseteq \V$ defined as a maximal set of reaction vectors such that $v_k \ne -v_i$ for any $v_k, v_i \in \widehat{\V}$. It is important to note that, once such a subset is identified, the specific ordering of the vectors in $\widehat{\V}$ does not affect the outcome of the algorithm. Moreover, the choice of $\widehat{\V}$ is uniquely determined up to the direction of the vectors and ensures that all directions allowed by the reactions are considered.}

\begin{remark}
    Note that the algorithm described here is for a specific ordering of the species associated with our specification of the SCRN.
%   
%     All possible \textcolor{black}{coclique} level structures for the SCRN may be obtained by considering all possible orderings of the species.

% If the ordering of the species is permuted, one can find equivalent \textcolor{black}{coclique} level structures 
    
\end{remark}

\subsection{Using \textcolor{black}{coclique} level structure to bound MFPTs}
\label{indeplevelstructure}
%\label{uplowbounds}

As mentioned in the introduction, obtaining an explicit analytical expression for the mean first passage time (MFPT) for one-dimensional birth-death processes is relatively straightforward by using the formula \eqref{MFPTsystem} (see SI - Section \ref{sec:Appendix1Dmodel}). However, it is typically very complicated to obtain an explicit expression for the MFPT for more complex Markov chains, especially in dimension greater than one. In this section, we will show how to obtain closed form formulas for upper and lower bounds for the MFPTs of continuous time Markov chains $\widecheck{X}$ associated with SCRNs having \textcolor{black}{a coclique} level structure.

%We first introduce the definition of $independent\; level\; structure$ and then determine formulas for upper and lower bounds for the MFPTs.

To this end, consider the projected continuous time Markov chain $\widecheck{X}$ with finite state space $\widecheck{\X} \subseteq \Z^{d-1}_{+}$, associated with a SCRN satisfying Assumption \ref{assumption:Unimolecular_change} and for which $\G$ is weakly connected
%as defined in Remark \ref{remark1},
(as defined in Section \ref{sec:IndLevelStructureApp}), %for $d\geq 2$,
and having infinitesimal generator $\widecheck{Q}$. Suppose $\widecheck{X}$ has \textcolor{black}{a coclique} level function 
%There is a finite set of possible transition directions denoted by $\{\widecheck{v}_1,\dots,\widecheck{v}_n\} \subseteq  \Z^{d-1} \setminus \{0\}$.
%
%Consider also the \textcolor{black}{coclique} level function 
$L: \Z^{d-1} \rightarrow \Z$, defined as in \eqref{indlevfunctiondef}, with (ordered) partition $\{\L_{\ell},\L_{\ell+1}, \dots, \L_{u-1}, \L_u \}$ of $\widecheck{\mathcal{X}}$, defined as in \eqref{indlevstructure} {\color{black} and where $\ell < u$}.
%Let $G_{+} :=  \{k:L(\widecheck{v}_k) =1\}$ and $G_{-} :=  \{k:L(\widecheck{v}_k) =-1\}$. 
Given such \textcolor{black}{a coclique} level structure for $\widecheck{X}$, we will determine analytical expressions for upper and lower bounds for both the MFPT for $\widecheck{X}$ from $\L_{\ell}$ to $\L_{u}$ and the MFPT for $\widecheck{X}$ from $\L_{u}$ to $\L_{\ell}$. We first focus on the MFPT from $\L_{\ell}$ to $\L_{u}$.

{\color{black}
Let
\begin{equation*}
    G_{+} \vcentcolon=  \{k:L(\widecheck{v}_k) =1\} \quad \text{ and } \quad G_{-} \vcentcolon=  \{k:L(\widecheck{v}_k) =-1\}.
\end{equation*}
}For $\ell < u$, $\ell \le z \le u$ and $\widecheck{x} \in \L_{z}$, define the ``rate of increase'' $\lambda_z(\widecheck{x})$ and the ``rate of decrease'' $\gamma_z(\widecheck{x})$ as follows: 
\begin{equation}\label{lambdagamma}\nonumber
	    \lambda_z(\widecheck{x})=\sum_{k \in G_{+}}\widecheck{Q}_{\widecheck{x},\widecheck{x}+\widecheck{v}_k},\;\;\;\;\gamma_z(\widecheck{x})=\sum_{k \in G_{-}}\widecheck{Q}_{\widecheck{x},\widecheck{x}+\widecheck{v}_k}.
\end{equation}
Then, for each $z\in \{\ell,\ell+1,\dots,u-1,u\}$, let
	\begin{equation}\label{minmaxrates}
	\begin{aligned}
	&\lambda^M_z=\max_{\widecheck{x}\in \L_z}\lambda_z(\widecheck{x}),\;\;\;\;\lambda^m_z=\min_{\widecheck{x}\in \L_z}\lambda_z(\widecheck{x}),\;\;\;\;\gamma^M_z=\max_{\widecheck{x}\in \L_z}\gamma_z(\widecheck{x}),\;\;\;\;\gamma^m_z=\min_{\widecheck{x}\in \L_z}\gamma_z(\widecheck{x}).
	\end{aligned}
	\end{equation}
    %Introduce continuous time Markov chains $\Breve{X}$ and $\invbreve{X}$ with infinitesimal generators $\Breve{Q}$ and $\invbreve{Q}$, respectively, such that, for $z\in \{\ell, \ell +1 ,\dots, u-1,u\}$, $\widecheck{x}\in \L_z$, $\Breve{Q}_{\widecheck{x},\widecheck{x}+\widecheck{v}_k}=\frac{\lambda^M_z}{|G_{+}|}$ for $k \in G_{+}$,  $\Breve{Q}_{\widecheck{x},\widecheck{x}+\widecheck{v}_k}=\frac{\gamma^m_z}{|G_{-}|}$ for $k \in G_{-}$, $\invbreve{Q}_{\widecheck{x},\widecheck{x}+\widecheck{v}_k}=\frac{\lambda^m_z}{|G_{+}|}$ for $k \in G_{+}$, and $\invbreve{Q}_{\widecheck{x},\widecheck{x}+\widecheck{v}_k}=\frac{\gamma^M_z}{|G_{-}|}$ for $k \in G_{-}$.
    {\color{black} 
    We will define continuous time Markov chains $\Breve{X}$ and $\invbreve{X}$ with the same state space as $\widecheck{X}$.
    For $\ell \le z \le u$ and $\widecheck{x} \in \L_{z}$, let
    \begin{equation}
\label{G+-(x)}
        G_{+} (\widecheck{x}) = \{k \in G_{+}: \widecheck{x}+\widecheck{v}_k \in \widecheck{\X} \} \quad \text{ and } \quad G_{-} (\widecheck{x}) = \{k \in G_{-}: \widecheck{x}+\widecheck{v}_k \in \widecheck{\X} \}.
    \end{equation}
    We assume that $G_{+} (\widecheck{x}) : \widecheck{x} \in \widecheck{\X} \setminus \L_u$ are all non-empty sets when $G_+$ is non-empty, and $G_{-} (\widecheck{x}) : \widecheck{x} \in \widecheck{\X} \setminus \L_\ell$ are all non-empty sets when $G_-$ is non-empty. Note that $G_{+} (\widecheck{x})$ is empty when $\widecheck{x} \in \L_u$ and $G_{-} (\widecheck{x})$ is empty when $\widecheck{x} \in \L_\ell$. Then, we define the infinitesimal generator $\Breve{Q}$ for $\Breve{X}$ such that the only positive entries of $\Breve{Q}$ are given by 
    \begin{equation*}
        \Breve{Q}_{\widecheck{x},\widecheck{x}+\widecheck{v}_k} = \frac{\lambda^M_z}{|G_{+} (\widecheck{x})|} \text{ for } k \in G_{+} (\widecheck{x}) \quad \text{ and } \quad \Breve{Q}_{\widecheck{x},\widecheck{x}+\widecheck{v}_k} = \frac{\gamma^m_z}{|G_{-}(\widecheck{x})|} \text{ for } k \in G_{-} (\widecheck{x}),
    \end{equation*}
    for $\ell \le z \le u$ and $\widecheck{x} \in \L_{z}$. Similarly, we define the infinitesimal generator $\invbreve{Q}$ for $\invbreve{X}$ such that the only positive entries of $\invbreve{Q}$ are given by 
    \begin{equation*}
        \invbreve{Q}_{\widecheck{x},\widecheck{x}+\widecheck{v}_k} = \frac{\lambda^m_z}{|G_{+} (\widecheck{x})|} \text{ for } k \in G_{+} (\widecheck{x}) \quad \text{ and } \quad \invbreve{Q}_{\widecheck{x},\widecheck{x}+\widecheck{v}_k}=\frac{\gamma^M_z}{|G_{-} (\widecheck{x})|} \text{ for } k \in G_{-} (\widecheck{x}),
    \end{equation*}
    for $\ell \le z \le u$ and $\widecheck{x} \in \L_{z}$. 
    }
 
Then, by the comparison theorem\textcolor{black}{s}, i.e., Theorems 3.3 \textcolor{black}{and 3.4} in \cite{Monotonicitypaper} \textcolor{black}{(see SI - Section \ref{compthmpaper_thm})}, using the matrix $A$ equal to $b^T$ (associated with the \textcolor{black}{coclique} level function $L$), we have that for $\widecheck{x}, \Breve{x}$ and $\invbreve{x}$ in $\L_{\ell}$, we can realize $\widecheck{X}$ and $\Breve{X}$ with infinitesimal generators $\widecheck{Q}$ and $\Breve{Q}$ (resp. $\widecheck{X}$ and $\invbreve{X}$ with infinitesimal generators $\widecheck{Q}$ and $\invbreve{Q}$) on the same probability space such that $L(\widecheck{X}) \le L(\Breve{X})$, $\widecheck{X}(0) = \widecheck{x}$, $\Breve{X}(0) = \Breve{x}$ a.s. (resp. $L(\invbreve{X}) \le L(\widecheck{X})$, $\invbreve{X}(0) = \invbreve{x}$, $\widecheck{X}(0) = \widecheck{x}$ a.s.).
Then, for 
$\widecheck{\tau}_{u} = \inf\{ t \geq 0:\: \widecheck{X}(t) \in \L_u \} = \inf\{ t \geq 0:\: L(\widecheck{X}(t)) = u \}$, $\Breve{\tau}_{u} = \inf\{ t \geq 0:\: \Breve{X}(t) \in \L_u \} = \inf\{ t \geq 0:\: L(\Breve{X}(t)) = u \}$, and $\invbreve{\tau}_{u} = \inf\{ t \geq 0:\: \invbreve{X}(t) \in \L_u \} = \inf\{ t \geq 0:\: L(\invbreve{X}(t)) = u \}$, we have 
\begin{equation}
\label{upandlowbounds}
    \E_{\Breve{x}}[\Breve{\tau}_{u}]\le
     \E_{\widecheck{x}}[\widecheck{\tau}_{u}] \le 
     \E_{\invbreve{x}}[\invbreve{\tau}_{u}].
\end{equation}
By the choice of the rates for $\Breve{X}$ and $\invbreve{X}$, we have that $L(\Breve{X})$ and $L(\invbreve X)$ are continuous time Markov chains. In fact, they are simple birth-death processes. It follows that the upper and lower bounds in \eqref{upandlowbounds} can be explicitly evaluated. For the lower bound, 
{\color{black} suppose that $\lambda^M_z$ is positive for $\ell \leq z \leq u-1$. We} 
replace $\lambda_i$ with $\lambda_z^M$ and $\gamma_i$ with $\gamma_z^m$, respectively, in formula \eqref{formula3}, to obtain
%{\color{teal} BEGIN OLD}
%\begin{equation*}
%\begin{aligned}
%\E_{\Breve{x}}[\Breve{\tau}_{u}]=\Breve{h}_{\ell,u}=\frac{\tilde r^m_{u-(\ell+1)}}{\lambda^M_{\ell}}\left(1+\sum_{j=1}^{u-(\ell+1)}\frac{1}{\tilde r^m_j}\right)+\sum_{i=2}^{u-(\ell+1)}\left[\frac{\tilde r^m_{i-1}}{\lambda^M_{u-i}}\left(1+\sum_{j=1}^{i-1}\frac{1}{\tilde r^m_j}\right)\right]+\frac{1}{\lambda^M_{u-1}},
%\end{aligned}
%\end{equation*}
%with $\tilde r^m_j=\frac{\gamma^m_{u-1}\gamma^m_{u-2}\dots\gamma^m_{u-j}}{\lambda^M_{u-1}\lambda^M_{u-2}\dots\lambda^M_{u-j}}$. 
%
%{\color{teal} END OLD}
%
%{\color{teal} BEGIN NEW}
{\color{black}
\begin{equation}
\label{formula3lowbound}
    \E_{\Breve{x}}[\Breve{\tau}_{u}]=\Breve{h}_{\ell,u}=\frac{1}{\lambda^M_{u-1}} + \sum_{i=\ell}^{u-2} \frac{1}{\lambda^M_{i}} \left( 1 + \sum_{j=i+1}^{u-1} \frac{\gamma^m_{i+1}\dots\gamma^m_{j}}{\lambda^M_{i+1}\dots\lambda^M_{j}} \right).
\end{equation}
}
%
%{\color{teal} END NEW}
% From the SI:
% \begin{eqnarray*}
%     h_{0,\upb}&=&\frac{1}{\lambda_{0}}\left(1+\frac{\gamma_1}{\lambda_1}+\frac{\gamma_1\gamma_2}{\lambda_1\lambda_2}+\dots+\tilde r_{\upb-1}\right)\\
%     && +\frac{1}{\lambda_{1}}\left(1+\frac{\gamma_2}{\lambda_2}+\frac{\gamma_2\gamma_3}{\lambda_2\lambda_3}+\dots+\tilde r_{\upb-2}\right)+\dots+\frac{1}{\lambda_{\upb-1}}\\
% \end{eqnarray*}
% and
% \begin{equation*}
%     \tilde r_j=\frac{\gamma_{\upb-1}\gamma_{\upb-2}\dots\gamma_{\upb-j}}{\lambda_{\upb-1}\lambda_{\upb-2}\dots\lambda_{\upb-j}}=\frac{\gamma_{\upb-j}\gamma_{\upb-j+1}\dots\gamma_{\upb-1}}{\lambda_{\upb-j}\lambda_{\upb-j+1}\dots\lambda_{\upb-1}}
% \end{equation*}
% Then, 
% \begin{eqnarray*}
%     h_{\ell,u}&=&\frac{1}{\lambda_{\ell}}\left(1+\frac{\gamma_{\ell+1}}{\lambda_{\ell+1}}+\frac{\gamma_{\ell+1}\gamma_{\ell+2}}{\lambda_{\ell+1}\lambda_{\ell+2}}+\dots+\frac{\gamma_{\ell+1}\gamma_{\ell+2}\dots\gamma_{u-1}}{\lambda_{\ell+1}\lambda_{\ell+2}\dots\lambda_{u-1}}\right)\\
%     && +\frac{1}{\lambda_{\ell+1}}\left(1+\frac{\gamma_{\ell+2}}{\lambda_{\ell+2}}+\frac{\gamma_{\ell+2}\gamma_{\ell+3}}{\lambda_{\ell+2}\lambda_{\ell+3}}+\dots+\frac{\gamma_{\ell+2}\gamma_{\ell+3}\dots\gamma_{u-1}}{\lambda_{\ell+2}\lambda_{\ell+3}\dots\lambda_{u-1}}\right)+\dots+\frac{1}{\lambda_{u-1}}
% \end{eqnarray*}
% Thus,
% \begin{equation*}
%     h_{\ell,u}=\frac{1}{\lambda_{u-1}} + \sum_{i=\ell}^{u-2} \frac{1}{\lambda_{i}} \left( 1 + \sum_{j=i+1}^{u-1} \frac{\gamma_{i+1}\dots\gamma_{j}}{\lambda_{i+1}\dots\lambda_{j}} \right)
% \end{equation*}
For the upper bound, 
{\color{black} suppose that $\lambda^m_z$ is positive for $\ell \leq z \leq u-1$. We} 
replace $\lambda_i$ with $\lambda_z^m$ and $\gamma_i$ with $\gamma_z^M$, respectively, in formula \eqref{formula3}, to obtain
%
%{\color{teal} BEGIN OLD}
%\begin{equation*}
%\begin{aligned}
%\E_{\invbreve{x}}[\invbreve{\tau}_{u}]=\invbreve{h}_{\ell,u}=\frac{\tilde r^M_{u-(\ell+1)}}{\lambda^m_{\ell}}\left(1+\sum_{j=1}^{u-(\ell+1)}\frac{1}{\tilde r^M_j}\right)+\sum_{i=2}^{u-(\ell+1)}\left[\frac{\tilde r^M_{i-1}}{\lambda^m_{u-i}}\left(1+\sum_{j=1}^{i-1}\frac{1}{\tilde r^M_j}\right)\right]+\frac{1}{\lambda^m_{u-1}}, 
%\end{aligned}
%\end{equation*}
%with $\tilde r^M_j=\frac{\gamma^M_{u-1}\gamma^M_{u-2}\dots\gamma^M_{u-j}}{\lambda^m_{u-1}\lambda^m_{u-2}\dots\lambda^m_{u-j}}$. 
%{\color{teal} If $\lambda^m_z=0$ for some $\ell \leq z \leq u-1$ or $\gamma^M_z=0$ for some $\ell+1 \leq z \leq u$, then one may observe that $\invbreve{h}_{\ell,u} = \infty$ in \eqref{formula3upbound}.}
%
%{\color{teal} END OLD}
%
%{\color{teal} BEGIN NEW}
%
{\color{black}
\begin{equation}
\label{formula3upbound}
    \E_{\invbreve{x}}[\invbreve{\tau}_{u}]=\invbreve{h}_{\ell,u}=\frac{1}{\lambda^m_{u-1}} + \sum_{i=\ell}^{u-2} \frac{1}{\lambda^m_{i}} \left( 1 + \sum_{j=i+1}^{u-1} \frac{\gamma^M_{i+1}\dots\gamma^M_{j}}{\lambda^m_{i+1}\dots\lambda^m_{j}} \right).
\end{equation}
}
%
%{\color{teal} END NEW}

	With a similar procedure, we can obtain lower and upper bounds for {\color{black} the MFPT for $\widecheck{X}$ from $\L_{u}$ to $\L_{\ell}$}. In particular, 
 for $\widecheck{\tau}_{\ell} = \inf\{ t \geq 0:\: \widecheck{X}(t) \in \L_{\ell} \} = \inf\{ t \geq 0:\: L(\widecheck{X}(t)) = \ell \}$ and $\widecheck{x}\in \L_u$, we have 
\begin{equation}
\label{upandlowbounds2}
    \invbreve{h}_{u,\ell}\le
     \E_{\widecheck{x}}[\widecheck{\tau}_{\ell}] \le 
     \Breve{h}_{u,\ell}
\end{equation}
 in which
%
 %{\color{teal} BEGIN OLD}
%	\begin{equation*}
%	\begin{aligned}
%	\invbreve{h}_{u,\ell}=\frac{r^m_{u-1}}{\gamma^M_{u}}\left(1+\sum_{j=\ell+1}^{u-1}\frac{1}{r^m_j}\right)+\sum_{i=\ell+2}^{u-1}\left[\frac{r^m_{i-1}}{\gamma^M_i}\left(1+\sum_{j=\ell+1}^{i-1}\frac{1}{r^m_j}\right)\right]+\frac{1}{\gamma^M_{\ell+1}}, 
%	\end{aligned}
%	\end{equation*}
%	\begin{equation*}
%	\begin{aligned}
%	\Breve{h}_{u,\ell}=\frac{r^M_{u-1}}{\gamma^m_{u}}\left(1+\sum_{j=\ell+1}^{u-1}\frac{1}{r^M_j}\right)+\sum_{i=\ell+2}^{u-1}\left[\frac{r^M_{i-1}}{\gamma^m_i}\left(1+\sum_{j=\ell+1}^{i-1}\frac{1}{r^M_j}\right)\right]+\frac{1}{\gamma^m_{\ell+1}}, 
%	\end{aligned}
%	\end{equation*}	
%	with $r^m_j=\frac{\lambda^m_{\ell+1} \lambda^m_{\ell+2}\dots\lambda^m_{j}}{\gamma^M_{\ell+1} \gamma^M_{\ell+2} \dots \gamma^M_{\ell}}$ and $r^M_j=\frac{\lambda^M_{\ell+1} \lambda^M_{\ell+2}\dots\lambda^M_{j}}{\gamma^m_{\ell+1} \gamma^m_{\ell+2} \dots \gamma^m_{\ell}}$.\\
%
%{\color{teal} END OLD}
%
%{\color{teal} BEGIN NEW}
%
{\color{black}
\begin{equation}
\label{formulaDto0lowbound}
\invbreve{h}_{u,\ell}=\frac{1}{\gamma^M_{\ell+1}} + \sum_{i=\ell+2}^{u} \frac{1}{\gamma^M_{i}} \left( 1 + \sum_{j=\ell+1}^{i-1} \frac{\lambda^m_{j}\dots\lambda^m_{i-1}}{\gamma^M_{j}\dots\gamma^M_{i-1}} \right),
\end{equation}	
\begin{equation}
\label{formulaDto0upbound}
\Breve{h}_{u,\ell}=\frac{1}{\gamma^m_{\ell+1}} + \sum_{i=\ell+2}^{u} \frac{1}{\gamma^m_{i}} \left( 1 + \sum_{j=\ell+1}^{i-1} \frac{\lambda^M_{j}\dots\lambda^M_{i-1}}{\gamma^m_{j}\dots\gamma^m_{i-1}} \right),
\end{equation}
provided that, for \eqref{formulaDto0lowbound}, $\gamma^M_z$ is positive for $\ell + 1 \leq z \leq u$, and for \eqref{formulaDto0upbound}, $\gamma^m_z$ is positive for $\ell + 1 \leq z \leq u$.\\
%Here, $\invbreve{h}_{u,\ell}$ in \eqref{formulaDto0lowbound} is well-defined provided that $\gamma^M_z$ is positive for $\ell + 1 \leq z \leq u$, and $\Breve{h}_{u,\ell}$ in \eqref{formulaDto0upbound} is well-defined provided that $\gamma^m_z$ is positive for $\ell + 1 \leq z \leq u$.\\
}
%
%{\color{teal} END NEW}

{\color{black}
\noindent \textbf{Application to Example \ref{ex:illustrativeexample}:} Consider the SCRN introduced in \eqref{reacs-3speciesexample}, with associated graph $\G$ shown in Figure \ref{fig:illexample}(b). We seek lower and upper bounds for the MFPT from $n_{\mathrm{Z}}=\Ntot$ to $n_{\mathrm{Y}}=\Ntot$. To this end, let us consider the coclique level function $L(x_1,x_2)=x_1+2x_2$ previously identified. The coclique level structure associated with it can be written as $\L_{\ell}, \dots, \L_u$, with $\L_z:=\{x\in \widecheck{\X}:\: L(x_1,x_2)=x_1 + 2x_2=z\}$ for $z=\ell,\dots,u$, and $\ell=0$, $u=2\Ntot$. This coclique level structure is such that $(0,0)^T$ (i.e., $n_{\mathrm{Z}}=\Ntot$) is the only state belonging to $\L_\ell$ and $(0,\Ntot)^T$ (i.e., $n_{\mathrm{Y}}=\Ntot$) is the only state belonging to $\L_u$ (Figure \ref{fig:illexample}(c)). This feature is critical in order to determine good lower and upper bounds for the MFPT from $n_{\mathrm{Z}}=\Ntot$ to $n_{\mathrm{Y}}=\Ntot$.

Now,
 \begin{equation}
\label{G+-ex00}
G_{+}=\{1,2\}\;\;\;\mathrm{and}\;\;\;G_{-}=\emptyset.
\end{equation}
%
%with $\widecheck{v}_1$ and $\widecheck{v}_2$ defined in \eqref{vectors3species}. 
The rate of increase $\lambda_z(\widecheck{x})$ and the rate of decrease $\gamma_z(\widecheck{x})$ can be written as 	
\begin{equation}\label{SimpleExampleuplowbounddrates}	\lambda_z(\widecheck{x})=f_1(\widecheck{x})+f_2(\widecheck{x})\;\;\;\mathrm{and}\;\;\;\gamma_z(\widecheck{x})=0, 
	\end{equation}
with $f_1(\widecheck{x})$, $f_2(\widecheck{x})$ defined in (\ref{rates3species}). 

The two continuous time Markov chains, $\Breve{X}$ and $\invbreve{X}$, are defined on the same state space as $\widecheck{X}$, with infinitesimal generators $\Breve{Q}$ and $\invbreve{Q}$, respectively, such that, for $z\in \{\ell, \ell +1 ,\dots, u-1,u\}$ and $\widecheck{x}\in \L_z$, $\Breve{Q}_{\widecheck{x},\widecheck{x}+\widecheck{v}_k} = \frac{\lambda^M_z}{|G_{+} (\widecheck{x})|}$ for $k \in G_{+} (\widecheck{x})$, 
$\invbreve{Q}_{\widecheck{x},\widecheck{x}+\widecheck{v}_k} = \frac{\lambda^m_z}{|G_{+} (\widecheck{x})|}$ for $k \in G_{+} (\widecheck{x})$, with $\lambda^{M}_z=\max_{\widecheck{x}\in \mathcal{L}_z}\lambda_z(\widecheck{x})$ and $\lambda^{m}_z=\min_{\widecheck{x}\in \mathcal{L}_z}\lambda_z(\widecheck{x})$, as defined in \eqref{minmaxrates},
and $G_{+} (\widecheck{x})$ defined as in \eqref{G+-(x)} where $G_{+}$ is given in \eqref{G+-ex00}. Note that in this example, 
$G_{-} =\emptyset$.
%$G_{-} (\widecheck{x})=\emptyset$ for all $\widecheck{x} \in \widecheck{\X}$.

 Then, as described in Section \ref{indeplevelstructure},
 %we can compare the Markov chain $\widecheck{X}$ with $\Breve X$ and $\invbreve{X}$, separately, by using Theorems 3.3 \textcolor{black}{and 3.4} in \cite{Monotonicitypaper} \textcolor{black}{(see SI - Section \ref{compthmpaper_thm})}, with the matrix $A=[1\; 2]$  associated with the \textcolor{black}{coclique} level function $L(x_1,x_2)=x_1+2x_2$. Then, the 
 we obtain analytical expressions for lower and upper bounds for the MFPT for $\widecheck{X}$ from $n_{\mathrm{Z}}=\Ntot$ to $n_{\mathrm{Y}}=\Ntot$ as
\begin{equation}\label{lowboundMFPTXSimpleExample}
	\begin{aligned}
		&\Breve{h}_{\ell,u}=\frac{1}{\lambda^{M}_{\ell}}+\frac{1}{\lambda^{M}_{\ell+1}}+\dots+\frac{1}{\lambda^{M}_{u-2}}+\frac{1}{\lambda^{M}_{u-1}},\;\; \invbreve{h}_{\ell,u}=\frac{1}{\lambda^{m}_{\ell}}+\frac{1}{\lambda^{m}_{\ell+1}}+\dots+\frac{1}{\lambda^{m}_{u-2}}+\frac{1}{\lambda^{m}_{u-1}},\nonumber
	\end{aligned}
\end{equation}
respectively. Given that $\lambda^{M}_z=\max_{\widecheck{x}\in \mathcal{L}_z} \left( \alpha(\Ntot -(x_1+x_2)) + \beta x_1\right)$ and $\lambda^{m}_z=\min_{\widecheck{x}\in \mathcal{L}_z} \left( \alpha(\Ntot -(x_1+x_2)) + \beta x_1\right)$, where the functions being maximized and minimized are increasing in $\alpha$, respectively $\beta$,
% {\color{purple}
% \begin{equation*}
%     \lambda^{M}_z=\max_{\widecheck{x}\in \mathcal{L}_z}\{ \alpha(\Ntot -(x_1+x_2)) + \beta x_1\}, %= \alpha(\Ntot - \frac{z}{2}) + \max_{\widecheck{x}\in \mathcal{L}_z}  (\beta -  \frac{\alpha}{2} )x_1 = \alpha(\Ntot - \frac{z}{2}) +   (\beta -  \frac{\alpha}{2} ) (z \wedge \Ntot) \one_{\{\beta > \frac{\alpha}{2}\}},
% \end{equation*}
% \begin{equation*}
%     \lambda^{m}_z=\min_{\widecheck{x}\in \mathcal{L}_z} \{ \alpha(\Ntot -(x_1+x_2)) + \beta x_1\}, %= \alpha(\Ntot - \frac{z}{2}) + \min_{\widecheck{x}\in \mathcal{L}_z}  (\beta -  \frac{\alpha}{2} )x_1,
% \end{equation*}
% }
%= \max_{\widecheck{x}\in \mathcal{L}_z} \alpha(\Ntot - \frac{z+x_1}{2}) + \beta x_1 
%$\lambda^{M}_z=\max_{\widecheck{x}\in \mathcal{L}_z}\lambda_z(\widecheck{x}) = \max_{\widecheck{x}\in \mathcal{L}_z}\{f_1(\widecheck{x})+f_2(\widecheck{x})\}$ and $\lambda^{m}_z=\min_{\widecheck{x}\in \mathcal{L}_z}\lambda_z(\widecheck{x})=\min_{\widecheck{x}\in \mathcal{L}_z}\{f_1(\widecheck{x})+f_2(\widecheck{x})\}$, 
we see that the lower and upper bounds $\Breve{h}_{\ell,u}$ and $\invbreve{h}_{\ell,u}$ decrease as either of the rate constants $\alpha$ or $\beta$ increases. 
%In this case, given the simplicity of the SCRN considered, we can also directly evaluate the MFPT from $n_{\mathrm{Z}}=\Ntot$ to $n_{\mathrm{Y}}=\Ntot$, which can be written as follows:
%
%\begin{equation}\label{MFPTexactformulasimpleexample}
%\begin{aligned}
%h_{Z=\Ntot,Y=\Ntot} &= \sum_{k=0}^{\Ntot - 1} \frac{1}{f_1(k, 0)} + \sum_{k=0}^{\Ntot - 1} \frac{1}{f_2(\Ntot - k, k)}\\
%&= \sum_{k=0}^{\Ntot - 1} \left(\frac{1}{\alpha(\Ntot - k)} + \frac{1}{\beta (\Ntot - k)}\right).
%\end{aligned}
%\end{equation}

%By inspecting the expression, it is possible to note that the dependence of the MFPT on the rate constants is consistent with what was observed through the study of its bounds.
%
For this example, due to the simplicity of the SCRN considered, it is also possible to derive a compact analytical expression for the MFPT. Each of the $\Ntot$ molecules acts independently and must undergo one $\ce{Z \to W}$ transition and one $\ce{W \to Y}$ transition to be converted from Z to Y. 
%
%(1) Each such reaction fires after an exponentially distributed waiting time, whose rate depends on the current state of the system. Let $n_{\mathrm{Z},k}$ and $n_{\mathrm{W},k}$ denote the number of Z and W molecules, respectively, at the time of the $k^{\text{th}}$ $\ce{Z \to W}$ and $\ce{W \to Y}$ transition. Then, the MFPT can be written as
%
%\begingroup
%\small
%\begin{equation}\label{MFPTexactformulasimpleexample}
%h_{Z=\Ntot,Y=\Ntot} = \sum_{k=1}^{\Ntot} \E\left[ \frac{1}{\alpha\, n_{\mathrm{Z},k}} \right] + \sum_{k=1}^{\Ntot} \E\left[ \frac{1}{\beta\, n_{\mathrm{W},k}} \right] = \frac{1}{\alpha}\sum_{k=1}^{\Ntot} \E\left[ \frac{1}{n_{\mathrm{Z},k}} \right] + \frac{1}{\beta} \sum_{k=1}^{\Ntot} \E\left[ \frac{1}{n_{\mathrm{W},k}} \right] ,
%\end{equation}
%\endgroup
%where the expectations are taken over the stochastic trajectory of the process. We can write the MFPT in this form because, in our case, the underlying continuous-time Markov chain is finite, irreducible on the transient states, and has a unique absorbing state corresponding to full conversion to species Y. Moreover, the process evolves through a sequence of exactly $2\Ntot$ transitions (i.e., $\Ntot$ $\ce{Z \to W}$ and $\Ntot$ $\ce{W \to Y}$ reactions), each governed by exponentially distributed waiting times with state-dependent rates. As a consequence, the MFPT corresponds to the expected total accumulated holding time before absorption, which can be expressed as the sum of expected sojourn times over all transition events \cite{Norris,Ross1996}. 
%
The time for one molecule of Z to be converted to Y is the sum of two independent exponential random variables with parameters $\alpha$ and $\beta$, respectively, which has a hypoexponential distribution with distribution function
\begin{equation*}
    F(x) = (F_1 * F_2) (x) = \int_0^x (1-e^{-\alpha(x-y)}) \beta e^{-\beta y} dy = \int_0^x (1-e^{-\beta(x-y)}) \alpha e^{-\alpha y} dy,
\end{equation*}
for $x>0$ where $F_1$ and $F_2$ are distribution functions for exponential random variables with parameters $\alpha$ and $\beta$, respectively.
% \begin{equation}
% \label{eqn:Fx}
% F(x) = 
% \begin{cases}
%     1 - \frac{\beta}{\beta-\alpha}e^{-\alpha x} - \frac{\alpha}{\alpha-\beta}e^{-\beta x} &\text{if } \alpha \neq \beta, \\
%     1-(1+\alpha x) e^{-\alpha x} &\text{if } \alpha = \beta.
% \end{cases}
% \end{equation}
% \begingroup
% \small
% \begin{equation}
% F(x) = 1 - \frac{\beta}{\alpha-\beta}e^{-\alpha x} - \frac{\alpha}{\alpha-\beta}e^{-\beta x}\;\;\mathrm{if}\;\; \alpha \neq \beta\;\;\mathrm{or}\;\;  \int_{0}^{x}\alpha^2ye^{-\alpha y} dy\;\; \mathrm{if}\;\; \alpha = \beta
% \end{equation}
% \endgroup
The MFPT we desire is the mean of the maximum of $\mathrm{N_{tot}}$ such independent hypoexponential random variables, which is given by $\int_{0}^{\infty}(1-F(x)^{\mathrm{N_{tot}}})dx$. Observe that for fixed $\beta>0$ and $0<y<x$, $1-e^{-\alpha(x-y)}$ is increasing with $\alpha$, and then $F(x)$ is increasing with $\alpha$. Thus, the MFPT decreases with increasing $\alpha$. Similarly, for fixed $\alpha>0$, we have that the MFPT decreases with increasing $\beta$.

Overall, this formulation explicitly gives the dependence of the MFPT on the rate constants $\alpha$ and $\beta$, and is fully consistent with the monotonic behavior that was observed through our study of its bounds.
}

{\color{black}
\subsection{Generalization to the non-weakly-connected case}
\label{generalization}

    In this section, we will consider SCRNs satisfying Assumption \ref{assumption:Unimolecular_change} and whose associated graphs do not necessarily have to be weakly connected. In this case, the graph $\G$ associated with a SCRN can be decomposed into finitely many weakly connected components, and we use $\wccNum$ to denote the number of weakly connected components in $\G$. For each $\wccIter = 1, \dots, \wccNum$, we use $\G^\wccIter$ to denote the $\wccIter^{th}$ weakly connected component of $\G$ and use $d_\wccIter$ to denote the number of vertices in $\G^\wccIter$. 

    \begin{example}
    \label{ex:Gen1}
        Consider a SCRN with four species ($\mathrm{S}_1$, $\mathrm{S}_2$, $\mathrm{S}_3$, $\mathrm{S}_4$) and the following three reactions under mass-action kinetics: 
        \begin{equation}
        \label{reacs-disconnected}
            {\large \textcircled{\small 1}}\;\ce{S_1 ->[$\kappa_1$] S_2}, \quad
            {\large \textcircled{\small 2}}\;\ce{S_2 ->[$\kappa_2$] S_1}, \quad
            {\large \textcircled{\small 3}}\;\ce{S_3 ->[$\kappa_3$] S_4},
        \end{equation}
        where the reaction rate constants $\kappa_1,\kappa_2,\kappa_3$ are positive. 
        %The associated reaction vectors are: $v_1 = (-1,1,0,0)^T, v_2 = (1,-1,0,0)^T, v_3 = (0,0,-1,1)^T$. 
        The associated graph $\G$ can be represented as follows:
  %      \begin{equation*}
   %     \text{(1)}\; 1 \rightleftarrows 2, \qquad \text{(2)}\; 3 \rightarrow 4.
    %    \end{equation*}
\begin{equation*}
\text{(1)}\; \mathrm{S_1} \xrightleftharpoons[e_2]{e_1} \mathrm{S_2}, \qquad \text{(2)}\; \mathrm{S_3} \xrightarrow{e_3} \mathrm{S_4}.
\end{equation*}
        In this case, $\G$ has two weakly connected components. The dynamics of $\{\mathrm{S}_1, \mathrm{S}_2\}$ evolve independently from those of $\{\mathrm{S}_3, \mathrm{S}_4\}$, and we can analyze the two subsystems independently. We will show, at the end of this section, how one can study the whole system by studying each subsystem independently.
    \end{example}

\begin{figure}[t]
\centering
\includegraphics[scale=0.39]{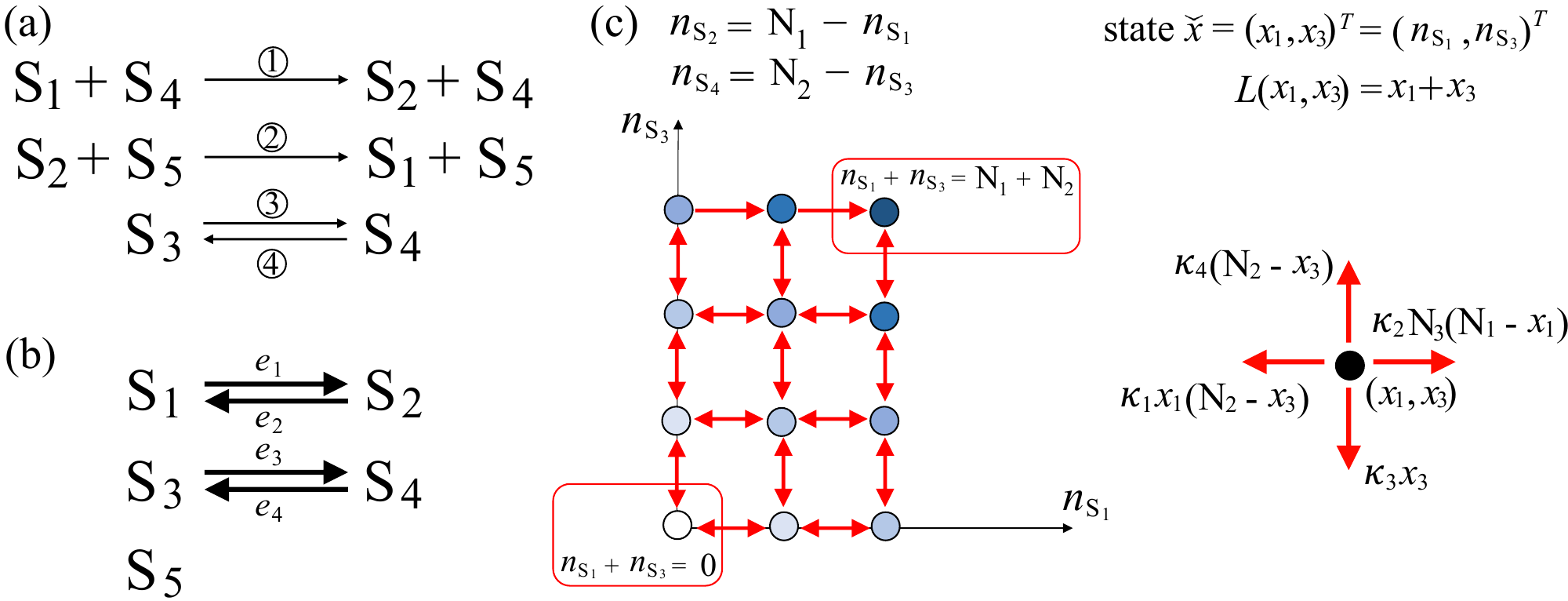}
\caption{\small {\color{black} { \bf Example of a SCRN whose associated graph $\G$ is not weakly connected (Example \ref{ex:Gen2}).}
(a) Chemical reaction system diagram.  The numbers on the arrows correspond to the reactions associated with the arrows as described in (\ref{reacs-example2}) in the main text. 
(b) Graph $\G$ associated with the chemical reaction system in panel (a).
(c) State space and transitions of the projected continuous time Markov chain $\widecheck{X}=\{
(X_1(t),X_3(t))^T:\: t \ge 0\}$, which keeps track of $(n_{\mathrm{S_1}},n_{\mathrm{S_3}})$ through time. Here, we consider $N_1=2$, $N_2=3$ and we use dots to represent the states, and red double-ended (single-ended) arrows to represent transitions in both directions (in a single direction). 
Additionally, we use shades of blue to distinguish the level to which each state belongs. The function $L$ associated with the coclique level structure is $L(x_1,x_3)=x_1+x_3$.
%
%(c)
%
On the right-hand side of the panel, we show the rates associated with the one-step transitions for the projected Markov chain $\widecheck{X}$. Note that, given an initial condition, the quantity of species $\mathrm{S}_5$ does not change over time and we denote this conserved quantity by $\mathrm{N_3}$.
            }}
   \label{fig:4.2example}
        \end{figure}    

    \begin{example}
    \label{ex:Gen2}
        For a reaction edge $\edge_k$ lying in some weakly connected component of $\G$, the infinitesimal rate $\Upsilon_k$ may depend on the quantities of some species where the graph vertices for those species may be in other weakly connected components of $\G$. For example, consider the following SCRN under mass-action kinetics:
        \begin{equation}
        \label{reacs-example2}
            \begin{aligned}
            & {\large \textcircled{\small 1}}\;\ce{S_1 + S_4 ->[$\kappa_1$] S_2 + S_4}, \quad
            {\large \textcircled{\small 3}}\;\ce{S_3 ->[$\kappa_3$] S_4}, \\
            & {\large \textcircled{\small 2}}\;\ce{S_2 + S_5 ->[$\kappa_2$] S_1 + S_5}, \quad {\large \textcircled{\small 4}}\;\ce{S_4 ->[$\kappa_4$] S_3},
            \end{aligned}
        \end{equation}
        where the reaction rate constants $\kappa_1, \kappa_2, \kappa_3, \kappa_4$ are positive. The diagram of this SCRN is shown in Figure \ref{fig:4.2example}(a) and the associated graph $\G$ is shown in Figure \ref{fig:4.2example}(b).
        %follows:
  %%      \begin{equation*}
   %%     \text{(1)}\; 1 \rightleftarrows 2, \qquad \text{(2)}\; 3 \rightarrow 4, \qquad \text{(3)}\; 5.
    %%    \end{equation*}
      %  \begin{equation*}
%\text{(1)}\; \mathrm{S_1} \xrightleftharpoons[e_2]{e_1} \mathrm{S_2}, \qquad 
%\text{(2)}\; \mathrm{S_3} \xrightleftharpoons[e_4]{e_3} \mathrm{S_4}, \qquad 
%\text{(3)}\; \mathrm{S_5}.
%\end{equation*}
        As a convention, we keep species $\mathrm{S}_5$ even though its quantity is conserved over time, and the vertex representing $\mathrm{S}_5$ in $\G$ is in its own weakly connected component.
        Here, $\G$ has three weakly connected components. 
        The dynamics of $\{\mathrm{S}_1, \mathrm{S}_2\}$ do depend on the varying quantity of $\mathrm{S}_4$, where $\mathrm{S}_4$ is in a different weakly connected component from that of $\{\mathrm{S}_1, \mathrm{S}_2\}$. In this case, we cannot analyze the first two subsystems independently. 
    \end{example}
    
    \medskip
    Now, for a SCRN satisfying Assumption \ref{assumption:Unimolecular_change} and whose associated graph $\G$ has $\wccNum$ weakly connected components, we shall and do relabel the species so that $j \in \{1,\dots,d \}$ belongs to the $\wccIter^{th}$ weakly connected component of $\G$ if and only if $\sum_{\wccIter'=1}^{\wccIter-1} d_{\wccIter'} < j \leq \sum_{\wccIter'=1}^{\wccIter} d_{\wccIter'}$ where a sum over an empty set is considered to be $0$. 
    With this relabeling, for a state $x \in \X$, $1 \leq \wccIter \leq \wccNum$ and $1 \leq i \leq d_\wccIter$, we shall use $x^\wccIter_i$ to denote the $\left(i+\sum_{\wccIter'=1}^{\wccIter-1} d_{\wccIter'} \right)^{th}$ entry of $x$. In other words, we have
    \begin{equation}
    \label{eqn:stateX_paper}
        x = (x^1, \dots, x^\wccNum)^T \in \Z^d_+ \text{ where } x^\wccIter = (x^\wccIter_1, \dots, x^\wccIter_{d_\wccIter})^T \in \Z^{d_\wccIter}_+ \text{ for } \wccIter = 1, \dots, \wccNum.
    \end{equation}
    For $\wccIter = 1, \dots, \wccNum$, each vertex in $\G^\wccIter$ corresponds to a species and each edge in $\G^\wccIter$ corresponds to a reaction, and in this sense, there is a stoichiometric matrix $S^\wccIter$ associated with $\G^\wccIter$.
    Since there is no edge between different weakly connected components of $\G$, the stoichiometric matrix $S$ associated with the SCRN has the form\footnote{\textcolor{black}{Here, $S^\wccIter$ is a $d_\wccIter \times r_\wccIter$ matrix, where $r_\wccIter$ is the number of edges in $\G^\wccIter$. A component $\G^\wccIter$ of $\G$ consisting of one vertex has no edges and the associated $S^\wccIter$ is degenerate.}}
    \begin{equation}
    \label{eqn:StoiMat}
        S = 
        \begin{bmatrix}
            S^1 & 0 & 0 \\
            0 & \ddots & 0 \\
            0 & 0 & S^\wccNum
        \end{bmatrix}.
    \end{equation}
    By applying Lemma \ref{lemmabis} to each weakly connected component of $\G$, we can show that the rank of the stoichiometric matrix $S$ is $d-\wccNum$ (which also implies that the chemical reaction network associated with $\G$, in the manner described in Remark \ref{rmk:deficiencyzero}, has deficiency zero) and there are $\wccNum$ linearly independent conservation vectors for this SCRN. 
    As a result, we shall consider a \textbf{projected continuous time Markov chain} $\widecheck{X}=\{\widecheck{X}(t):\: t \ge 0\}$ in which the state $\widecheck x \in \widecheck{\mathcal{X}} \subset \Z_+^{d-\wccNum}$ tracks, for each $\wccIter=1,\dots,\wccNum$, the number of molecules of each species in $\G^\wccIter$, except for the last species in $\G^\wccIter$, that is, 
    \begin{equation}
    \label{eqn:projectX_paper}
        \widecheck x = (\widecheck x^1, \dots, \widecheck x^\wccNum)^T \text{ where } \widecheck x^\wccIter = (x^\wccIter_1, \dots, x^\wccIter_{d_\wccIter-1})^T \text{ for } \wccIter = 1, \dots, \wccNum.
    \end{equation}
    Here, we slightly abuse notation, where $\widecheck{x}^\wccIter$ is a zero-dimensional vector if $|\G^\wccIter| = 1$. Note that the choice to express $x^\wccIter_{d_\wccIter}$ as a function of $x^\wccIter_i : 1 \leq i \leq d_\wccIter-1$ is without loss of generality, since the species can always be relabeled so that a species chosen to be expressed as a function of the others is the last one. 
    The process $\widecheck{X}$ is a continuous time Markov chain\footnote{\textcolor{black}{With this reduction, coordinates corresponding to graph components consisting of a single vertex are removed, although the transition rates of the Markov chain $\widecheck{X}$ may depend on the constant amount of any species associated with such an isolated vertex.}} defined on the finite state space 
    \begin{equation*}
        \begin{aligned}
            \widecheck{\mathcal{X}}&\vcentcolon=\left\{\widecheck x = ( x^\wccIter_i)_{1 \leq i \leq d_\wccIter - 1, 1 \leq \wccIter \leq \wccNum} \in \Z_{+}^{d-\wccNum}: \phantom{\sum_{i=1}^{d_\wccIter-1}} \right.\\
            & \qquad \qquad \quad \left. ( x^\wccIter_i)_{1 \leq i \leq d_\wccIter, 1 \leq \wccIter \leq \wccNum} \in \X \text{ where }  x^\wccIter_{d_\wccIter} = x^\wccIter_{\mathrm{tot}} - \sum_{i=1}^{d_\wccIter-1}x^\wccIter_i \text{ for } \wccIter = 1, \dots, \wccNum \right\}\\
            &\;\subset\left\{\widecheck x = ( x^\wccIter_i)_{1 \leq i \leq d_\wccIter - 1, 1 \leq \wccIter \leq \wccNum} \in \Z_{+}^{d-\wccNum}: x^\wccIter_1+\dots+x^\wccIter_{d_\wccIter-1} \le x^\wccIter_{\mathrm{tot}} \text{ for } \wccIter = 1, \dots, \wccNum \right\},   
        \end{aligned}
    \end{equation*}
    where $x^\wccIter_{\mathrm{tot}}=\sum_{i=1}^{d_\wccIter}X^\wccIter_i(0)$ for each $\wccIter = 1, \dots, \wccNum$.
    We will assume that $|\widecheck{\X}| > 1$, and the infinitesimal generator of $\widecheck{X}$ will be denoted by $\widecheck Q$. 

    A \textbf{coclique level function} for $\widecheck{X}$ is a linear function $L: \Z^{d-\wccNum} \to \Z$ such that for each $k= 1,\dots,n$,
    \begin{equation}
    \label{indlevfunctiondef2_paper}
        L(\widecheck{v}_k) \in \{-1,+1\},
    \end{equation}
    where $\widecheck{v}_k\in \Z^{d-\wccNum}$ is the vector obtained from $v_k$ by removing the last element in each of the weakly connected components in $\G$, in a similar manner to how we obtain \eqref{eqn:projectX_paper} from \eqref{eqn:stateX_paper}.
    If such an $L$ exists, it can be written as 
    \begin{equation}
    \label{ILFformula2_paper}
        L(x) = b^T x\;\; \mathrm{for}\; x \in \Z^{d-\wccNum}\; \mathrm{and\; some}\; b \in \Z^{d-\wccNum},
    \end{equation}
    where, upon partitioning the set of edges of the associated graph $\G$ into two disjoint subsets $\Eu=\{\edge_k:L(\widecheck{v}_k)=1\}$ and $\Ed=\{\edge_k:L(\widecheck{v}_k)=-1\}$ (where one of these may be empty), the vector $b=(b_1,\dots,b_{d-\wccNum})^T$ solves the system of equations
    \begin{equation}
    \label{syst2_paper}
	\sum_{i=1}^{d-\wccNum} b_{i} \widecheck{v}_k(i)=
        \begin{cases}
    	+1\;\;\;\text{if }\; \edge_k \in \Eu\\
    	-1\;\;\;\text{if }\; \edge_k \in \Ed\\
	\end{cases} \;\;\;\;\; \text{for } k=1,\dots,n.
    \end{equation}
    For a coclique level function $L$, a \textbf{coclique level structure} for $\widecheck{X}$ can be defined similar to that in the case where $\G$ is weakly connected.

    Now, we are going to characterize all coclique level structures for $\widecheck{X}$. In the following theorem, if the $\wccIter^{th}$ weakly connected component $\G^\wccIter$ of $\G$ has more than one vertex, i.e., $|\G^\wccIter| > 1$, then we will consider a SCRN that consists of all the species and reactions associated with this weakly connected component.
    If $|\G^\wccIter| = 1$, then $\widecheck{x}^\wccIter$ is a zero-dimensional vector and we do not need to consider this weakly connected component when constructing coclique level functions as in \eqref{eqn:genCLF}.
    
    \begin{theorem}
    \label{thm:Gen}
        Consider a SCRN satisfying Assumption \ref{assumption:Unimolecular_change} and $|\widecheck{\X}| > 1$.
        A coclique level function for the projected continuous time Markov chain $\widecheck X$ exists if and only if for every weakly connected component $\G^\wccIter$ of the associated graph $\G$ satisfying $|\G^\wccIter|>1$ where $\wccIter \in \{1,\dots,\wccNum\}$, an associated coclique level function exists.
        The set of all coclique level functions for $\widecheck{X}$ is the set of all functions of the form
        \begin{equation}
        \label{eqn:genCLF}
            L(\widecheck x) = L (\widecheck x^1, \dots, \widecheck x^\wccNum) = \sum_{\substack{\wccIter=1,\dots,\wccNum: \\ |\G^\wccIter|>1}} L^\wccIter (\widecheck x^\wccIter)
        \end{equation}
        where $L^\wccIter: \Z^{d_\wccIter-1} \rightarrow \Z$ is a coclique level function associated with $\G^\wccIter$ for $\wccIter=1,\dots,\wccNum$ and $|\G^\wccIter| > 1$.
        In particular, a coclique level function for $\widecheck{X}$ exists if and only if $\G$ is bipartite.
    \end{theorem}
    
    The proof of Theorem \ref{thm:Gen} can be found in SI - Section \ref{sec:thmGenProof}. For this, we can break down the problem to look at each weakly connected component
    separately, which is made possible by the block structure of the stoichiometric matrix shown in \eqref{eqn:StoiMat}. 
    
    Suppose the projected continuous time Markov chain $\widecheck{X}$ with finite state space $\widecheck{\X} \subseteq \Z^{d-\wccNum}_{+}$ has a coclique level function $L: \Z^{d-\wccNum} \rightarrow \Z$, defined as in \eqref{indlevfunctiondef2_paper}, with coclique level structure $\L_{\ell},\L_{\ell+1}, \dots, \L_{u-1}, \L_u$ defined as in \eqref{indlevstructure}--\eqref{lowhighletters} with $\ell < u$.
    Similar to the reasoning in Remark \ref{correlationpartitionstructure}, for a coclique level function $L$ as in \eqref{eqn:genCLF}, the sets $\L_{\ell}, \dots, \L_u$ in the coclique level structures associated with $L$ and $-L$ are the same, while the orderings of those sets in the two coclique level structure partitions are opposite. As a convention, we consider the coclique level structures associated with $L$ and $-L$ to be the same.
    Similar to the discussion in Section \ref{indeplevelstructure}, we can determine analytical expressions for upper and lower bounds for both the MFPT for $\widecheck{X}$ from $\L_{\ell}$ to $\L_{u}$ and the MFPT for $\widecheck{X}$ from $\L_{u}$ to $\L_{\ell}$. 

    \medskip

    \noindent \textbf{Application to Example \ref{ex:Gen1}:} In this example, $\G$ has two weakly connected components, and $m_1 = (1,1,0,0)^T$ and $m_2 = (0,0,1,1)^T$ are two linearly independent conservation vectors for this system. Therefore, we can consider a projected continuous time Markov chain $\widecheck{X}=\{({X}_1(t), {X}_3(t))^T : t \geq 0\}$, which tracks the number of $\mathrm{S}_1$ and $\mathrm{S}_3$ through time. The state space for $\widecheck{X}$ is $\{ (x_1,x_3)^T \in \Z^2_+: 0 \leq x_1 \leq \mathrm{N_1}, 0 \leq x_3 \leq \mathrm{N_2} \}$, where $\mathrm{N_1} = X_1(0) + X_2(0)$ and $\mathrm{N_2} = X_3(0) + X_4(0)$. We shall use Theorem \ref{thm:Gen} to identify all of the coclique level structures for $\widecheck{X}$. For the first weakly connected component $\G^1$ of $\G$, the only possible useful partitions of edges of $\G^1$ are $\Eu^1=\{\edge_2\},\Ed^1=\{\edge_1\}$ and $\Eu^1=\{\edge_1\},\Ed^1=\{\edge_2\}$, by Remark \ref{rmk:reversible} and the fact that reactions ${\large \textcircled{\small 1}}$ and ${\large \textcircled{\small 2}}$ are a pair of reversible reactions. One then can verify that $L^1 (\widecheck{x}^1) = x_1^1 = x_1$ and $L^1 (\widecheck{x}^1) = - x_1^1 = - x_1$ are the two coclique level functions associated with $\G^1$. Similarly, $L^2 (\widecheck{x}^2) = x_1^2 = x_3$ and $L^2 (\widecheck{x}^2) = - x_1^2 = - x_3$ are the only coclique level functions associated with $\G^2$. By Theorem \ref{thm:Gen}, the coclique level functions for $\widecheck{X}$ are $L(\widecheck{x})=x_1+x_3$, $L(\widecheck{x})=x_1-x_3$, $L(\widecheck{x})=-x_1+x_3$ and $L(\widecheck{x})=-x_1-x_3$, where the first two give the two distinct coclique level structures for $\widecheck{X}$ by similar reasoning to that in Remark \ref{correlationpartitionstructure}.

     %We can determine lower and upper bounds for the MFPT from $\widecheck{x}=(x_1,x_3)^T=(0,0)^T$ to $(\mathrm{N_1},\mathrm{N_2})^T$ by exploiting the coclique level structure associated with the coclique level function $L(\widecheck{x})=x_1+x_3$.

     Because the dynamics of $\{\mathrm{S}_1, \mathrm{S}_2\}$ evolve independently from those of $\{\mathrm{S}_3, \mathrm{S}_4\}$, we can choose to analyze the two subsystems independently. For example, $\widecheck{X}^1=\{X_1(t) : t \geq 0\}$ is itself a continuous time Markov chain that tracks the number of $\mathrm{S}_1$ and we have shown that $L^1 (\widecheck{x}^1) = x_1^1 = x_1$ is a coclique level function associated with $\G^1$. Thus, we may study the MFPT from $\widecheck{X}^1 = 0$ to $\widecheck{X}^1 = \mathrm{N_1}$. However, it is not possible to study this MFPT using a coclique function for the whole system, because $L(\widecheck{x})=x_1$ is not a coclique level function for $\widecheck{X}$. In conclusion, when the dynamics of each subsystem in the SCRN evolve independently, it is beneficial to apply our theory in Sections \ref{suffcond}--\ref{indeplevelstructure} to each subsystem separately.
     
    \noindent \textbf{Application to Example \ref{ex:Gen2}:} In this example, $\G$ has three weakly connected components, and $m_1 = (1,1,0,0,0)^T$, $m_2 = (0,0,1,1,0)^T$ and $m_3=(0,0,0,0,1)^T$ are three linearly independent conservation vectors for this system.  Therefore, we can consider a projected continuous time Markov chain $\widecheck{X}=\{(X_1(t), X_3(t))^T : t \geq 0\}$, which tracks the number of $\mathrm{S}_1$ and $\mathrm{S}_3$ through time (Figure \ref{fig:4.2example}(c)). The state space for $\widecheck{X}$ is $\{ (x_1,x_3)^T \in \Z^2_+: 0 \leq x_1 \leq \mathrm{N_1}, 0 \leq x_3 \leq \mathrm{N_2} \}$, where $\mathrm{N_1} = X_1(0) + X_2(0)$ and $\mathrm{N_2} = X_3(0) + X_4(0)$. Note that, given an initial condition, the quantity of species $\mathrm{S}_5$ does not change over time and we shall denote this conserved quantity by $\mathrm{N_3}$. Accordingly, the trivial dynamics of the system associated with the third weakly connected component of $\G$ is omitted in the projected chain $\widecheck{X}$. Similar to the previous example, we can use Theorem \ref{thm:Gen} to identify all of the coclique level structures for $\widecheck{X}$, which are the coclique level structures associated with the coclique level functions $L(\widecheck{x})=x_1+x_3$ and $L(\widecheck{x})=x_1-x_3$. 
    %This allows us to determine lower and upper bounds for the MFPT from $\widecheck{x}=(x_1,x_3)^T=(0,0)^T$ to $(\mathrm{N_1},\mathrm{N_2})^T$ and the MFPT from $(\mathrm{N_1},\mathrm{N_2})^T$ to $(0,0)^T$, by exploiting the coclique level structure associated with the coclique level function $L(\widecheck{x})=x_1+x_3$. We can also determine lower and upper bounds for the MFPT from $(0,\mathrm{N_2})^T$ to $(\mathrm{N_1},0)^T$ and the MFPT from $(\mathrm{N_1},0)^T$ to $(0,\mathrm{N_2})^T$, by exploiting the coclique level structure associated with the coclique level function $L(\widecheck{x})=x_1-x_3$. 
    
    As an example, we can use the coclique level structure associated with the coclique level function $L(\widecheck{x})=x_1+x_3$ to bound the MFPT from $\widecheck{x}=(x_1,x_3)^T=(\mathrm{N_1},\mathrm{N_2})^T$ to $(0,0)^T$ above and below by \eqref{formulaDto0lowbound} and \eqref{formulaDto0upbound}, respectively, where the parameters in these two expressions are given by $\ell=0$, $ u= \mathrm{N_1}+\mathrm{N_2}$ and for $\ell \leq z \leq u$,
    \begin{eqnarray*}
        \lambda^m_z &=& \min_{\substack{0 \leq x_1 \leq \mathrm{N_1}, \\0 \leq x_3 \leq \mathrm{N_2},\\x_1+x_3=z}} \left( \kappa_2 (\mathrm{N_1} - x_1) \mathrm{N_3} + \kappa_4 (\mathrm{N_2} - x_3) \right), \\
        \lambda^M_z &=& \max_{\substack{0 \leq x_1 \leq \mathrm{N_1}, \\0 \leq x_3 \leq \mathrm{N_2},\\x_1+x_3=z}} \left( \kappa_2 (\mathrm{N_1} - x_1) \mathrm{N_3} + \kappa_4 (\mathrm{N_2} - x_3) \right), \\ 
        \gamma^m_z &=& \min_{\substack{0 \leq x_1 \leq \mathrm{N_1}, \\0 \leq x_3 \leq \mathrm{N_2},\\x_1+x_3=z}} \left( \kappa_1 x_1 (\mathrm{N_2} - x_3) + \kappa_3 x_3 \right), \quad \gamma^M_z = \max_{\substack{0 \leq x_1 \leq \mathrm{N_1}, \\0 \leq x_3 \leq \mathrm{N_2},\\x_1+x_3=z}} \left( \kappa_1 x_1 (\mathrm{N_2} - x_3) + \kappa_3 x_3 \right),
    \end{eqnarray*}
    in which one can verify that $\gamma^M_z \geq \gamma^m_z > 0$ for each $1 \leq z \leq \mathrm{N_1}+\mathrm{N_2}$.
    
    %{\color{teal} TO ADD?? for Example \ref{ex:Gen1} (if we look at the choice of independent level function for (1)+(2) then it is $L(x)=x_1 \pm x_3$; if we look at them separately, then we have $L^1(x)=x_1$ and $L^2(x)=x_3$. In the latter case, we can study the MFPT from no $\mathrm{S}_1$ to all $\mathrm{S}_1$, but in the first case, can only study the MFPT for example, from ``no $\mathrm{S}_1$ and no $\mathrm{S}_3$" to ``all $\mathrm{S}_1$ and all $\mathrm{S}_3$", from ``no $\mathrm{S}_1$ and all $\mathrm{S}_3$" to ``and $\mathrm{S}_1$ and no $\mathrm{S}_3$")... for Example \ref{ex:Gen2} ...}
}

%%%%%%%%%%%%%%%%%%%%%%%%%%%%%%%%%%%%%%%%%%%%%%%%%%%%%

\section{Examples}
%\section{Four Models For Chromatin Modification Circuits}
\label{sec:Applications}

In this section, we consider three examples and study their stochastic behavior, in terms of MFPT, by exploiting the theoretical tools developed in this paper. The examples come from biological areas of epigenetics, neurobiology and ecology \cite{bib:BWD2022,RMilo2002,FBrglez1989}. All the continuous time Markov chains associated with the SCRNs considered in our examples have a finite state space. 
%Furthermore, for all the SCRNs considered,  it is straightforward to verify that Assumption \ref{assumption:Unimolecular_change} is satisfied.
%
%The choice of matrix $A$ in each example is based on the specific monotonicity relationship of interest. While for simpler cases the choice of $A$ is straightforward, for more complicated systems the choice can be more subtle. In many cases, in order to study the monotonicity properties for the stochastic behavior of our system, we can rely on Theorem \ref{thm:InnerProductTheorem}, which provides a reasonable approach to narrow down the choices for suitable $A$. The approach consists in solving, for each row $i$, the system of equations $\sum_{k=1}^d A_{ik}(v_j)_k=b_{ij}$, with $b_{ij}$ equal to $1,-1,$ or $0$ depending, based on the monotonicity relationship of interest, whether we expect that the Markov chain transition in the direction $v_j$ leads inside, outside, or is parallel to the boundary of the region $K_A+x$. Finally, it is worth noticing that, while all the following examples compare two identical reaction networks with different rate constants, our theory can also be applied to compare two different reaction networks as long as they have the same net reaction vectors $\{v_j\}_{j=1}^n$.
%
%

	\subsection{Chromatin modification circuit including only histone modifications}\label{exp2}

\begin{figure}[t]
\centering
\includegraphics[scale=0.39]{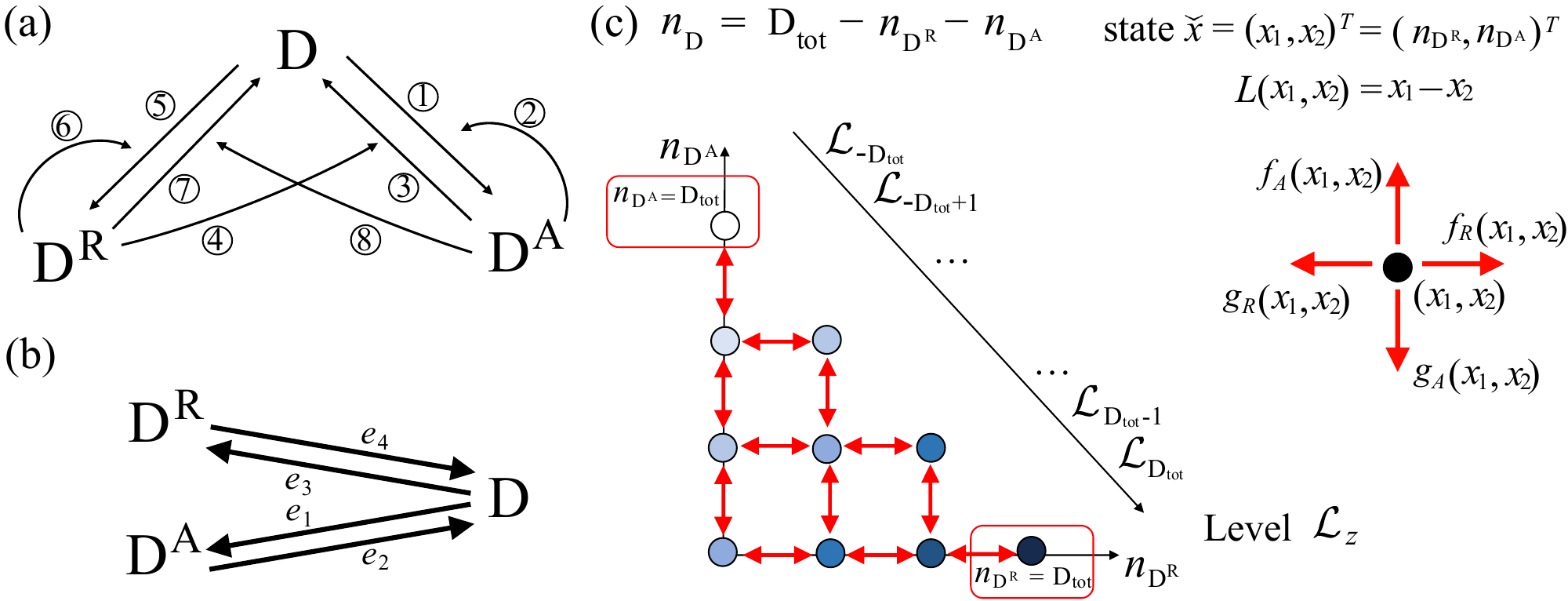}
\caption{\small { \bf Chromatin modification circuit including only histone modifications: reaction diagram, graph $\G$ and associated Markov chain.}
(a) Chemical reaction system diagram.  The numbers on the arrows correspond to the reactions associated with the arrows as described in (\ref{reacs2D}) in the main text. 
(b) Graph $\G$ associated with the chemical reaction system in panel (a).
(c) State space and transitions of the projected continuous time Markov chain $\widecheck{X}=\{
(X_1(t),X_2(t))^T:\: t \ge 0\}$, which keeps track of $(n_{\mathrm{D^R}},n_{\mathrm{D^A}})$ through time. Here, we consider $\mathrm{D_{tot}}=3$ and we use dots to represent the states, and red double-ended arrows to represent transitions in both directions. 
Additionally, we use shades of blue to distinguish the level to which each state belongs. The function $L$ associated with the coclique level structure is $L(x_1,x_2)=x_1-x_2$.
%
%(c)
%
The rates associated with the one-step transitions for the projected Markov chain $\widecheck{X}$ are given in (\ref{rates2D}).
            }
   \label{fig:2Dmodel}
        \end{figure}

	Epigenetic regulation is the modification of the DNA structure, due to chromatin modifications, that determines if a gene is active or repressed. Various chromatin modifications affect the structure of DNA. In this example, we will consider only histone modifications, while in the next one, we will study a more complex model including also DNA methylation. More precisely, in this example we analyze a well-established model for a histone modification circuit \cite{2007DoddCellPaper,bib:BWD2022,ECC2022}, which involves three species: unmodified nucleosome, denoted by D; nucleosome modified with repressive histone modifications, denoted by $\mathrm{D^R}$; nucleosome modified with activating histone modifications, denoted by $\mathrm{D^A}$. In this model, each histone modification catalyzes its own establishment on unmodified nucleosomes and catalyzes the erasure of the opposite modification \cite{2007DoddCellPaper,bib:BWD2022}. 
 %
 %The graph $\G$ associated with the chemical reaction system considered is depicted in Fig. \ref{fig:2Dmodel}(a). 
 The quantity of each species is denoted by $n_{\mathrm{D}}$, $n_{\mathrm{D^R}}$ and $n_{\mathrm{D^A}}$, respectively. Their sum remains constant, i.e., $n_{\mathrm{D}}+n_{\mathrm{D^R}}+n_{\mathrm{D^A}}=\Dtot$, where $\Dtot$ denotes the total number of nucleosomes within the gene. Then, the chemical reaction system, whose diagram is shown in Figure \ref{fig:2Dmodel}(a), can be written as
    
    \begingroup
    \small
	\begin{equation}\label{reacs2D}
	\begin{aligned}
	&{\large \textcircled{\small 1}}\;\ce{D ->[$k^A_{W0}+k^A_W$] D^A },\;\;{\large \textcircled{\small 2}}\;\ce{D +D^A ->[$k^A_M$] D^A + D^A},\;\;{\large \textcircled{\small 3}}\;\ce{D^A ->[$\delta +\bar k^A_E$] D},\\
 &{\large \textcircled{\small 4}}\;\ce{D^A + D^{R} ->[$k^A_E$] D + D^{R}},\;\;{\large \textcircled{\small 5}}\;\ce{D ->[$k^R_{W0}+k^R_{W}$] D^R },\;\;{\large \textcircled{\small 6}}\;\ce{D +D^R ->[$k^R_M$] D^R + D^R},\\
 &{\large \textcircled{\small 7}}\;\ce{D^R ->[$\delta+\bar k^R_E$] D},\;\;{\large \textcircled{\small 8}}\;\ce{D^R + D^{A} ->[$k^R_E$] D + D^{A}},\\
	\end{aligned}
	\end{equation}
	\endgroup
    in which $k^A_{W0},k^A_W,k^A_M,\delta, \bar k^A_E,  k^A_E, k^R_{W0}, k^R_{W}, k^R_M, \bar k^R_E, k^R_E > 0$. The expression of the reaction rate constants is because we combined reactions sharing the same reactants and products. Now, denoting the reaction volume by $V$, let us introduce $\eps := \frac{\delta + \bar k_E^A}{k_M^A(\Dtot/V)} = \frac{\delta_A}{k_M^A(\Dtot/V)}$, with $\delta_A:=\delta + \bar k_E^A$, and $\mu := \frac{k^R_E}{k^A_E}$. Furthermore, let us consider the constant $\tilde{b}$ such that $\mu \tilde{b} = \frac{\delta_{R}}{\delta_{A}}$, where $\delta_R:=\delta + \bar k_E^R$. Then, $\delta_R = \delta_A \mu \tilde{b} = \eps \frac{k_M^A\Dtot}{V} \mu \tilde{b}$. 
    %So, as $\eps \to 0$, both $\delta_A$ and $\delta_R$ go to $0$ with $\Dtot, \frac{k_M^A\Dtot}{V}, \mu$ and $b$ fixed. 

Considering $x=(n_{\mathrm{D^R}},n_{\mathrm{D^A}},n_{\mathrm{D}})^T$, the reaction vectors associated with \eqref{reacs2D} are $v_1=(0,1,-1)^T$, $v_2=(0,-1,1)^T$, $v_3=(1,0,-1)^T$, and $v_4=(-1,0,1)^T$. By examining them, it is possible to verify that Assumption \ref{assumption:Unimolecular_change} is satisfied. The graph $\G$ associated with the chemical reaction system \eqref{reacs2D} can then be represented as in Figure \ref{fig:2Dmodel}(b).

By inspecting $\G$, one can verify that the underlying undirected graph is connected and that $\G$ is bipartite. 
%Given that Assumption \ref{assumption:Unimolecular_change} is satisfied, we can then apply Lemma \ref{lemmabis}. This implies 
By Lemma \ref{lemmabis}, our SCRN has a unique conservation vector $m = (1,1,1)^T$ and then we can introduce a projected continuous time Markov chain $\widecheck{X}=\{
(X_1(t),X_2(t))^T:\: t \ge 0\}$, which keeps track of $(n_{\mathrm{D^R}},n_{\mathrm{D^A}})$ through time.
Since the total number of nucleosomes $\Dtot$ is conserved, the state space is $\widecheck{\X}= \{\widecheck{x}=(x_1,x_2)^T \in \Z_+^2 :\: x_1 + x_2 \leq \Dtot \}$. The potential one-step transitions for $\widecheck X$ from $ \widecheck x \in \widecheck{\X}$ are shown in Figure \ref{fig:2Dmodel}(c), where the associated transition vectors are  $\widecheck{v}_1=-\widecheck{v}_2=(0,1)^T$ and $\widecheck{v}_3=-\widecheck{v}_4=(1,0)^T$, and the infinitesimal transition rates (in which we assume mass-action kinetics with the usual rate constant volume scaling) are
   \begin{equation}\label{rates2D}
    \begin{aligned}
&\widecheck{Q}_{\widecheck{x},\widecheck{x}+\widecheck{v}_1}=f_A(\widecheck{x}) = (\Dtot -(x_1+x_2))\left(k_{W0}^A+k_{W}^A + \frac{k_{M}^A}{V}x_2\right),\\        &\widecheck{Q}_{\widecheck{x},\widecheck{x}+\widecheck{v}_2}=g_A(\widecheck{x}) = x_2\left(\eps \frac{k_{M}^A}{V}\Dtot + x_1\frac{k^A_E}{V}\right),\\
&\widecheck{Q}_{\widecheck{x},\widecheck{x}+\widecheck{v}_3}=f_R(\widecheck{x}) = (\Dtot -(x_1+x_2))\left(k_{W0}^R+k_{W}^R + \frac{k_{M}^R}{V}x_1\right),\\        &\widecheck{Q}_{\widecheck{x},\widecheck{x}+\widecheck{v}_4}=g_R(\widecheck{x}) = x_1\mu\left(\eps \frac{k_{M}^A}{V}\Dtot \tilde{b} + x_2\frac{k^A_E}{V}\right).
       % g^{\eps}_R(x_1,x_2) &= i\left(\eps \frac{k_{M}^A}{V}\Dtot\mu b + j\frac{k^R_E}{V}\right)
    \end{aligned}
   \end{equation}
   
Let us now focus on determining explicit analytical expressions for upper and lower bounds of MFPTs. To this end, let us apply Theorem \ref{thm2b} to determine the \textcolor{black}{coclique} level structures for $\widecheck{X}$. We can apply the theorem because Assumption \ref{assumption:Unimolecular_change} is satisfied
%, the species of our SCRN are $d=3>2$, 
and the associated graph $\G$ is weakly connected. %and  the rank of the stoichiometric matrix $S$ is $2=d-1$.

Consider all of the possible partitions $\{ \Eu, \Ed \}$ of edges of $\G$ that could allow us to determine \textcolor{black}{a coclique} level structure. These partitions are the following:
\begin{align}
    &\Eu=\{\edge_1,\edge_3\},\Ed=\{\edge_2,\edge_4\}\;\;\;\;\mathrm{and}\;\;\;\;\textcolor{black}{\Eu=\{\edge_2,\edge_3\},\Ed=\{\edge_1,\edge_4\}}.
\end{align}
We did not consider the partition $\Eu=\{\edge_2,\edge_4\},$ $\Ed=\{\edge_1,\edge_3\}$ or the partition \textcolor{black}{$\Eu=\{\edge_1,\edge_4\},$ $\Ed=\{\edge_2,\edge_3\}$} because, as explained in Remark \ref{correlationpartitionstructure}, the associated functions $L$ would be the opposite of the ones obtained for the partitions considered above and the resulting \textcolor{black}{coclique} level structures are considered to be the same.
Furthermore, we did not consider the partitions in which two edges in the same element of a partition are associated with two reaction vectors $v_k, v_{k'} \in \V$ such that $v_k = -v_{k'}$, because, as stated in Corollary \ref{thm2bis}, these partitions would not lead to \textcolor{black}{a coclique} level structure. For each partition, the system of equations in \eqref{syst}
%, obtaining
%\begin{align}
%b_1 =1,\;-b_1=-1,\;b_2=1,\;-b_2=-1
%\end{align}
%%
%and
%%
%\begin{align}
%b_1 =1,\;-b_1=-1,\;b_2=-1,\;-b_2=1,
%\end{align}
%%
%respectively. Both systems admit a unique solution $b \in \Z ^2$,
has a unique solution, these being $(b_1,b_2)^T=(1,1)^T$ and $(b_1,b_2)^T=(1,-1)^T$, respectively. Then, by applying Theorem \ref{thm2b}, we can conclude that the projected Markov chain $\widecheck{X}$ has two \textcolor{black}{coclique} level structures associated with \textcolor{black}{coclique} level functions $L(x_1,x_2)=x_1+x_2$ and $L(x_1,x_2)=x_1-x_2$.

Let us consider the \textcolor{black}{coclique} level function $L(x_1,x_2)=x_1-x_2$. The \textcolor{black}{coclique} level structure associated with it can be written as $\L_{\ell}, \dots, \L_u$, with $\L_z:=\{x\in \widecheck{\X}:\: L(x_1,x_2)=x_1 - x_2=z\}$ for $z=\ell,\dots,u$, with $\ell=-\mathrm{D_{tot}}$ and $u=\mathrm{D_{tot}}$. This \textcolor{black}{coclique} level structure is such that $a=(0,\mathrm{D_{tot}})^T$ is the only state belonging to $\L_\ell$ and $r=(\mathrm{D_{tot}},0)^T$ is the only state belonging to $\L_u$ (see Figure \ref{fig:2Dmodel}(c)). As shown below, this feature of the \textcolor{black}{coclique} level structure is critical in order to determine lower and upper bounds for $h_{a,r}$, the MFPT from $a$ to $r$, and $h_{r,a}$, the MFPT from $r$ to $a$, and this is the reason why we consider the \textcolor{black}{coclique} level structure associated with the function $L(x_1,x_2)=x_1-x_2$ and not the one associated with the function $L(x_1,x_2)=x_1+x_2$. 

Let us now determine the lower and upper bounds for the MFPT from $a$ to $r$, and vice versa, following the approach described in the previous section. In particular, here we have 
%$G_{+}=\{\widecheck{v}_2,\widecheck{v}_3\}$ and $G_{-}=\{\widecheck{v}_1,\widecheck{v}_4\}$,
\begin{equation}
\label{G+-ex1}
G_{+}=\{2,3\}\;\;\;\mathrm{and}\;\;\;G_{-}=\{1,4\},
\end{equation}
and the rate of increase $\lambda_z(\widecheck{x})$ and the rate of decrease $\gamma_z(\widecheck{x})$ can be written as 	
\begin{equation}\label{2Duplowbounddrates}	\lambda_z(\widecheck{x})=f_R(\widecheck{x})+g_A(\widecheck{x})\;\;\;\mathrm{and}\;\;\;\gamma_z(\widecheck{x})=f_A(\widecheck{x})+g_R(\widecheck{x}), 
	\end{equation}
with $f_R(\widecheck{x})$, $g_A(\widecheck{x})$, $f_A(\widecheck{x})$, $g_R(\widecheck{x})$ defined in (\ref{rates2D}). 
The two continuous time Markov chains $\Breve{X}$ and $\invbreve{X}$, defined on the same state space as $\widecheck{X}$ have infinitesimal generators $\Breve{Q}$ and $\invbreve{Q}$, respectively, such that, for $z\in \{\ell, \ell +1 ,\dots, u-1,u\}$ and $\widecheck{x}\in \L_z$, \textcolor{black}{$\Breve{Q}_{\widecheck{x},\widecheck{x}+\widecheck{v}_k} = \frac{\lambda^M_z}{|G_{+} (\widecheck{x})|}$ for $k \in G_{+} (\widecheck{x})$, $\Breve{Q}_{\widecheck{x},\widecheck{x}+\widecheck{v}_k} = \frac{\gamma^m_z}{|G_{-}(\widecheck{x})|}$ for $k \in G_{-} (\widecheck{x})$,
$\invbreve{Q}_{\widecheck{x},\widecheck{x}+\widecheck{v}_k} = \frac{\lambda^m_z}{|G_{+} (\widecheck{x})|}$ for $k \in G_{+} (\widecheck{x})$, and
$\invbreve{Q}_{\widecheck{x},\widecheck{x}+\widecheck{v}_k}=\frac{\gamma^M_z}{|G_{-} (\widecheck{x})|}$ for $k \in G_{-} (\widecheck{x})$, with $\lambda^{M}_z$, $\lambda^{m}_z$, $\gamma^{M}_z$, and $\gamma^{m}_z$ defined as in \eqref{minmaxrates}, and non-empty $G_{+} (\widecheck{x})$ and $G_{-} (\widecheck{x})$ defined as in \eqref{G+-(x)} where $G_{+}$ and $G_{-}$ are given in \eqref{G+-ex1}.}
 %
%$\Breve{Q}_{\widecheck{x},\widecheck{x}+\widecheck{v}_k}=\frac{\lambda^{M}_z}{2}$ for $k \in G_{+}$,
%$\Breve{Q}_{\widecheck{x},\widecheck{x}+\widecheck{v}_k}=\frac{\gamma^{m}_z}{2}$ for $k \in G_{-}$, $\invbreve{Q}_{\widecheck{x},\widecheck{x}+\widecheck{v}_k}=\frac{\lambda^{m}_z}{2}$ for $k \in G_{+}$, and
%$\invbreve{Q}_{\widecheck{x},\widecheck{x}+\widecheck{v}_k}=\frac{\gamma^{M}_z}{2}$ for $k \in G_{-}$, with $\lambda^{M}_z$, $\lambda^{m}_z$, $\gamma^{M}_z$, and $\gamma^{m}_z$ defined as in \eqref{minmaxrates}.
 %with a superscript $\eps$ appended to each $\lambda$ and $\gamma$ there.

 Then, as described in Section \ref{indeplevelstructure}, we can compare the Markov chain $\widecheck{X}$ with $\Breve X$ and $\invbreve{X}$, separately,
 %by using Theorems 3.3 \textcolor{black}{and 3.4} in \cite{Monotonicitypaper} \textcolor{black}{(see SI - Section \ref{compthmpaper_thm})}, with the matrix $A=[1\; -1]$ associated with the \textcolor{black}{coclique} level function $L(x_1,x_2)=x_1-x_2$. We then 
 to conclude that
\begin{equation}
\label{upandlowbounds2D}
    \Breve{h}_{\ell,u}\le h_{a,r}\le \invbreve{h}_{\ell,u}\;\;\;\mathrm{and}\;\;\;\invbreve{h}_{u,\ell}\le
     h_{r,a} \le 
     \Breve{h}_{u,\ell}
\end{equation}
where $\Breve{h}_{\ell,u}$, $\invbreve{h}_{\ell,u}$, $\invbreve{h}_{u,\ell}$, and $\Breve{h}_{u,\ell}$ can be written as in \eqref{formula3lowbound}, \eqref{formula3upbound}, \eqref{formulaDto0lowbound}, and \eqref{formulaDto0upbound}, respectively, with the quantities in \eqref{minmaxrates} replaced by those in \eqref{2Duplowbounddrates}.
%with a superscript of $\eps$ appended.
	
{\color{black}
These analytical expressions allow us to study the effect of the parameter $\eps$ on the MFPT from $a$ to $r$ (average time to memory loss of the active state) and the MFPT from $r$ to $a$ (average time to memory loss of the repressed state). Specifically, given that the only $O(\eps)$ rates are $\gamma^{M}_u=\gamma^{m}_u=g_R(\mathrm{D_{tot}},0)$ and $\lambda^{M}_{\ell}=\lambda^{m}_{\ell}=g_A(0,\mathrm{D_{tot}})$, with the other rates being $O(1)$, we can conclude that the upper bounds $\invbreve{h}_{\ell,u}$, $\Breve{h}_{u,\ell}$ and the lower bounds $\Breve{h}_{\ell,u}$, $\invbreve{h}_{u,\ell}$ in \eqref{upandlowbounds2D} are all $O(\eps^{-1})$. This implies that the average time to memory loss of both the repressed and active states are $O(\eps^{-1})$, and as $\eps$ approaches $0$, the average time to memory loss of both the repressed and active states tends to infinity.
}
%These analytical expressions allow us to determine the effect of system parameters on the bounds for the MFPT. Specifically, given that the only rates that are $O(\eps)$ are $\gamma^{m}_u=\gamma^{M}_u=g_R(\mathrm{D_{tot}},0)$ and $\lambda^{m}_{\ell}=\lambda^{M}_{\ell}=g_A(0,\mathrm{D_{tot}})$, then $\invbreve{h}_{u,\ell}$ and $\Breve{h}_{\ell,u}$, that are the lower bounds for the time to memory loss of the repressed and active state, respectively, are $O(\eps^{-1})$. This means that, as $\eps$ approaches 0, the bounds tend to infinity, implying that the time to memory loss of the repressed and active state also tends to infinity as $\eps \rightarrow 0$.

\subsection{Full chromatin modification circuit}\label{SS4}

\begin{figure}[t!]
	\centering
	\includegraphics[scale=0.39]{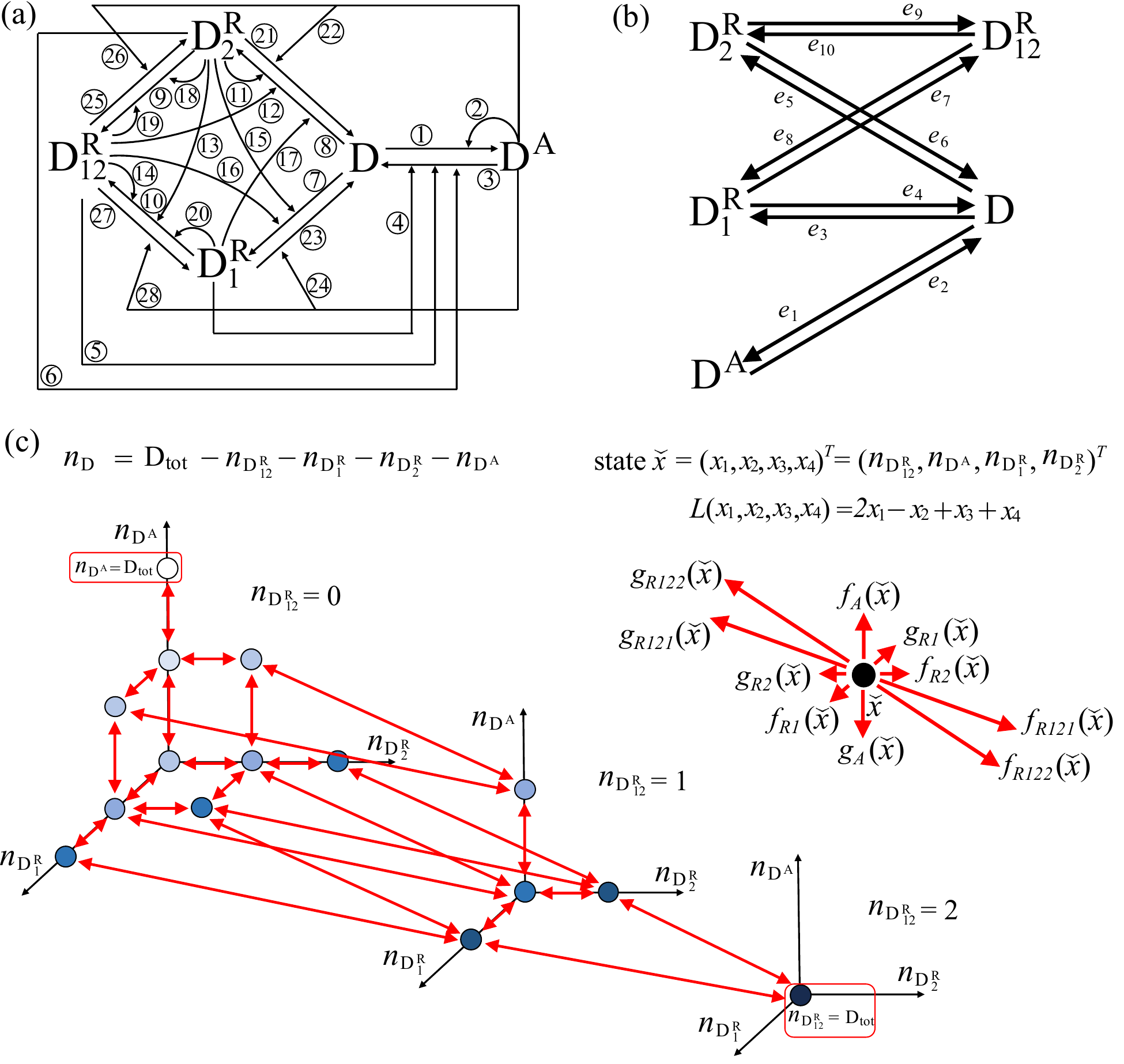}
	\caption{\small { \bf Full chromatin modification circuit: reaction diagram, graph $\G$ and associated Markov chain.}
(a) Chemical reaction system diagram. The numbers on the arrows correspond to the reactions associated with the arrows as described in (\ref{reacs4D}) in the main text.
(b) Graph $\G$ associated with the chemical reaction system in panel (a).
\textcolor{black}{(c)} State space and transitions of the projected continuous time Markov chain $\widecheck{X}=\{
(X_1(t),X_2(t),X_3(t),X_4(t))^T:\: t \ge 0\}$, which keeps track of $(n_{\mathrm{D^R_{12}}},n_{\mathrm{D^A}},n_{\mathrm{D^R_1}},n_{\mathrm{D^R_2}})$ through time. Here, we consider $\mathrm{D_{tot}}=2$ and we use dots to represent the states, and red double-ended arrows to represent transitions in both directions. 
Additionally, we use shades of blue to distinguish the level to which each state belongs. The function $L$ associated with the coclique level structure is $L(x_1,x_2,x_3,x_4)=2x_1-x_2+x_3+x_4$.
The rates associated with the one-step transitions for the projected Markov chain $\widecheck{X}$ are given in (\ref{rates4D}).
%(c) Directions of the possible one-step transitions for the projected Markov chain $\widecheck{X}$ starting from a state $\widecheck x=(x_1,x_2,x_3,x_4)^T$, whose rates are defined in (\ref{rates4D}).
}       
	\label{fig:4Dmodel1}
\end{figure}

In this example, we consider a chromatin modification circuit model that includes not only  histone modifications, but also DNA methylation \cite{bib:BWD2022,BrunoSontagPaper2023}. The model involves five species: unmodified nucleosome, denoted by D; nucleosome with CpGme only, denoted by $\mathrm{D^R_1}$; nucleosome with H3K9me3 only, denoted by $\mathrm{D^R_2}$; nucleosome with both H3K9me3 and CpGme, denoted by $\mathrm{D^R_{12}}$; nucleosome with an activating histone modification, denoted by $\mathrm{D^A}$. 
 In terms of molecular interactions, DNA methylation catalyzes the establishment of repressive histone modifications (and vice versa), while enhancing the erasure of activating marks (and vice versa) \cite{bib:BWD2022}. The quantity of each species is denoted by $n_{\mathrm{D}}$,  $n_{\mathrm{D^A}}$, $n_{\mathrm{D^R_1}}$, $n_{\mathrm{D^R_2}}$, and $n_{\mathrm{D^R_{12}}}$, respectively. Their sum remains constant, i.e., $n_{\mathrm{D}}+n_{\mathrm{D^A}}+n_{\mathrm{D^R_1}}+n_{\mathrm{D^R_2}}+n_{\mathrm{D^R_{12}}}=\Dtot$. The chemical reaction system, whose diagram is shown in Figure \ref{fig:4Dmodel1}(a), can be written as

\begingroup
\small 
\begin{align}
	\notag &{\large \textcircled{\small 1}}\;\ce{D ->[$k^A_{W0}+k^A_W$] D^A },\;\;\;{\large \textcircled{\small 2}}\;\ce{D + $\mathrm{D^A}$ ->[$k^A_M$] $\mathrm{D^A}$ + $\mathrm{D^A}$ },\;\;\;{\large \textcircled{\small 3}}\;\ce{$\mathrm{D^A}$ ->[$\delta + \bar k^A_{E}$] D },\\
	\notag &{\large \textcircled{\small 4}}\;\ce{$\mathrm{D^A}$ + $\mathrm{D^R_1}$ ->[$k^A_E$] D + $\mathrm{D^R_1}$ },\;\;\;{\large \textcircled{\small 5}}\;\ce{$\mathrm{D^A}$ + $\mathrm{D^R_{12}}$ ->[2 $k^A_E$] D + $\mathrm{D^R_{12}}$ },\;\;\;{\large \textcircled{\small 6}}\;\ce{$\mathrm{D^A}$ + $\mathrm{D^R_2}$ ->[$k^A_E$] D + $\mathrm{D^R_2}$ },\\
	\notag &{\large \textcircled{\small 7}}\;\ce{D ->[$k^1_{W0}+k^1_{W}$] $\mathrm{D^R_1}$ },\;\;\;{\large \textcircled{\small 8}}\;\ce{D ->[$k^2_{W0}+k^2_{W}$] $\mathrm{D^R_2}$ },\;\;\;{\large \textcircled{\small 9}}\;\ce{$\mathrm{D^R_2}$ ->[$k^1_{W0}$] $\mathrm{D^R_{12}}$ },\;\;\;{\large \textcircled{\small 10}}\;\ce{$\mathrm{D^R_1}$ ->[$k^2_{W0}$] $\mathrm{D^R_{12}}$ },
	\\
	\notag & {\large \textcircled{\small 11}}\;\ce{D + $\mathrm{D^R_2}$ ->[$k_M$] $\mathrm{D^R_2}$ + $\mathrm{D^R_2}$ },\;\;\;{\large \textcircled{\small 12}}\;\ce{D + $\mathrm{D^R_{12}}$ ->[$k_M + \bar k_M$] $\mathrm{D^R_2}$ + $\mathrm{D^R_{12}}$ },\\
	\label{reacs4D} &{\large \textcircled{\small 13}}\;\ce{$\mathrm{D^R_1}$ + $\mathrm{D^R_2}$ ->[$k_M$] $\mathrm{D^R_{12}}$ + $\mathrm{D^R_2}$ },\;\;\;{\large \textcircled{\small 14}}\;\ce{$\mathrm{D^R_1}$ + $\mathrm{D^R_{12}}$ ->[$k_M + \bar k_M$] $\mathrm{D^R_{12}}$ + $\mathrm{D^R_{12}}$ },\\
	\notag &{\large \textcircled{\small 15}}\;\ce{D + $\mathrm{D^R_2}$ ->[$k^{'}_M$] $\mathrm{D^R_1}$ + $\mathrm{D^R_2}$ },\;\;\;{\large \textcircled{\small 16}}\;\ce{D + $\mathrm{D^R_{12}}$ ->[$k^{'}_M$] $\mathrm{D^R_1}$ + $\mathrm{D^R_{12}}$ },\;\;\;{\large \textcircled{\small 17}}\;\ce{D + $\mathrm{D^R_1}$ ->[$\bar k_M$] $\mathrm{D^R_2}$ + $\mathrm{D^R_1}$ }\\
	\notag &{\large \textcircled{\small 18}}\;\ce{$\mathrm{D^R_2}$ + $\mathrm{D^R_2}$ ->[$k^{'}_M$] $\mathrm{D^R_{12}}$ + $\mathrm{D^R_2}$ },\;\;\;{\large \textcircled{\small 19}}\;\ce{$\mathrm{D^R_2}$ + $\mathrm{D^R_{12}}$ ->[$k^{'}_M$] $\mathrm{D^R_{12}}$ + $\mathrm{D^R_{12}}$ },\\
	\notag &{\large \textcircled{\small 20}}\;\ce{$\mathrm{D^R_1}$ + $\mathrm{D^R_1}$ ->[$\bar k_M$] $\mathrm{D^R_{12}}$ + $\mathrm{D^R_1}$ },\;\;\;
	{\large \textcircled{\small 21}}\;\ce{$\mathrm{D^R_2}$ ->[$\delta+\bar k^R_{E}$] D },\;\;\;{\large \textcircled{\small 22}}\;\ce{$\mathrm{D^R_2}$ + $\mathrm{D^A}$ ->[$k^R_E$] D + $\mathrm{D^A}$ },\\
	\notag &{\large \textcircled{\small 23}}\;\ce{$\mathrm{D^R_1}$ ->[$\delta^{'}+k^{'}_{T}$] D },\;\;\;{\large \textcircled{\small 24}}\;\ce{$\mathrm{D^R_1}$ + $\mathrm{D^A}$ ->[$k^{'*}_T$] D + $\mathrm{D^A}$ },\;\;\;
	{\large \textcircled{\small 25}}\;\ce{$\mathrm{D^R_{12}}$ ->[$\delta^{'}+k^{'}_{T}$] $\mathrm{D^R_2}$ },\\
	\notag &{\large \textcircled{\small 26}}\;\ce{$\mathrm{D^R_{12}}$ + $\mathrm{D^A}$ ->[$k^{'*}_T$] $\mathrm{D^R_2}$ + $\mathrm{D^A}$ },\;\;\;{\large \textcircled{\small 27}}\;\ce{$\mathrm{D^R_{12}}$ ->[$\delta+\bar k^R_{E}$] $\mathrm{D^R_1}$ },\;\;\;
	{\large \textcircled{\small 28}}\;\ce{$\mathrm{D^R_{12}}$ + $\mathrm{D^A}$ ->[$k^R_E$] $\mathrm{D^R_1}$ + $\mathrm{D^A}$ },
\end{align}
\endgroup
\normalsize
in which $k^A_{W0},k^A_W,k^A_M,\delta, \bar k^A_E,  k^A_E, k^1_{W0}, k^1_{W}, k^2_{W0}, k^2_{W}, k'_M, \bar k_M, k_M, \delta', k'_T, k^{'*}_T, \bar k^R_E, k^R_E > 0$ and the expression of the reaction rate constants is because we combined reactions sharing the same reactants and products. As we did for Example \ref{exp2}, let us introduce parameters $\eps=\frac{\delta+\bar{k}^A_E}{\frac{k^A_M}{V}\Dtot}$ and $\mu=\frac{k^R_E}{k^A_E}$, with a constant $\tilde{b}$ such that $\frac{\delta+\bar k^R_E}{\delta + \bar k^A_E} = \tilde{b}\mu$. Furthermore, since in this model we also have DNA methylation, we also introduce $\mu'=\frac{k^{'*}_T}{k^A_E}$ and a constant $\beta$ such that $\frac{\delta^{'}+k^{'}_T}{\delta + \bar k^A_E} = \beta \mu'$. The parameter $\mu'$ quantifies the relative speed between the erasure rate of DNA methylation and the erasure rate of activating histone modifications.

Considering $x=(n_{\mathrm{D^R_{12}}},n_{\mathrm{D^A}},n_{\mathrm{D^R_1}},n_{\mathrm{D^R_2}},n_{\mathrm{D}})^T$, the reaction vectors associated with \eqref{reacs4D} are $v_1=(0,1,0,0,-1)^T$, $v_2=(0,-1,0,0,1)^T$, $v_3=(0,0,1,0,-1)^T$, $v_4=(0,0,-1,0,1)^T$, $v_5=(0,0,0,1,-1)^T$, $v_6=(0,0,0,-1,1)^T$, $v_7=(1,0,-1,0,0)^T$, $v_8=(-1,0,1,0,0)^T$, $v_9=(1,0,0,-1,0)^T$, and $v_{10}=(-1,0,0,1,0)^T$. By examining them, one can verify that Assumption \ref{assumption:Unimolecular_change} is satisfied. The graph $\G$ associated with the chemical reaction system \eqref{reacs4D} can then be represented as in Figure \ref{fig:4Dmodel1}(b).
As done for Example \ref{exp2}, by inspecting $\G$, one can verify that the underlying undirected graph is connected and that $\G$ is bipartite. 
%Given that Assumption \ref{assumption:Unimolecular_change} is satisfied, we can then apply Lemma \ref{lemmabis}. This implies 
By Lemma \ref{lemmabis}, our SCRN has a unique conservation vector $m = (1,\ldots,1)^T$ and we can introduce the projected continuous time Markov chain $\widecheck{X}=\{
(X_1(t),X_2(t),X_3(t),X_4(t))^T:\: t \ge 0\}$, which keeps track of $(n_{\mathrm{D^R_{12}}},n_{\mathrm{D^A}},n_{\mathrm{D^R_1}},n_{\mathrm{D^R_2}})$ through time.
Since the total number of nucleosomes $\Dtot$ is conserved, the state space is $\widecheck{\X}= \{\widecheck{x}=(x_1,x_2,x_3,x_4)^T \in \Z_+^4 :\: x_1 + x_2 + x_3 + x_4\leq \Dtot \}$. The potential one-step transitions for $\widecheck X$ from $\widecheck x \in \widecheck{\X}$ are shown in Figure \ref{fig:4Dmodel1}(c), where the associated transition vectors are
$\widecheck{v}_1=-\widecheck{v}_2=(0,1,0,0)^T$, $\widecheck{v}_3=-\widecheck{v}_4=(0,0,1,0)^T$, $\widecheck{v}_5=-\widecheck{v}_6=(0,0,0,1)^T$, $\widecheck{v}_7=-\widecheck{v}_8=(1,0,-1,0)^T$, and $\widecheck{v}_9=-\widecheck{v}_{10}=(1,0,0,-1)^T$, and the infinitesimal transition rates (in which we assume mass-action kinetics with the usual rate constant volume scaling) are

\begingroup
\small 
\begin{align}
	\notag &\widecheck{Q}_{\widecheck{x},\widecheck{x}+\widecheck{v}_1}=f_A(\widecheck{x}) = (\Dtot -(x_1+x_2+x_3+x_4))\left(k_{W0}^A+k_{W}^A + \frac{k_{M}^A}{V}x_2\right),\\
	\notag &\widecheck{Q}_{\widecheck{x},\widecheck{x}+\widecheck{v}_2}=g_A(\widecheck{x}) = x_2\left(\eps \frac{k_{M}^A}{V}\Dtot + \frac{k^A_E}{V}(x_3+x_4+2x_1)\right),\\
	\notag &\widecheck{Q}_{\widecheck{x},\widecheck{x}+\widecheck{v}_3}=f_{R1}(\widecheck{x}) = (\Dtot -(x_1+x_2+x_3+x_4))\left(k^1_{W0}+k^1_{W} + \frac{k^{'}_{M}}{V}(x_1+x_4)\right),\\
	\notag &\widecheck{Q}_{\widecheck{x},\widecheck{x}+\widecheck{v}_4}=g_{R1}(\widecheck{x}) = x_3\mu'\left(\eps \frac{k_{M}^A}{V}\Dtot \beta + x_2\frac{k^A_E}{V}\right),\\
	\notag &\widecheck{Q}_{\widecheck{x},\widecheck{x}+\widecheck{v}_5}=f_{R2}(\widecheck{x}) = (\Dtot -(x_1+x_2+x_3+x_4))\left(k^2_{W0}+k^2_{W} +  \frac{k_{M}}{V}(x_1+x_4)  +  \frac{\bar k_{M}}{V}(x_1+x_3)\right),\\
	\label{rates4D} &Q_{\widecheck{x},\widecheck{x}+\widecheck{v}_6}=g_{R2}(\widecheck{x}) = x_4\mu\left(\eps \frac{k_{M}^A}{V}\Dtot \tilde{b} + x_2\frac{k^A_E}{V}\right),\\
	\notag &\widecheck{Q}_{\widecheck{x},\widecheck{x}+\widecheck{v}_7}=f_{R121}(\widecheck{x}) = x_3\left(k^2_{W0}+  \frac{k_{M}}{V}(x_1+x_4)  +  \frac{\bar k_{M}}{V}\left(x_1+\frac{x_3-1}{2}\right)\right),\\
	\notag &\widecheck{Q}_{\widecheck{x},\widecheck{x}+\widecheck{v}_8}=g_{R121}(\widecheck{x}) = x_1\mu\left(\eps \frac{k_{M}^A}{V}\Dtot b + x_2\frac{k^A_E}{V}\right), \\
	\notag &\widecheck{Q}_{\widecheck{x},\widecheck{x}+\widecheck{v}_9}=f_{R122}(\widecheck{x}) = x_4\left(k^1_{W0} + \frac{k^{'}_{M}}{V}\left(x_1+\frac{x_4-1}{2}\right)\right),\\
	\notag &\widecheck{Q}_{\widecheck{x},\widecheck{x}+\widecheck{v}_{10}}=g_{R122}(\widecheck{x}) = x_1\mu'\left(\eps \frac{k_{M}^A}{V}\Dtot \beta + x_2\frac{k^A_E}{V}\right).
\end{align}
\endgroup
\normalsize

 A representation of the Markov chain graph for $\Dtot=2$ is given in Figure \ref{fig:4Dmodel1}(c).
 We now focus on determining explicit analytical expressions for upper and lower bounds of the MFPT from the active state $a=(0,\mathrm{D_{tot}},0,0)^T$ to the repressed state $r=(\mathrm{D_{tot}},0,0,0)^T$, i.e., $h_{a,r} $, and vice versa, i.e., $h_{r,a}$, in order to understand how the parameters $\eps$ and $\mu'$ affect them.
 %
  %For this study, we will consider $k^A_W=k^1_W=k^2_W=0$ (i.e., there are no external transcription factors enhancing the establishment of chromatin modifications). This assumption will not change the qualitative nature of the results focused on studying the effect of $\eps$ and $\mu'$ on the MFPT.
 
 To calculate the MFPT upper and lower bounds, let us first apply Theorem \ref{thm2b} in order to determine the \textcolor{black}{coclique} level structures for $\widecheck{X}$. We can apply the theorem because Assumption \ref{assumption:Unimolecular_change} is satisfied
 %, the species of our SCRN are $d=5>2$,
 and the associated graph $\G$ is weakly connected. The rank of the stoichiometric matrix $S$ is $4=d-1$.
 Now, for this example, there are several possible partitions $\{ \Eu, \Ed \}$ of edges of $\G$ that could allow us to determine \textcolor{black}{a coclique} level structure. Let us consider the following one: $\Eu=\{\edge_2,\edge_3,\edge_5,\edge_7,\edge_9\},\Ed=\{\edge_1,\edge_4,\edge_6,\edge_8,\edge_{10}\}$. The reason for this choice is that it is the only one that, as we will see later, allows us to determine \textcolor{black}{a coclique} level structure in which the active state $a=(0,\mathrm{D_{tot}},0,0)^T$ and the repressed state $r=(\mathrm{D_{tot}},0,0,0)^T$ are the two extremum levels. For this partition, 
 %let us now write the system of equations as the one in \eqref{syst}, obtaining
 %\begin{align}
  %     b_2&=-1,\;-b_2=1,\;b_3=1,\;-b_3=-1,\;b_4=1,\;-b_4=-1,\\
   % \notag   -b_1+b_3&=-1,\;b_1-b_3=1,\;-b_1+b_4=-1,\;b_1-b_4=1.
 %\end{align}
 %%
%The system of equations written above admits a solution $b \in \Z ^4$, that is, $(b_1,b_2,b_3,b_4)=(2,-1,1,1)$, respectively.
the system of equations in \eqref{syst} admits a unique solution $(b_1,b_2,b_3,b_4)^T=(2,-1,1,1)^T$.

Then, by applying Theorem \ref{thm2b}, we can conclude that the projected Markov chain $\widecheck{X}$ has \textcolor{black}{a coclique} level structure associated with the \textcolor{black}{coclique} level function $L(x_1,x_2,x_3,x_4)=2x_1-x_2+x_3+x_4$.
 This \textcolor{black}{coclique} level structure can be written as $\L_{\ell}, \dots, \L_u$, with $\L_z:=\{x\in \widecheck{\X}:\: L(x)=2x_1-x_2+x_3+x_4=z\}$ for $z=\ell,\dots,u$, with $\ell=-\mathrm{D_{tot}}$ and $u=2\mathrm{D_{tot}}$ (Figure \ref{fig:4Dmodel1}(c)). This \textcolor{black}{coclique} level structure is such that $a=(0,\mathrm{D_{tot}},0,0)^T$ is the only state belonging to $\L_\ell$ and $r=(\mathrm{D_{tot}},0,0,0)^T$ is the only state belonging to $\L_u$. Here, 
 \begin{equation}
\label{G+-ex2}
G_{+}=\{2,3,5,7,9\}\;\;\;\mathrm{and}\;\;\;G_{-}=\{1,4,6,8,10\},
\end{equation}
and the rate of increase $\lambda_z(\widecheck{x})$ and the rate of decrease $\gamma_z(\widecheck{x})$ can then be written as  
\begin{equation}\label{4Duplowbounddrates}
	\begin{aligned}	&\lambda_z(\widecheck{x})=f_{R121}(\widecheck{x})+f_{R122}(\widecheck{x})+g_A(\widecheck{x})+f_{R1}(\widecheck{x})+f_{R2}(\widecheck{x}),\\
		&\gamma_z(\widecheck{x})=f_A(\widecheck{x})+g_{R121}(\widecheck{x})+g_{R122}(\widecheck{x})+g_{R1}(\widecheck{x})+g_{R2}(\widecheck{x}),
	\end{aligned}
\end{equation}
respectively, with $f_{R121}(\widecheck{x})$, $f_{R122}(\widecheck{x})$, $g_A(\widecheck{x})$, $f_{R1}(\widecheck{x})$, $f_{R2}(\widecheck{x})$, $f_A(\widecheck{x})$, $g_{R121}(\widecheck{x})$, $g_{R122}(\widecheck{x})$, $g_{R1}(\widecheck{x})$, $g_{R2}(\widecheck{x})$ defined in (\ref{rates4D}).

The two continuous time Markov chains $\Breve{X}$ and $\invbreve{X}$ are defined on the same state space as $\widecheck{X}$ and have infinitesimal generators $\Breve{Q}$ and $\invbreve{Q}$, respectively, such that, for $z\in \{\ell, \ell +1 ,\dots, u-1,u\}$ and $\widecheck{x}\in \L_z$, \textcolor{black}{$\Breve{Q}_{\widecheck{x},\widecheck{x}+\widecheck{v}_k} = \frac{\lambda^M_z}{|G_{+} (\widecheck{x})|}$ for $k \in G_{+} (\widecheck{x})$, $\Breve{Q}_{\widecheck{x},\widecheck{x}+\widecheck{v}_k} = \frac{\gamma^m_z}{|G_{-}(\widecheck{x})|}$ for $k \in G_{-} (\widecheck{x})$,
$\invbreve{Q}_{\widecheck{x},\widecheck{x}+\widecheck{v}_k} = \frac{\lambda^m_z}{|G_{+} (\widecheck{x})|}$ for $k \in G_{+} (\widecheck{x})$, and
$\invbreve{Q}_{\widecheck{x},\widecheck{x}+\widecheck{v}_k}=\frac{\gamma^M_z}{|G_{-} (\widecheck{x})|}$ for $k \in G_{-} (\widecheck{x})$, with $\lambda^{M}_z=\max_{\widecheck{x}\in \mathcal{L}_z}\lambda_z(\widecheck{x})$, $\lambda^{m}_z=\min_{\widecheck{x}\in \mathcal{L}_z}\lambda_z(\widecheck{x}),$ $\gamma^{M}_z=\max_{\widecheck{x}\in \mathcal{L}_z}\gamma_z(\widecheck{x})$, and $\gamma^{m}_z=\min_{\widecheck{x}\in \mathcal{L}_z}\gamma_z(\widecheck{x})$, as defined in \eqref{minmaxrates}, and non-empty $G_{+} (\widecheck{x})$ and $G_{-} (\widecheck{x})$ defined as in \eqref{G+-(x)} where $G_{+}$ and $G_{-}$ are given in \eqref{G+-ex2}.}
%
%$\Breve{Q}_{\widecheck{x},\widecheck{x}+\widecheck{v}_k}=\frac{\lambda^{M}_z}{5}$ for $k \in G_{+}$,
%$\Breve{Q}_{\widecheck{x},\widecheck{x}+\widecheck{v}_k}=\frac{\gamma^{m}_z}{5}$ for $k \in G_{-}$,
%$\invbreve{Q}_{\widecheck{x},\widecheck{x}+\widecheck{v}_k}=\frac{\lambda^{m}_z}{5}$ for $k \in G_{+}$, and
%$\invbreve{Q}_{\widecheck{x},\widecheck{x}+\widecheck{v}_k}=\frac{\gamma^{M}_z}{5}$ for $k \in G_{-}$, with $\lambda^{M}_z=\max_{\widecheck{x}\in \mathcal{L}_z}\lambda_z(\widecheck{x})$, $\lambda^{m}_z=\min_{\widecheck{x}\in \mathcal{L}_z}\lambda_z(\widecheck{x}),$ $\gamma^{M}_z=\max_{\widecheck{x}\in \mathcal{L}_z}\gamma_z(\widecheck{x})$, and $\gamma^{m}_z=\min_{\widecheck{x}\in \mathcal{L}_z}\gamma_z(\widecheck{x})$, as defined in \eqref{minmaxrates}.

 Then, as described in Section \ref{indeplevelstructure}, we can compare the Markov chain $\widecheck{X}$ with $\Breve X$ and $\invbreve{X}$, separately,
 %by using Theorems 3.3 \textcolor{black}{and 3.4} in \cite{Monotonicitypaper} \textcolor{black}{(see SI - Section \ref{compthmpaper_thm})}, with the matrix $A=[2\; -1\;\;\; 1\;\;\; 1]$  associated with the \textcolor{black}{coclique} level function $L(x_1,x_2,x_3,x_4)=2x_1-x_2+x_3+x_4$. Then, 
 to obtain analytical expressions for lower and upper bounds, respectively, for the MFPT $h_{a,r}$ from the fully active state to the fully repressed state:
%\begin{equation}\label{lowboundMFPTactstate4D}
%	\begin{aligned}
%		&\Breve{h}_{\ell,u}=\frac{\tilde r^{m}_{u-(\ell+1)}}{\lambda^{M}_{\ell}}\left(1+\sum_{j=1}^{u-(\ell+1)}\frac{1}{\tilde r^{m}_j}\right)+\sum_{i=2}^{u-(\ell+1)}\left[\frac{\tilde r^{m}_{i-1}}{\lambda^{M}_{u-i}}\left(1+\sum_{j=1}^{i-1}\frac{1}{\tilde r^{m}_j}\right)\right]+\frac{1}{\lambda^{M}_{u-1}},\\
%		&\invbreve{h}_{\ell,u}=\frac{\tilde r^{M}_{u-(\ell+1)}}{\lambda^{m}_{\ell}}\left(1+\sum_{j=1}^{u-(\ell+1)}\frac{1}{\tilde r^{M}_j}\right)+\sum_{i=2}^{u-(\ell+1)}\left[\frac{\tilde r^{M}_{i-1}}{\lambda^{m}_{u-i}}\left(1+\sum_{j=1}^{i-1}\frac{1}{\tilde r^{M}_j}\right)\right]+\frac{1}{\lambda^{m}_{u-1}},\nonumber
%	\end{aligned}
%\end{equation}
%respectively, with $\tilde r^{m}_j=\frac{\gamma^{m}_{u-1}\gamma^{m}_{u-2}\dots\gamma^{m}_{u-j}}{\lambda^{M}_{u-1}\lambda^{M}_{u-2}\dots\lambda^{M}_{u-j}}$ and $\tilde r^{M}_j=\frac{\gamma^{M}_{u-1}\gamma^{M}_{u-2}\dots\gamma^{M}_{u-j}}{\lambda^{m}_{u-1}\lambda^{m}_{u-2}\dots\lambda^{m}_{u-j}}$, 
{\color{black}
\begin{equation}\label{lowboundMFPTactstate4D}
	\begin{aligned}
&\Breve{h}_{\ell,u}=\frac{1}{\lambda^M_{u-1}} + \sum_{i=\ell}^{u-2} \frac{1}{\lambda^M_{i}} \left( 1 + \sum_{j=i+1}^{u-1} \frac{\gamma^m_{i+1}\dots\gamma^m_{j}}{\lambda^M_{i+1}\dots\lambda^M_{j}} \right),\\
&\invbreve{h}_{\ell,u}=\frac{1}{\lambda^m_{u-1}} + \sum_{i=\ell}^{u-2} \frac{1}{\lambda^m_{i}} \left( 1 + \sum_{j=i+1}^{u-1} \frac{\gamma^M_{i+1}\dots\gamma^M_{j}}{\lambda^m_{i+1}\dots\lambda^m_{j}} \right),\nonumber
	\end{aligned}
\end{equation}}
and expressions for lower and upper bounds for the MFPT $h_{r,a}$ from the fully repressed state to the fully active state:
%\begin{equation}\label{lowupboundMFPTrepstate4D}
%	\begin{aligned}
%		&\invbreve{h}_{u,\ell} =\frac{r^{m}_{u-1}}{\gamma^{M}_{u}}\left(1+\sum_{j=\ell+1}^{u-1}\frac{1}{r^{m}_j}\right)+\sum_{i=\ell+2}^{u-1}\left[\frac{r^{m}_{i-1}}{\gamma^{M}_i}\left(1+\sum_{j=\ell+1}^{i-1}\frac{1}{r^{m}_j}\right)\right]+\frac{1}{\gamma^{M}_{\ell+1}},\\ 
%		&\Breve{h}_{u,\ell} =\frac{r^{M}_{u-1}}{\gamma^{m}_{u}}\left(1+\sum_{j=\ell+1}^{u-1}\frac{1}{r^{M}_j}\right)+\sum_{i=\ell+2}^{u-1}\left[\frac{r^{M}_{i-1}}{\gamma^{m}_i}\left(1+\sum_{j=\ell+1}^{i-1}\frac{1}{r^{M}_j}\right)\right]+\frac{1}{\gamma^{m}_{\ell+1}},\nonumber
%	\end{aligned}
%\end{equation}	
%
%respectively, with $r^{M}_j=\frac{\lambda^{M}_{\ell+1} \lambda^{M}_{\ell+2}\dots\lambda^{M}_{j}}{\gamma^{m}_{\ell+1} \gamma^{m}_{\ell+2} \dots \gamma^{m}_{j}}$ and $r^{m}_j=\frac{\lambda^{m}_{\ell+1} \lambda^{m}_{\ell+2}\dots\lambda^{m}_{j}}{\gamma^{M}_{\ell+1} \gamma^{M}_{\ell+2} \dots \gamma^{M}_{j}}$
{\color{black}
\begin{equation}\label{lowupboundMFPTrepstate4D}
	\begin{aligned}
    & \invbreve{h}_{u,\ell}=\frac{1}{\gamma^M_{\ell+1}} + \sum_{i=\ell+2}^{u} \frac{1}{\gamma^M_{i}} \left( 1 + \sum_{j=\ell+1}^{i-1} \frac{\lambda^m_{j}\dots\lambda^m_{i-1}}{\gamma^M_{j}\dots\gamma^M_{i-1}} \right),\\
    &\Breve{h}_{u,\ell}=\frac{1}{\gamma^m_{\ell+1}} + \sum_{i=\ell+2}^{u} \frac{1}{\gamma^m_{i}} \left( 1 + \sum_{j=\ell+1}^{i-1} \frac{\lambda^M_{j}\dots\lambda^M_{i-1}}{\gamma^m_{j}\dots\gamma^m_{i-1}} \right). \nonumber
	\end{aligned}
\end{equation}}
{\color{black}
Given that the only $O(\eps)$ rates are $\lambda^{M}_{\ell}$, $\lambda^{m}_{\ell}$, $\gamma^{M}_u$, $\gamma^{M}_{u-1}$, and $\gamma^{m}_{z}:\Dtot \leq z \leq u$, with the other rates being $O(1)$, we can conclude that both $\Breve{h}_{\ell,u}$ and $\invbreve{h}_{\ell,u}$ are $O(\eps^{-1})$, and thus the MFPT from $a$ to $r$ (average time to memory loss of the active state) is $O(\eps^{-1})$. Furthermore, we have $\invbreve{h}_{u,\ell}$ is $O(\eps^{-2})$ and $\Breve{h}_{u,\ell}$ is $O(\eps^{-\Dtot})$, and thus the MFPT from $r$ to $a$ (average time to memory loss of the repressed state) is at least $O(\eps^{-2})$. In cases like this, alternative approaches are needed to identify the precise scalings of MFPTs, such as the method we developed in our recent work \cite{epifinite}, which allows us to show that average time to memory loss of the repressed state is indeed $O(\eps^{-2})$.
}
%
%Given that $\lambda^{M}_z$ is $O(1)$ for all $z$ except $z=\ell$ (in which case $\lambda^{M}_{\ell}=g_A(0,\mathrm{D_{tot}},0,0)=O(\eps)$) and $\gamma^{M}_z$ is $O(1)$ for all $z$ except $z=u$ (in which case $\gamma^{M}_u=g_{R121}(\mathrm{D_{tot}},0,0,0)+g_{R122}(\mathrm{D_{tot}},0,0,0)=O(\eps)$) and $z=u-1$ (in which case $\gamma^{M}_{u-1}=\max_{\widecheck{x}\in \L_{u-1}}\gamma_z(\widecheck{x})=\max\{g_{R121}(\mathrm{D_{tot}}-1,0,1,0)+g_{R122}(\mathrm{D_{tot}}-1,0,1,0)+g_{R1}(\mathrm{D_{tot}}-1,0,1,0),g_{R121}(\mathrm{D_{tot}}-1,0,0,1)+g_{R122}(\mathrm{D_{tot}}-1,0,0,1)+g_{R1}(\mathrm{D_{tot}}-1,0,0,1)\}=O(\eps)$), we can conclude that $\Breve{h}_{\ell,u}$ and $\invbreve{h}_{u,\ell}$, i.e., the lower bounds for the two MFPTs under consideration, are $O(\eps^{-1})$ and $O(\eps^{-2})$, respectively.
%
These results suggest that decreasing $\eps$ extends the memory of both the active and repressed chromatin states, but with a more pronounced impact on the repressed state. This difference can be attributed to the cooperation of repressive chromatin marks, i.e., DNA methylation and repressive histone modifications, which introduces a structural bias in the chromatin modification circuit towards a repressed chromatin state. 

Let us now determine the effect of $\mu'$, i.e., the parameter quantifying the relative speed between the DNA methylation erasure rate and the activating histone modification erasure rate, on the time to memory loss. Since the rates $g_{R1}(\widecheck{x})$ and $g_{R121}(\widecheck{x})$ are linear in $\mu'$ and are the only transition rates depending on $\mu'$ (see \eqref{rates4D}), then, based on the definition in (\ref{4Duplowbounddrates}), $\gamma_z(\widecheck{x})$ increases for lower values of $\mu'$. This implies that increasing $\mu'$ leads to higher $\invbreve{h}_{\ell,u}$ and $\Breve{h}_{\ell,u}$ and lower $\Breve{h}_{u,\ell}$ and $\invbreve{h}_{u,\ell}$. The opposite happens when $\mu'$ decreases.

	\subsection{Bi-parallel network motif}\label{exp3}

\begin{figure}[t]
            \centering
            \includegraphics[scale=0.39]{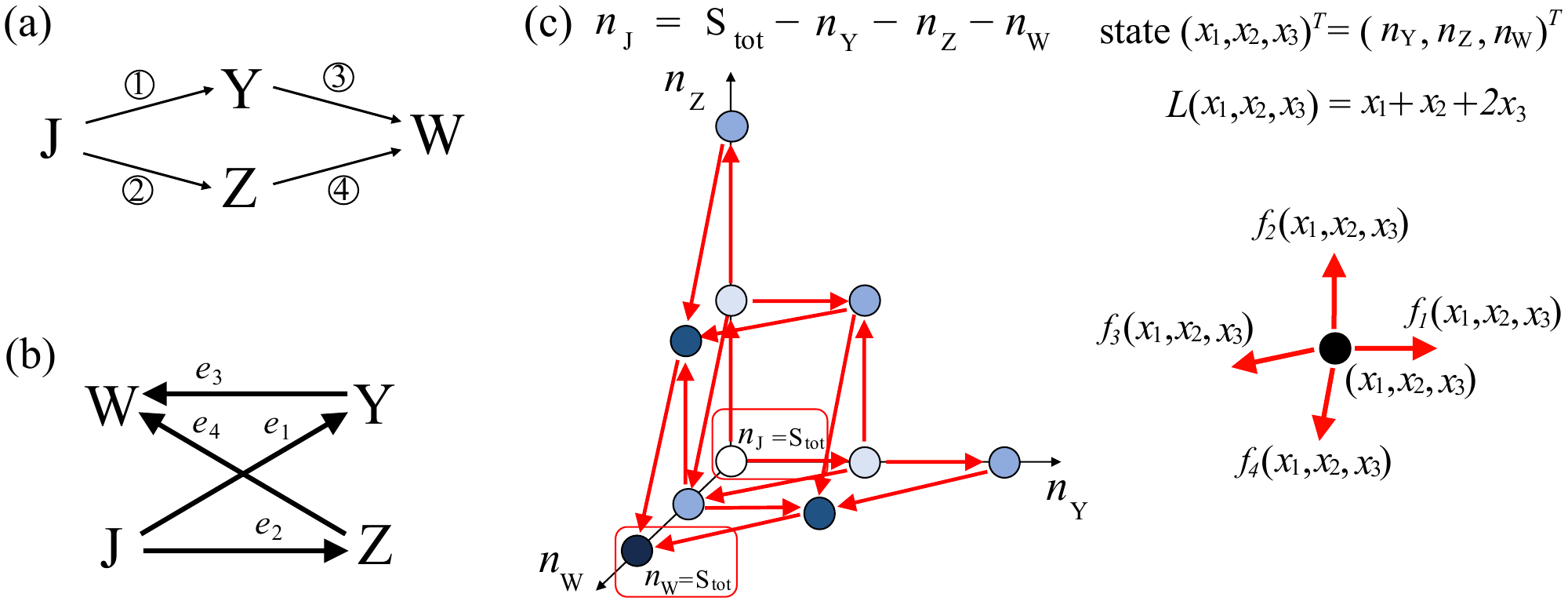}
            \caption{\small { \bf Bi-parallel network motif: reaction diagram, graph $\G$ and associated Markov chain.}
(a) Chemical reaction system diagram.  The numbers on the arrows correspond to the reactions associated with the arrows as described in (\ref{reacs2D}) in the main text. 
(b) Graph $\G$ associated with the chemical reaction system in panel (a).
(c) State space and transitions for the projected continuous time Markov chain $\widecheck{X}=\{
(X_1(t),X_2(t),X_3(t))^T:\: t \ge 0\}$, which keeps track of $(n_{\mathrm{Y}},n_{\mathrm{Z}},n_{\mathrm{W}})$ through time. Here, we consider $\mathrm{S_{tot}}=2$ and we use dots to represent the states, and red double-ended arrows to represent transitions in both directions. 
Additionally, we use shades of blue to distinguish the coclique level to which each state belongs. The function $L(x_1,x_2,x_3)$ associated with the level structure is $L(x_1,x_2,x_3)=x_1+x_2+2x_3$.
%
%(c)
%
The rates associated with the one-step transitions for the projected Markov chain $\widecheck{X}$ are given in (\ref{rates2D}).
            }
   \label{fig:3rdexp}
        \end{figure}

In this example, we analyze a bi-parallel network motif, which is a typical building block found in complex networks \cite{RMilo2002}, such as the neuronal connectivity network of the nematode Caenorhabditis elegans or food web networks \cite{RMilo2002,FBrglez1989}. The model includes four species, which are J, Y, Z, and W, and the quantity of each species is denoted by $n_{\mathrm{J}}$, $n_{\mathrm{Y}}$, $n_{\mathrm{Z}}$, and $n_{\mathrm{W}}$, respectively. Their sum remains constant, that is $n_{\mathrm{J}}+n_{\mathrm{Y}}+n_{\mathrm{Z}}+n_{\mathrm{W}}=\Stot$. Then, the chemical reaction system, whose diagram is shown in Figure \ref{fig:3rdexp}(a), can be written as
	\begin{equation}\label{reacs3rdexp}
	\begin{aligned}
	&{\large \textcircled{\small 1}}\;\ce{J ->[$k_1$] Y },\;\;{\large \textcircled{\small 2}}\;\ce{J ->[$k_2$] Z},\;\;{\large \textcircled{\small 3}}\;\ce{Y ->[$k_3$] W},\;\;{\large \textcircled{\small 4}}\;\ce{Z ->[$k_4$] W},
	\end{aligned}
	\end{equation}
in which $k_1,k_2,k_3,k_4 > 0$. 

Considering $x=(n_{\mathrm{Y}},n_{\mathrm{Z}},n_{\mathrm{W}},n_{\mathrm{J}})$, the reaction vectors associated with \eqref{reacs3rdexp} are $v_1=(1,0,0,-1)^T$, $v_2=(0,1,0,-1)^T$, $v_3=(-1,0,1,0)^T$, and $v_4=(0,-1,1,0)^T$. By examining them, we see that Assumption \ref{assumption:Unimolecular_change} is satisfied. The graph $\G$ associated with the chemical reaction system \eqref{reacs3rdexp} can be represented as in Figure \ref{fig:3rdexp}(b).
By inspecting $\G$, we see that the underlying undirected graph is connected and that $\G$ is bipartite. 
%Given that Assumption \ref{assumption:Unimolecular_change} is satisfied, we can then apply Lemma \ref{lemmabis}. This implies 
By Lemma \ref{lemmabis}, our SCRN has a unique conservation vector $m=(1,1,1,1)^T$ and then we can introduce a projected continuous time Markov chain $\widecheck{X}=\{
(X_1(t),X_2(t),X_3(t))^T:\: t \ge 0\}$, which keeps track of $(n_{\mathrm{Y}},n_{\mathrm{Z}},n_{\mathrm{W}})$ through time.
Since $\Stot$ is conserved, the state space is $\widecheck{\X}= \{\widecheck{x}=(x_1,x_2,x_3)^T \in \Z_+^3 :\: x_1 + x_2 + x_3 \leq \Stot \}$. The potential one-step transitions for $\widecheck X$ from $\widecheck x \in \widecheck{\X}$ are shown in Figure \ref{fig:3rdexp}(c), where the associated transition vectors are $\widecheck{v}_1=(1,0,0)^T$, $\widecheck{v}_2=(0,1,0)^T$, $\widecheck{v}_3=(-1,0,1)^T$, and $\widecheck{v}_4=(0,-1,1)^T$, and the infinitesimal transition rates (in which we assume mass-action kinetics) are
   \begin{equation}\label{rates3rdexp}
    \begin{aligned}
&\widecheck{Q}_{\widecheck{x},\widecheck{x}+\widecheck{v}_1}=f_1(\widecheck{x}) = k_1(\Stot -(x_1+x_2+x_3)),\\
&\widecheck{Q}_{\widecheck{x},\widecheck{x}+\widecheck{v}_2}=f_2(\widecheck{x}) = k_2(\Stot -(x_1+x_2+x_3)),\\
&\widecheck{Q}_{\widecheck{x},\widecheck{x}+\widecheck{v}_3}=f_3(\widecheck{x}) = k_3 x_1,\;\; \widecheck{Q}_{\widecheck{x},\widecheck{x}+\widecheck{v}_4}=f_4(\widecheck{x}) = k_4x_2.
       % g^{\eps}_R(x_1,x_2) &= i\left(\eps \frac{k_{M}^A}{V}\Dtot\mu b + j\frac{k^R_E}{V}\right)
    \end{aligned}
   \end{equation}
Let us now focus on determining explicit analytical expressions for upper and lower bounds of MFPTs. To this end, we apply Theorem \ref{thm2b} to determine the \textcolor{black}{coclique} level structures for $\widecheck{X}$. We can apply the theorem because Assumption \ref{assumption:Unimolecular_change} is satisfied,
%the species of our SCRN are $d=4>2$,
  and the associated graph $\G$ is weakly connected. The rank of the stoichiometric matrix $S$ is $3=d-1$.

Now, consider all of the possible partitions $\{ \Eu, \Ed \}$ of edges of $\G$ that could allow us to determine \textcolor{black}{a coclique} level structure. These partitions are the following:
\begin{align}
   \notag &\Eu=\{\edge_1,\edge_2,\edge_3,\edge_4\},\Ed=\emptyset;\;\;\;\Eu=\{\edge_1\},\Ed=\{\edge_2,\edge_3,\edge_4\};\\
    \notag &\Eu=\{\edge_2\},\Ed=\{\edge_1,\edge_3,\edge_4\};\;\;\;\Eu=\{\edge_3\},\Ed=\{\edge_1,\edge_2,\edge_4\};\\
    \notag &\Eu=\{\edge_4\},\Ed=\{\edge_1,\edge_2,\edge_3\};\;\;\;\Eu=\{\edge_1,\edge_2\},\Ed=\{\edge_3,\edge_4\};\\
    &\Eu=\{\edge_1,\edge_3\},\Ed=\{\edge_2,\edge_4\};\;\;\;\Eu=\{\edge_1,\edge_4\},\Ed=\{\edge_2,\edge_3\}.
\end{align}
As was done for Example \ref{exp2}, we did not consider the partitions obtained by switching the labels of the subsets $\Eu$, $\Ed$ in the partitions listed above because, as explained in Remark \ref{correlationpartitionstructure}, the associated functions $L$ would be the opposite of the ones obtained for the partitions considered above. Therefore, the resulting \textcolor{black}{coclique} level structures may be considered to be the same.
%
%Furthermore, we did not consider the partitions in which two edges in the same element of a partition are associated with two reaction vectors $\widecheck{v}_k, \widecheck{v}_i \in \V$ such that $\widecheck{v}_k = -\widecheck{v}_i$, because, as stated in Corollary \ref{thm2bis}, these partitions would not lead to \textcolor{black}{a coclique} level structure. 

For each partition, we can write the system of equations as in \eqref{syst}. The only systems that admit a solution $b \in \Z ^3$ are the ones associated with the first, sixth, and  \textcolor{black}{eighth partitions}. The solutions are $(b_1,b_2,b_3)^T=(1,1,2)^T$, $(b_1,b_2,b_3)^T=(1,1,0)^T$,  \textcolor{black}{and $(b_1,b_2,b_3)^T=(1,-1,0)^T$,} respectively. 
%obtaining
%\begin{align}
%b_1 =1,\;-b_1=-1,\;b_2=1,\;-b_2=-1
%\end{align}
%%
%and
%%
%\begin{align}
%b_1 =1,\;-b_1=-1,\;b_2=-1,\;-b_2=1,
%\end{align}
%%
%respectively. Both systems admit a solution \textcolor{black}{$b \in \Z ^2$}, that is, $(b_1,b_2)=(1,1)$ and $(b_1,b_2)=(1,-1)$, respectively. 
Then, by applying Theorem \ref{thm2b}, we can conclude that the projected Markov chain $\widecheck{X}$ has \textcolor{black}{three} \textcolor{black}{coclique} level structures associated with the \textcolor{black}{coclique} level functions $L(x_1,x_2,x_3)=x_1+x_2+2x_3$, $L(x_1,x_2,x_3)=x_1+x_2$,  \textcolor{black}{and $L(x_1,x_2,x_3)=x_1-x_2$}.

We consider the \textcolor{black}{coclique} level function $L(x_1,x_2,x_3)=x_1+x_2+2x_3$. The \textcolor{black}{coclique} level structure associated with it can be written as $\L_{\ell}, \dots, \L_u$, with $\L_z:=\{x\in \widecheck{\X}:\: L(x_1,x_2,x_3)=x_1 + x_2+2x_3=z\}$ for $z=\ell,\dots,u$, with $\ell=0$ and $u=2\Stot$. This \textcolor{black}{coclique} level structure is such that $(0,0,0)^T$ (i.e., $n_{\mathrm{J}}=\Stot$) is the only state belonging to $\L_\ell$ and $(0,0,\Stot)^T$ (i.e., $n_{\mathrm{W}}=\Stot$) is the only state belonging to $\L_u$ (Figure \ref{fig:3rdexp}(c)). This feature, as shown in the previous examples, is critical in order to determine good lower and upper bounds for the MFPT from $n_{\mathrm{J}}=\Stot$ to $n_{\mathrm{W}}=\Stot$.
%, and vice versa and this is the reason why we consider the \textcolor{black}{coclique} level structure associated with the function $L(x_1,x_2)=x_1-x_2$ and not the one associated with the function $L(x_1,x_2)=x_1+x_2$. 

Let us now determine the lower and upper bounds for the MFPT from $n_{\mathrm{J}}=\Stot$ to $n_{\mathrm{W}}=\Stot$. In particular, here we have 
 \begin{equation}
\label{G+-ex3}
G_{+}=\{1,2,3,4\}\;\;\;\mathrm{and}\;\;\;G_{-}=\emptyset,
\end{equation}
and the rate of increase $\lambda_z(\widecheck{x})$ and the rate of decrease $\gamma_z(\widecheck{x})$ can be written as 	
\begin{equation}\label{BiParalleluplowbounddrates}	\lambda_z(\widecheck{x})=f_1(\widecheck{x})+f_2(\widecheck{x})+f_3(\widecheck{x})+f_4(\widecheck{x})\;\;\;\mathrm{and}\;\;\;\gamma_z(\widecheck{x})=0, 
	\end{equation}
with $f_1(\widecheck{x})$, $f_2(\widecheck{x})$, $f_3(\widecheck{x})$, $f_4(\widecheck{x})$ defined in (\ref{rates3rdexp}). 

The two continuous time Markov chains $\Breve{X}$ and $\invbreve{X}$, are defined on the same state space as $\widecheck{X}$ and have infinitesimal generators $\Breve{Q}$ and $\invbreve{Q}$, respectively, such that, for $z\in \{\ell, \ell +1 ,\dots, u-1,u\}$ and $\widecheck{x}\in \L_z$, \textcolor{black}{$\Breve{Q}_{\widecheck{x},\widecheck{x}+\widecheck{v}_k} = \frac{\lambda^M_z}{|G_{+} (\widecheck{x})|}$ for $k \in G_{+} (\widecheck{x})$ and 
$\invbreve{Q}_{\widecheck{x},\widecheck{x}+\widecheck{v}_k} = \frac{\lambda^m_z}{|G_{+} (\widecheck{x})|}$ for $k \in G_{+} (\widecheck{x})$, with $\lambda^{M}_z=\max_{\widecheck{x}\in \mathcal{L}_z}\lambda_z(\widecheck{x})$ and $\lambda^{m}_z=\min_{\widecheck{x}\in \mathcal{L}_z}\lambda_z(\widecheck{x})$, as defined in \eqref{minmaxrates},
and non-empty $G_{+} (\widecheck{x})$ defined as in \eqref{G+-(x)} where $G_{+}$ is given in \eqref{G+-ex3}.}
%
%$\Breve{Q}_{\widecheck{x},\widecheck{x}+\widecheck{v}_k}=\frac{\lambda^{M}_z}{4}$ for $k \in G_{+}$,
% $\Breve{Q}_{\widecheck{x},\widecheck{x}+\widecheck{v}_k}=0$ for $k \in G_{-}$,
%$\invbreve{Q}_{\widecheck{x},\widecheck{x}+\widecheck{v}_k}=\frac{\lambda^{m}_z}{4}$ for $k \in G_{+}$, and
% $\invbreve{Q}_{\widecheck{x},\widecheck{x}+\widecheck{v}_k}=0$ for $k \in G_{-}$, with $\lambda^{M}_z=\max_{\widecheck{x}\in \mathcal{L}_z}\lambda_z(\widecheck{x})$, and $\lambda^{m}_z=\min_{\widecheck{x}\in \mathcal{L}_z}\lambda_z(\widecheck{x})$, as defined in \eqref{minmaxrates}.

 Then, as described in Section \ref{indeplevelstructure}, we can compare the Markov chain $\widecheck{X}$ with $\Breve X$ and $\invbreve{X}$, separately,
 %by using Theorems 3.3 \textcolor{black}{and 3.4} in \cite{Monotonicitypaper} \textcolor{black}{(see SI - Section \ref{compthmpaper_thm})}, with the matrix $A=[1\; 1\; 2]$  associated with the \textcolor{black}{coclique} level function $L(x_1,x_2,x_3)=x_1+x_2+2x_3$. Then, the 
 to obtain analytical expressions for lower and upper bounds for the MFPT from $n_{\mathrm{J}}=\Stot$ to $n_{\mathrm{W}}=\Stot$, %$h_{J=\Stot,W=\Stot}$,
 which can be written as
\begin{equation}\label{lowboundMFPTXW3rdexp}
	\begin{aligned}
		&\Breve{h}_{\ell,u}=\frac{1}{\lambda^{M}_{\ell}}+\frac{1}{\lambda^{M}_{\ell+1}}+\dots+\frac{1}{\lambda^{M}_{u-2}}+\frac{1}{\lambda^{M}_{u-1}},\;\; \invbreve{h}_{\ell,u}=\frac{1}{\lambda^{m}_{\ell}}+\frac{1}{\lambda^{m}_{\ell+1}}+\dots+\frac{1}{\lambda^{m}_{u-2}}+\frac{1}{\lambda^{m}_{u-1}},\nonumber
	\end{aligned}
\end{equation}
respectively. Since $\lambda^{M}_z=\max_{\widecheck{x}\in \mathcal{L}_z}\lambda_z(\widecheck{x}) = \max_{\widecheck{x}\in \mathcal{L}_z} \left( f_1(\widecheck{x})+f_2(\widecheck{x})+f_3(\widecheck{x})+f_4(\widecheck{x})\right)$, and $\lambda^{m}_z=\min_{\widecheck{x}\in \mathcal{L}_z}\lambda_z(\widecheck{x})=\min_{\widecheck{x}\in \mathcal{L}_z}\left( f_1(\widecheck{x})+f_2(\widecheck{x})+f_3(\widecheck{x})+f_4(\widecheck{x})\right)$, we observe that lower and upper bounds $\Breve{h}_{\ell,u}$ and $\invbreve{h}_{\ell,u}$ for the MFPT under consideration decrease as any of the rate constants $k_1,k_2,k_3,k_4$ increases. This result suggests that, if any reaction associated with either pathway in the bi-parallel network becomes faster, then the upper and lower bounds for the mean time needed to transform all J into all W (i.e., the MFPT from $n_{\mathrm{J}}=\mathrm{S_{tot}}$ to $n_{\mathrm{W}}=\mathrm{S_{tot}}$) decrease.

We note that, by following the same approach used above, these results can be generalized to any network having two species (J and W) connected by n-parallel pathways (beyond just two).
%%%%%%%%%%%%%%%%%%%%%%%%%%%%%%%%%%%%%%%%%%%%%%%%%%%%%%%

%\section{Methods}\label{sec11}
%
%Topical subheadings are allowed. Authors must ensure that their Methods section includes adequate experimental and characterization data necessary for others in the field to reproduce their work. Authors are encouraged to include RIIDs where appropriate. 

%%%%%%%%%%%%%%%%%%%%%%%%%%%%%%%%%%%%%%%%%%%%%%%%%%%%%%%%

\section{Conclusion}
\label{sec:conclusion}
	
%summary of results

In this paper, we started by providing a description of Stochastic Chemical Reaction Networks (SCRNs), which are a class of continuous time Markov chain models commonly used to describe the stochastic behavior of chemical reaction systems (Section \ref{sec:BasicDefinitions}).
We then introduced the notion of \textit{\textcolor{black}{coclique} level structure} (Section \ref{sec:IndLevelStructureApp}) and developed theoretical tools for identifying such \textcolor{black}{coclique} level structures for continuous time Markov chains associated with SCRNs, where each reaction involves the consumption of one molecule of a given species and the production of one molecule of another species, and an associated graph $\G$ is weakly connected (Section \ref{suffcond}). Additionally, we provided conditions for identifying if a SCRN does or does not admit \textcolor{black}{a coclique} level structure (Section \ref{suffcond}).
Finally, we derived analytical expressions for upper and lower bounds for  MFPTs of  SCRNs having \textcolor{black}{a coclique} level structure (Section \ref{indeplevelstructure}).

Following this, we provided illustrative examples to demonstrate the utility of our theoretical tools in studying the stochastic behavior of SCRNs (Section \ref{sec:Applications}). More precisely, we focused on models describing the main interactions among histone modifications alone, and in combination with DNA methylation \cite{bib:BWD2022}, as well as on a bi-parallel network motif, a typical building block found in complex networks \cite{RMilo2002}. Through these examples, we demonstrated that our algorithm for identifying \textcolor{black}{coclique} level structures is easy to apply, and the analytical expressions for upper and lower bounds of MFPTs obtained with our theoretical tools provide mechanistic insights into how system parameters affect the stochastic behavior of SCRNs.
This mechanistic insight is particularly valuable for applications where it is crucial to understand which biological parameters must be tuned to modulate system dynamics in a specific manner.

The mathematical results and theoretical tools developed in this paper can be applied to all stochastic models that meet the considered assumptions. Future work will be focused on generalizing these results by relaxing some of these assumptions, such as allowing the Markov chain to have countably many states and allowing SCRNs to have more than one molecule consumed and produced per reaction.
\textcolor{black}{Another valuable direction for future work is the implementation of the algorithm introduced in Section \ref{algIndLevFac}. While in this paper we focused on providing theoretical guidance, a practical implementation of the algorithm would broaden the applicability of our approach and facilitate its integration into simulation pipelines.
}

%%%%%%%%%%%%%%%%%%%%%%%%%%%%%%%%%%%%%%%%%%%%%%%%%%%%

\backmatter

\bmhead{Supplementary information}

%{\bf Supplementary information (SI) file:}

{\color{black} File containing proofs of a supporting lemma and one theorem, as well as the mathematical derivation of the analytical expression for the mean first passage times of a one-dimensional finite state birth-death process. Statements of two theorems from reference \cite{Monotonicitypaper} are also provided.}

\bmhead{Acknowledgements}

The work of the first and fourth authors was partially supported by NSF Collaborative Research grant MCB-2027949 (PI: D.D.V.). The work of the second, third, and fifth authors was partially supported by NSF Collaborative Research grant MCB-2027947 (PI: R.J.W.) and by the Charles Lee Powell Foundation (PI: R.J.W.).
Furthermore, the version of record of this article, first published in \textit{Journal of Mathematical Biology}, is available online at Publisher’s website: \href{https://doi.org/10.1007/s00285-025-02261-6}{https://doi.org/10.1007/s00285-025-02261-6}.\\ %We are grateful to the anonymous referees for several very helpful comments.\\ 

\section*{Declarations}

%Some journals require declarations to be submitted in a standardised format. Please check the Instructions for Authors of the journal to which you are submitting to see if you need to complete this section. If yes, your manuscript must contain the following sections under the heading `Declarations':

\begin{itemize}
\item Conflict of interest/Competing interests: the authors declare that they have no conflicts of interest.
\item Ethical statement: this research was conducted in adherence to the highest standards of academic integrity and ethical conduct. The research did not involve any studies with human participants or animals.
%\item Consent for publication
\item Data availability: data sharing not applicable to this article as no datasets were generated or analysed during the current study.
\end{itemize}

%\noindent
%If any of the sections are not relevant to your manuscript, please include the heading and write `Not applicable' for that section. 

	\newpage

\appendixpageoff
\appendixtitleoff
\renewcommand{\appendixtocname}{Supplementary Information}
\begin{appendices}
  \setcounter{section}{19}
  \crefalias{section}{supp}
  \setcounter{figure}{0}
\renewcommand{\thefigure}{S.\arabic{figure}}

  \setcounter{equation}{0}
\renewcommand{\theequation}{S.\arabic{equation}}

\pagenumbering{arabic}% resets `page` counter to 1
\renewcommand*{\thepage}{\arabic{page}}

\begin{center}
{\LARGE \textcolor{black}{Coclique} Level Structure for Stochastic Chemical Reaction Networks}\\
\end{center}

\begin{center}
{\large Simone Bruno$^{1,2,*}$, Yi Fu$^{3,*}$, Felipe A. Campos$^{4}$, Domitilla Del Vecchio$^{2}$, and Ruth J. Williams$^{4}$}\\
\end{center}

\begin{center}
$^1$\textit{\small Department of Data Science, Dana-Farber Cancer Institute, 450 Brookline Avenue, Boston, MA 02115. Email: {\tt\small sbruno@ds.dfci.harvard.edu}}

$^{2}$\textit{\small Department of Mechanical Engineering, Massachusetts Institute of Technology, 77 Massachusetts Avenue, Cambridge, MA 02139. Email: {\tt\small ddv@mit.edu}}

$^{3}$\textit{\small Bioinformatics and Systems Biology Program, University of California, San Diego, 9500 Gilman Drive, La Jolla CA 92093-0112. Email: {\tt\small yif064@ucsd.edu}}

$^{4}$\textit{\small Department of Mathematics, University of California, San Diego, 9500 Gilman Drive, La Jolla CA 92093-0112. Email: {\tt\small (fcamposv,rjwilliams)@ucsd.edu}}

$^{*}$\textit{\small These authors contributed equally: S. Bruno and Y. Fu}
\end{center}

\section*{Supplementary Information (SI)}

% {\color{black}
% \subsection{Example of a system whose associated graph contains more than one weakly connected component}
% \label{SI:disconnected-comp-example}

% The following example shows why, in this paper, we focus on SCRNs whose associated graphs $\G$ are weakly connected. More general systems with disconnected graphs can be treated as the union of their weakly connected components, each of which can be analyzed separately.

% Consider the following SCRN with three species ($A$, $B$, and $C$) and three reactions: 
% \begin{equation}\label{reacs-disconnected}
% \begin{aligned}
% &{\large \textcircled{\small 1}}\;\ce{A ->[$\alpha$] B}, \quad
% {\large \textcircled{\small 2}}\;\ce{B ->[$\beta$] A}, \quad
% {\large \textcircled{\small 3}}\;\ce{C ->[$\gamma$] D},
% \end{aligned}
% \end{equation}
% %
% in which $\alpha, \beta, \gamma >0$. The associated reaction vectors are: $v_1 = (-1,1,0,0)^T, v_2 = (1,-1,0,0)^T, v_3 = (0,0,-1,1)^T$. The associated graph $\G$ can be represented as follows:
% \begin{equation}
% \text{(1)}\; A \rightleftarrows B, \qquad \text{(2)}\; C \rightarrow D.
% \end{equation}

% In this case, $\G$ has two weakly connected components: component (1) involving species $A$ and $B$, and component (2) involving species $C$ and $D$. There are no reactions connecting the species in the two components. Therefore, the dynamics of $\{A, B\}$ evolve separately from those of $\{C, D\}$. Studying the full system does not provide additional insight beyond that obtained by analyzing the two subsystems separately.}

\subsection{Proof of Lemma \ref{lemmabis}}
\label{SI:proofs}

Before introducing the proof of Lemma \ref{lemmabis}, we provide some definitions used in the proof.
The \textbf{degree} of a vertex is the number of edges that are incident to the vertex, where we count both incoming and outgoing edges. A \textbf{weakly directed tree} is a directed graph whose underlying undirected graph is a tree. Given a directed graph with $d$ vertices, a \textbf{weakly directed spanning tree} is a subgraph of the graph with all $d$ vertices and such that it is a weakly directed tree, and it will have $d-1$ edges. We abbreviate weakly directed spanning tree as wd-spanning tree.
 Please note that a weakly connected graph $\G$ has a wd-spanning tree.

\begin{proof}
Consider a wd-spanning tree of the graph $\G$. Let $\widehat{\V}_{st}=\{\widehat v_1,\widehat v_2,\dots,\widehat v_{d-1}\}$ denote the set of $d-1$ reaction vectors associated with the wd-spanning tree edges.  We first show that the stoichiometric matrix associated with $\widehat{\V}_{st}$, $S_{d,d-1}=[\widehat v_1,\widehat v_2,\dots,\widehat v_{d-1}]\in \Z^{d \times (d-1) }$ has rank $d - 1$. We prove this by induction. First consider $d$ = 2. In this case the only possible wd-spanning tree is given by the two vertices and one edge connecting them. This means that $\rank(S_{2,1})$ = 1. Then, let us assume that the result is true for any wd-spanning tree with $d-1$ vertices (i.e., $\rank (S_{d-1,d-2})$ = $d-2$) for some $d \ge 3$, and consider a wd-spanning tree with $d$ vertices. This wd-spanning tree always has a degree-one vertex, that we define as $vertex\;1$, and then we can rearrange the rows and columns of $S_{d,d-1}$ so that %such that the first column is associated with this vertex, obtaining:
	\begin{equation*}%\label{matrix}
	S_{d,d-1}=\begin{bmatrix}
	\pm 1 & 0\\
	* & S_{d-1,d-2}
	\end{bmatrix}.
	\end{equation*}
	Then, $S_{d-1,d-2}$ is the stoichiometric matrix associated with the wd-spanning tree with the $vertex\;1$ and associated edge removed. By the induction hypothesis, we can conclude that $S_{d-1,d-2}$ has rank $d-2$, and then $\rank (S_{d,d-1})$ = $d-1$. 
  Since $S_{d,d-1}$ is obtained by removing columns associated with the reaction vectors $v_k\in \V \setminus \widehat{\V}_{st}$ from $S$, this result implies that $\rank(S) \ge d-1$. 
  
  Under Assumption \ref{assumption:Unimolecular_change}, $\one^Tv_k = 0$ for $k=1,\dots,n$, and thus $\one^TS=0$. This means that $\one\in \ker (S^T)$ and then $\dim \ker (S^T) \ge 1$. By the rank-nullity theorem (see for example \cite{LinearAlgebra}), we have that $\rank(S)= \rank(S^T) = d - \dim \ker (S^T) \le d-1$. Putting together the results obtained, we can conclude that $\rank(S) = d-1$.
  
  Furthermore, given that $\one^T S=0$ and $\dim \ker (S^T)=1$,
%. Thus, $\one$ is a conservation vector. Since  by Lemma \ref{lemmabis}, the nullspace of $S^T$ is one-dimensional by the rank-nullity theorem \cite{LinearAlgebra}. 
then $\one$ is the only (up to scalar multiplication) conservation vector such that $S^T\one=0$. 
%We conclude that $\one$ is the unique conservation vector of the SCRN considered (up to scalar multiplication).
\end{proof}

{\color{black}
\subsection{Proof of Theorem \ref{thm:Gen}}
\label{sec:thmGenProof}

        Let $\widecheck{S}^\wccIter \in \Z^{(d_\wccIter-1) \times n}$ be the first $(d_\wccIter-1)$ rows of the stoichiometric matrix $S^\wccIter$. Then, similar to \eqref{systMATRIX} and since the stoichiometric matrix for the SCRN has the form \eqref{eqn:StoiMat}, the system \eqref{syst2_paper} can be re-written in matrix-vector form as
        \begin{equation}
        \label{eqn:StoiMatProjected}
            \begin{bmatrix}
                \left( \widecheck{S}^1 \right)^T & 0 & 0 \\
                0 & \ddots & 0 \\
                0 & 0 & \left( \widecheck{S}^\wccNum \right)^T
            \end{bmatrix} b = w,
        \end{equation}
        which is equivalent to
        \begin{equation*}
            \left( \widecheck{S}^\wccIter \right)^T b^\wccIter = w^\wccIter \quad \text{ for every } \wccIter = 1, \dots, \wccNum \text{ such that } |\G^\wccIter| > 1,
        \end{equation*}
        where for $\wccIter = 1, \dots, \wccNum$, the vectors $b^\wccIter, w^\wccIter \in \Z^{d_\wccIter-1}$ are the $\wccIter^{th}$ entries of $b = (b^1, \dots, b^\wccNum)^T$ and $w = (w^1, \dots, w^\wccNum)^T$, respectively. 
        Note that for $\wccIter = 1, \dots, \wccNum$, if $|\G^\wccIter| = 1$, then $\widecheck{S}^\wccIter \in \Z^{(d_\wccIter-1) \times n}$ is a $0 \times n$ matrix, which does not appear in \eqref{eqn:StoiMatProjected}, and $\widecheck{x}^\wccIter$ is a zero-dimensional vector, and so we do not need to consider 
        %the $\wccIter^{th}$ weakly connected 
        that component in any coclique level function, as in \eqref{eqn:genCLF}.
        Consider $\wccIter = 1, \dots, \wccNum$ such that $|\G^\wccIter| > 1$. As noted before Theorem \ref{thm:Gen}, there is a SCRN associated with each $\G^\wccIter$. The stoichiometric matrix for this SCRN is $S^\wccIter$. By Theorem \ref{thm2b}, $\left( \widecheck{S}^\wccIter \right)^T b^\wccIter = w^\wccIter$ has a solution $b^\wccIter \in \Z^{d_\wccIter-1}$ if and only if $L^\wccIter: \Z^{d_{\wccIter}-1} \rightarrow \Z$ given by $L^\wccIter (\widecheck{x}^\wccIter) = (b^\wccIter)^T \widecheck{x}^\wccIter$ is a coclique level function for the SCRN associated with $\G^\wccIter$. Thus, 
        $L$ is a coclique level function for the whole SCRN if and only if
        \begin{equation*}
            L(\widecheck x) = b^T \widecheck x = (b^1,\dots,b^\wccNum)^T (\widecheck x^1, \dots, \widecheck x^\wccNum) = \sum_{\substack{\wccIter=1,\dots,\wccNum: \\ |\G^\wccIter|>1}} \left(b^\wccIter \right)^T \widecheck x^\wccIter = \sum_{\substack{\wccIter=1,\dots,\wccNum: \\ |\G^\wccIter|>1}} L^\wccIter (\widecheck x^\wccIter),
        \end{equation*}
        where $L^\wccIter$ is a coclique level function for the SCRN associated with $\G^\wccIter$ where $|\G^\wccIter| > 1$.
        
        For each weakly connected component $\G^q$ such that $|\G^q|>1$, by Theorem \ref{thm}, there exists a coclique level structure for the SCRN associated with $\G^\wccIter$ if and only if $\G^q$ is bipartite. Since $\G$ is bipartite if and only if each $\G^q$ with $|\G^q|>1$ is bipartite, we conclude that a coclique level function for $\widecheck{X}$ exists if and only if $\G$ is bipartite.
        \qed
        %
        %which completes the proof.
        %
        %We note that the system \eqref{syst2_paper} can be written in matrix-vector form as $\widecheck{S}^T b = w$ where $\widecheck{S}$ is obtained by taking rows of the stoichiometric matrix $S$ given in \eqref{eqn:StoiMat} and $S$ has a block-matrix structure. Therefore, the set of all coclique level functions for $\widecheck X$ may be obtained by looking at each weakly connected component separately, finding a coclique function for each component and then combining these.

}
%%%%%%%%%%%%%%%%%%%%%%%%%%%%%%%%%%%%%%%%%%%%%%%%%%

 \subsection{One-dimensional birth-death process: Mean first passage time.}
 \label{sec:Appendix1Dmodel}

   Let us consider a one-dimensional irreducible finite state continuous time Markov chain in which the state space $\X =\{0,1,\ldots, \upb \}$ and the off-diagonal entries of the infinitesimal generator $Q$ are all zero except for the following positive rates:
    \begin{equation}
    \label{1DQmatrix}
     \begin{aligned}
     Q_{x,x+1} &= \lambda_x  & \text{ if } x \in \{0,\ldots,\upb-1\}, \\
     Q_{x,x-1} &= \gamma_x  & \text{ if } x \in \{1,\ldots,\upb\}.       
    \end{aligned}   
    \end{equation}
In other words, $Q$ is the infinitesimal generator for a finite state birth-death process.

We will determine an analytical expression for the MFPT from $x=\upb$ to $x=0$ and from $x=0$ to $x=\upb$ for this chain. 
To this end, it is important to note that $X$ can be equivalently characterized using holding times with exponential parameters $\{q_x\}_{x \in \mathcal{X}}$ and a transition probability matrix $P$ for the \textit{embedded discrete time Markov chain}. More precisely, for each $x \in \mathcal{X}$, $q_x = - Q_{x,x} \neq 0$, since $X$ is irreducible, and for all $x,y \in \X$, $P_{x,x} = 0$, $P_{x,y} = \frac{Q_{x,y}}{q_x}$, for $y \neq x$ in $\X$. Note that $Q=\diag(q)(P-I)$. Defining $\B$ as a nonempty subset of $\widecheck{\X}$ such that $\B \neq \widecheck{\X}$ and using first step analysis (see (3.1) in \cite{Norris}), we obtain that the MFPT from $x$ to $\B$ can be written as
	\begin{equation}
	h_{x,\B} =\begin{cases}\label{MFPTsystem}
	0 & \mbox{if } x \in \B \\ 
	\frac{1}{q_x}+\sum_{y\in \mathcal{X}}P_{x,y} h_{y,\B}  & \mbox{if } x \in \B^c.
	\end{cases}
	\end{equation}
 
 Now, let us first focus on the MFPT from $x=\upb$ to $x=0$. In this case $\B=\{0\}$ and then \eqref{MFPTsystem} can be rewritten as
 %
 \iffalse
 We first focus on the former. For this, we exploit first step analysis (see Equation 3.1 of \cite{Norris2}), proceeding in a similar manner to that for \eqref{MFPTsystem}, to obtain
\fi
%
 \begin{equation}
    \label{1stepanalysisbis}
	\begin{cases}
	h_{0,0}=0, \\ 
	h_{x,0}=\frac{1}{\lambda_x+\gamma_x}+\frac{\lambda_x}{\lambda_x+\gamma_x}h_{x+1,0}+\frac{\gamma_x}{\lambda_x+\gamma_x}h_{x-1,0} & \mbox{if } x\in \{1,\ldots,\upb-1\}, \\
	h_{\upb,0}=\frac{1}{\gamma_{\upb}}+h_{\upb-1,0},\\
	\end{cases}
	\end{equation}
 where for $x,y \in \X$, $h_{x,y}=\E_x [\tau_y]$, $\tau_y = \inf \{t \geq 0:\: X(t)=y\}$, $X$ is the continuous time Markov chain with infinitesimal generator given by \eqref{1DQmatrix}.
Now, defining $\Delta h_{x,x-1}=h_{x,0}-h_{x-1,0}$ for $x\in \{1,\dots,\upb\}$, we can rewrite (\ref{1stepanalysisbis}) in the following way:
	\begin{equation}\label{1stepanalysistris}
	\begin{cases}
	h_{0,0}=0, \\ 
	\Delta h_{x,x-1} = \frac{1}{\gamma_x}+\frac{\lambda_x}{\gamma_x}\Delta h_{x+1,x} & \mbox{if } x\in  \{1,\ldots,\upb-1\},\\
	\Delta h_{\upb,\upb-1}=\frac{1}{\gamma_{\upb}}.\\
	\end{cases}
	\end{equation}
 From (\ref{1stepanalysistris}), we have an explicit formula for $\Delta h_{\upb,\upb-1}$ and any $\Delta h_{x,x-1}$ can be expressed as a function of $\Delta h_{x+1,x}$. Furthermore, if we sum the $\Delta h_{x,x-1}$ for $x= 1,\ldots,\upb$, we obtain
	\begin{equation}\label{repTTML1D}
	\begin{aligned}
	h_{\upb,0}=h_{\upb,0}-h_{0,0}=\sum_{x=1}^{\upb}\left(\Delta h_{x,x-1}\right)&= \Delta h_{1,0} + \Delta h_{2,1} + \dots + \Delta h_{\upb-1,\upb-2} +\Delta h_{\upb,\upb-1}.
	\end{aligned}
	\end{equation}
	Thus, to evaluate the MFPT from $x=\upb$ to $x=0$, we can calculate $\Delta h_{x,x-1}$ for $x=\upb,\upb-1,\dots,1$ and then sum all of the terms. 
 %   Defining $r_j=\frac{\lambda_1 \lambda_2\dots\lambda_j}{\gamma_1 \gamma_2 \dots \gamma_j}$, for $j=1,\dots,\upb$, 
    We then obtain
\begin{align}
\notag h_{\upb,0} &=\frac{1}{\gamma_{\upb}}\left(1+\frac{\lambda_{\upb-1}}{\gamma_{\upb-1}}+\frac{\lambda_{\upb-1}\lambda_{\upb-2}}{\gamma_{\upb-1}\gamma_{\upb-2}}+\dots+\frac{\lambda_{\upb-1}\dots\lambda_{1}}{\gamma_{\upb-1}\dots\gamma_{1}}\right)\\
\notag &\qquad+\frac{1}{\gamma_{\upb-1}}\left(1+\frac{\lambda_{\upb-2}}{\gamma_{\upb-2}}+\frac{\lambda_{\upb-2}\lambda_{\upb-3}}{\gamma_{\upb-2}\gamma_{\upb-3}}+\dots+\frac{\lambda_{\upb-2}\dots\lambda_{1}}{\gamma_{\upb-2}\dots\gamma_{1}}\right)+\dots+\frac{1}{\gamma_{1}}\\
%&=\frac{r_{\upb-1}}{\gamma_{\upb}}\left(1+\sum_{i=1}^{\upb-1}\frac{1}{r_i}\right)+\sum_{i=2}^{\upb-1}\left[\frac{r_{i-1}}{\gamma_i}\left(1+\sum_{j=1}^{i-1}\frac{1}{r_j}\right)\right]+\frac{1}{\gamma_1}. 
\label{formulaDto0} &= {\color{black}\frac{1}{\gamma_{1}} + \sum_{i=2}^{\upb} \frac{1}{\gamma_{i}} \left( 1 + \sum_{j=1}^{i-1} \frac{\lambda_{j}\dots\lambda_{i-1}}{\gamma_{j}\dots\gamma_{i-1}} \right).}
\end{align}
With a similar procedure, we can obtain the MFPT from $x=0$ to $x=\upb$. More precisely, 
 %   defining $\tilde r_j=\frac{\gamma_{\upb-1}\gamma_{\upb-2}\dots\gamma_{\upb-j}}{\lambda_{\upb-1}\lambda_{\upb-2}\dots\lambda_{\upb-j}}$, 
we have
\begin{align}
\notag	h_{0,\upb}&=\frac{1}{\lambda_{0}}\left(1+\frac{\gamma_1}{\lambda_1}+\frac{\gamma_1\gamma_2}{\lambda_1\lambda_2}+\dots+\frac{\gamma_{1}\dots\gamma_{\upb-1}}{\lambda_{1}\dots\lambda_{\upb-1}}\right)\\
\notag	&\qquad+\frac{1}{\lambda_{1}}\left(1+\frac{\gamma_2}{\lambda_2}+\frac{\gamma_2\gamma_3}{\lambda_2\lambda_3}+\dots+\frac{\gamma_{2}\dots\gamma_{\upb-1}}{\lambda_{2}\dots\lambda_{\upb-1}}\right)+\dots+\frac{1}{\lambda_{\upb-1}}\\
%\notag	&=\frac{\tilde r_{\upb-1}}{\lambda_{0}}\left(1+\sum_{j=1}^{\upb-1}\frac{1}{\tilde r_i}\right)+\sum_{i=2}^{\upb-1}\left[\frac{\tilde r_{i-1}}{\lambda_{\upb-i}}\left(1+\sum_{j=1}^{i-1}\frac{1}{\tilde r_j}\right)\right]+\frac{1}{\lambda_{\upb-1}}. \\
\label{formula3} &={\color{black} \frac{1}{\lambda_{\upb-1}} + \sum_{i=0}^{\upb-2} \frac{1}{\lambda_{i}} \left( 1 + \sum_{j=i+1}^{\upb-1} \frac{\gamma_{i+1}\dots\gamma_{j}}{\lambda_{i+1}\dots\lambda_{j}} \right).}
\end{align}
	A more detailed derivation of the $h_{0,\upb}$ and $h_{\upb,0}$ is given in \cite{bib:BWD20222}.

%	Let us consider the one-dimensional continuous time Markov chain introduced in Section \ref{SS1}, with state space $\X =  \{0,1,\ldots,\mathrm{D_{tot}}\}$ and infinitesimal transition rates that can be written as in \eqref{1Drates}. Since all of the transition rates are $O(1)$, except for $\lambda^\eps_0$ and $\gamma^\eps_{\mathrm{D_{tot}}}$ which are $O(\eps)$, then both $h_{\mathrm{D_{tot}},0}(\eps)$ and $h_{0,\mathrm{D_{tot}}}(\eps)$ are $O(1/ \eps)$. This means that in the limit as $\eps \rightarrow 0$, $h_{\mathrm{D_{tot}},0}(\eps)$ and $h_{0,\mathrm{D_{tot}}}(\eps)$, which correspond to the time to memory loss of the repressed and active states, respectively, tend to infinity. Substituting parameters in \eqref{formulaDto0} and \eqref{formula3} yields \eqref{formulaDto0INTRO} and \eqref{formula3INTRO}, respectively.

%%%%%%%%%%%%%%%%%%%%%%%%%%%%%%%%%%%%%%%%%%%%%%%%%%%%%%

 %\subsection{2D Model: additional mathematical details}
 %\label{sec:Appendix2Dmodel}

 %%%%%%%%%%%%%%%%%%%%%%%%%%%%%%%%%%%%%%%%%%%%%%%%%%%%%

 %\subsection{3D Model: additional mathematical details}
 %\label{sec:Appendix3Dmodel}

 %%%%%%%%%%%%%%%%%%%%%%%%%%%%%%%%%%%%%%%%%%%%%%%%%%%%%

 %\subsection{4D Model: additional mathematical details}
 %\label{sec:Appendix4Dmodel}

 %%%%%%%%%%%%%%%%%%%%%%%%%%%%%%%%%%%%%%%%%%%%%%%%%%%%%

%%%%%%%%%%%%%%%%%%%%%%%%%%%%%%%%%%%%%%%%%%%%%%%%%%

{\color{black} 

 \subsection{Theorems 3.3 \textcolor{black}{and 3.4} from \cite{Monotonicitypaper2}}
 \label{compthmpaper_thm}

Let $A$ be an $m \times d$ matrix, with no rows identically zero,  and $K_A = \{ y \in \R^d : Ay \geq 0\}$. For $x,y \in \R^d$, we say that $x \preccurlyeq_A y$ whenever $A(y-x) \geq 0$.
For a non-empty set $\Gamma \subseteq \X \subseteq \Z_+^d$, we say that a set $\Gamma$ is increasing in $\X$ with respect to $\preccurlyeq_A$ if for every $x \in \Gamma$ and $y \in \X$, $x \preccurlyeq_A y$ implies that $y \in \Gamma$. Moreover, we say that a set $\Gamma \subseteq \X$ is decreasing in $\X$ with respect to $\preccurlyeq_A$ if for every $x \in \Gamma$ and $y \in \X$, $y \preccurlyeq_A x$ implies that $y \in \Gamma$.
Furthermore, for $x \in \R^d$, let $K_A +x = \{ y \in \R^d : x \preccurlyeq_A y \}$ and $\partial_i(K_A+x) := \{ y \in K_A +x : \inn{A_{i\bullet},y} = \inn{A_{i\bullet},x} \}$\footnote{\textcolor{black}{Here, for convenience of notation, let $A_{i\bullet}$ denote the row vector corresponding to the $i$-th row of $A$, for $1 \leq i \leq m$. In this article we will adopt the convention of considering the inner product $\inn{\cdot,\cdot}$ as a function of a row vector in its first entry and as a function of a column vector in the second entry. In particular, $\inn{A_{i\bullet},x} = \sum_{k=1}^d A_{ik}x_k$.}} for each $1 \leq i \leq m$. We can then characterize the boundary of $K_A + x$ as follows:
\begin{equation}
\label{eq:BoundaryKAx}
\partial(K_A +x) = \bigcup_{i=1}^{m}  \partial_i(K_A+x).  
\end{equation}
Finally, we introduce the concept of usual stochastic order $\preccurlyeq_{st}$ for two random variables $Y,Z$: we say that $Y$ is smaller than $Z$ in the usual stochastic order, that is, $Y \preccurlyeq_{st} Z$, if $F_Y(t) \geq F_Z(t)$ for every $t \in \R$, where $F_Y$ and $F_Z$ are the cumulative distribution functions for $Y$ and $Z$, respectively.
 %   Now, let us introduce the following theorem.

 %When there are multiple vectors $v_j$ with a common value for $Av_j$, the pointwise comparison in $j$, for $1 \leq j \leq n$, in conditions \eqref{eq:CouplingConditions_InnerProduct_I} and \eqref{eq:CouplingConditions_InnerProduct_II} in Theorem \ref{thm:InnerProductTheorem}, can be weakened. To this end, let us introduce 
 
    In the following theorem, we consider the set of distinct vectors $\{\eta^1,\dots,\eta^s\}$ formed by $Av_j$, for $1 \leq j \leq n$, where $s$ denotes the cardinality of this set, and we consider the subsets of indices
    \begin{equation*}
        G^{k} :=  \{ j  : 1\leq j \leq n \hbox{ and } Av_j =\eta^k \}, \quad \text{for } 1 \leq k \leq s.
    \end{equation*}    
    
    The following theorem applies even if $\X$ is countably infinite, although in this paper all of our state spaces are finite. \\

   \noindent \textbf{Theorem S.1 (immediate consequence of Theorem 3.3 in \cite{Monotonicitypaper2}\footnote{\textcolor{black}{Compared to Theorem 3.3 in \cite{Monotonicitypaper2}, Theorem S.1 includes some additional clarifications inserted in parentheses to improve clarity and completeness in the context of this paper.}}).} \textit{Consider a non-empty set $\X \subseteq \Z_+^d$, a collection of distinct vectors $v_1,\ldots,v_n$ in $\Z^d \setminus \{0\}$ and two collections of non-negative (intensity) functions on $\X$, $\rate=(\rate_1, \dots,\rate_n)$ and $\Breve{\rate}= (\Breve{\rate}_1, \dots,\Breve{\rate}_n)$ such that if $x +v_j \notin \X$, then $\rate_j(x) = \Breve{\rate}_j(x)= 0$, and assume the associated continuous time Markov chains (with intensity functions given by $\rate$ and $\Breve{\rate}$, respectively, for the transition directions $v_1,\dots,v_n$) do not explode in finite time. Consider a matrix $A \in \Z^{m \times d}$ with non-zero rows and suppose that both of the following conditions hold:
    \begin{enumerate}
        \item[(i)]
        For each $1\leq j \leq n$, the vector $Av_j$ has entries in $\{-1,0,1\}$ only.
        \item[(ii)]
        For each $x \in \X$, $1 \leq i \leq m$ and $y \in \partial_i(K_A+x) \cap \X$ we have that
        \begin{equation}
        \label{eq:CouplingConditionBIS32}
        \sum_{j \in G^{k}} \Breve{\rate}_j(y) \leq \sum_{j \in G^{k}} \rate_j(x), \quad \text{for each } k \text{ such that } \eta^k_i <0,
        \end{equation}
        and
        \begin{equation}
        \label{eq:CouplingConditionBIS31}
        \sum_{j \in G^{k} } \Breve{\rate}_j(y) \geq \sum_{j \in G^{k}} \rate_j(x), \quad \text{for each } k \text{ such that } \eta^k_i >0.
        \end{equation}
    \end{enumerate}
    Then, for each pair $\initialx,\initialxbreve \in \X$ such that $\initialx \preccurlyeq_A \initialxbreve$, there exists a probability space $(\Omega,\F,\PP)$ with two continuous time Markov chains $X = \{X(t) : t \geq 0\}$ and $\Breve{X}=\{\Breve{X}(t) : t \geq 0\}$ defined there, each having state space $\X \subseteq \Z^d_+$, with infinitesimal generators $Q$ and $\Breve{Q}$, associated with $\rate$ and $\Breve{\rate}$ respectively (as in \eqref{eq:TransitionMatrixQ}), with initial conditions $X(0)=\initialx$ and $\Breve{X}(0)=\initialxbreve$ and such that:
        \begin{equation}
        \label{eq:CouplingConditionBIS}
        \PP\left[X(t) \preccurlyeq_A \Breve{X}(t) \text{ for every } t \geq 0 \right]=1.   
        \end{equation}}\\

   \noindent \textbf{Theorem S.2 (immediate consequence of Theorem 3.4 in \cite{Monotonicitypaper2}\footnote{\textcolor{black}{Compared to Theorem 3.4 in \cite{Monotonicitypaper2}, Theorem S.2 includes some additional clarifications inserted in parentheses to improve clarity and completeness in the context of this paper.}}).} \textit{Consider a non-empty set $\X \subseteq \Z_+^d$, a collection of distinct vectors $v_1,\ldots,v_n$ in $\Z^d \setminus \{0\}$ and two collections of non-negative (intensity) functions on $\X$, $\rate=(\rate_1, \dots,\rate_n)$ and  $\Breve{\rate}= (\Breve{\rate}_1, \dots,\Breve{\rate}_n)$, such that if $x +v_j \notin \X$, then $\rate_j(x) = \Breve{\rate}_j(x)= 0$, and assume the associated continuous time Markov chains (with intensity functions given by $\rate$ and $\Breve{\rate}$, respectively, for the transition directions $v_1,\dots,v_n$) do not explode in finite time. Consider a matrix $A \in \Z^{m \times d}$ with non-zero rows 
   %and suppose that the matrix $A$ has integer-valued entries 
   and conditions $(i)$ and $(ii)$ in Theorem S.1 are satisfied.}
   
    \textit{Let $\initialx,\initialxbreve \in \X$ be such that $\initialx \preccurlyeq_A \initialxbreve$ and let $X = \{X(t) : t \geq 0\}$ and $\Breve{X}=\{\Breve{X}(t) : t \geq 0\}$ be two continuous time Markov chains (possibly defined on different probability spaces), each having state space $\X \subseteq \Z^d_+$, with infinitesimal generators $Q$ and $\Breve{Q}$, associated with $\rate$ and $\Breve{\rate}$ respectively, and with initial conditions $X(0)=\initialx$ and $\Breve{X}(0)=\initialxbreve$. For a non-empty set $\Gamma \subseteq \X$, consider $T_{\Gamma} := \inf\{ t \geq 0 : X(t) \in \Gamma \}$ and $\Breve{T}_{\Gamma} := \inf\{ t \geq 0 : \Breve{X}(t) \in \Gamma \}$. If $\Gamma$ is increasing in $\X$ with respect to the relation $\preccurlyeq_A$, then
        \begin{equation}
        \label{eq:StochasticOrderingFPT_increasing}
            \Breve{T}_{\Gamma} \preccurlyeq_{st} T_\Gamma,
        \end{equation}
   and the mean first passage time of $\Breve{X}$ from $\initialxbreve$ to $\Gamma$ is dominated by the mean first passage time of $X$ from $\initialx$ to $\Gamma$. If $\Gamma$ is decreasing in $\X$ with respect to the relation $\preccurlyeq_A$, then
     \begin{equation}
      \label{eq:StochasticOrderingFPT_decreasing}
           T_\Gamma \preccurlyeq_{st} \Breve{T}_{\Gamma},
    \end{equation}
    and the mean first passage time of $X$ from $\initialx$ to $\Gamma$ is dominated by the mean first passage time of $\Breve{X}$ from $\initialxbreve$ to $\Gamma$.}\\

    The proof of these theorems can be found in Sections 5.3 and 3.3 of \cite{Monotonicitypaper2}, respectively.}

\vfill\break

\end{appendices}



\clearpage

	
\begin{thebibliography}{9}


% \bibitem{Allis}
% {\sc Allis, C. D., Caparros, M.-L., Jenuwein, T. and  Reinberg, D.} (2015).
% {\em Epigenetics, 2nd Edition.} Cold Spring Harbor Laboratory Press.


% \bibitem{Kurtz:11}
% {\sc Anderson, D.F. and Kurtz, T.G.} (2011).
% {\em Continuous Time Markov Chain Models for Chemical Reaction Networks.} Design and Analysis of Biomolecular Circuits: Engineering Approaches to Systems and Synthetic Biology, pp. 3-42.


% \bibitem{AAN2004}
% {\sc Altman, E. and Avrachenkov, K.E. and Núñez-Queija, R.} (2004).
% {\em Perturbation analysis for denumerable markov chains with application to queueing models} Advances in Applied Probability, vol. 36(3): 839–853.


\bibitem{AndersonKurtzBook}
{\sc Anderson, D.F. and Kurtz, T.G.} (2015).
{\em Stochastic Analysis of Biochemical Systems.} Springer International Publishing.

%\bibitem{AsmussenQueues2nd}
%{\sc Asmussen, S.} (2003).
%{\em Applied Probability and Queues, 2nd edition.} Stochastic Modelling and Applied Probability, Springer, New York, NY.


% \bibitem{AvrachenkovHaviv}
% {\sc Avrachenkov K.E. and Haviv, M.} (2004).
% {\em The first laurent series coefficients for singularly perturbed stochastic matrices.} Linear Algebra and Its Applications, vol. 386: 243–259.


% \bibitem{AvrachenkovFilarHowlett}
% {\sc Avrachenkov, K.E., Filar, J.A. and Howlett, P.G.} (2013).
% {\em Analytic Perturbation Theory and Its Applications.} SIAM, Philadelphia.

\bibitem{MBackenkohler2020}
{\sc Backenk{\"o}hler, M., Bortolussi, L., Wolf, V.} (2020).
{\em Bounding Mean First Passage Times in Population Continuous-Time Markov Chains} QEST 2020: 17th International Conference on Quantitative Evaluation of Systems, vol. 12289: 155–174.






% \bibitem{BeltranLandim}
% {\sc Beltrán J., Landim C.} (2010).
% {\em Tunneling and Metastability of Continuous Time Markov Chains} Journal of Statistical Physics, vol. 140(6): 1065–1114.

% \bibitem{BeltranLandim2}
% {\sc Beltrán J., Landim C.} (2012).
% {\em Tunneling and metastability of continuous time Markov Chains II, the nonreversible case} Journal of Statistical Physics, vol. 149(4): 598–618.

%\bibitem{EBertram1983}
%{\sc Bertram, E. A.} (1983).
%{\em Some applications of graph theory to finite groups} Discrete Mathematics, vol. 44(1): 31–43.



% \bibitem{BlumenthalGetoor}
% {\sc Blumenthal, R.M. and Getoor, R.K.} (1968).
% {\em Markov Processes and Potential Theory, 2nd edition.} Academic Press, New York and London.

\bibitem{FBrglez1989}
{\sc Brglez, F., Bryan, D., Koiminski, K.} (1989). 
Combinational profiles of sequential benchmark circuits. {\em IEEE International Symposium on Circuits and Systems}, Public Library of Science, vol. 3: 1929 - 1934.

\bibitem{epifinite}
{\sc Bruno, S., Campos, F.A., Fu, Y., Del Vecchio, D., Williams, R.J.} (2024).
{\em Analysis of Singularly Perturbed Stochastic Chemical Reaction Networks Motivated by Applications to Epigenetic Cell Memory}. SIAM Journal on Applied Dynamical Systems, {\bf 23}(4), 2695–2731.



\bibitem{bib:BWD2022}
{\sc Bruno, S., Williams, R.J. and Del Vecchio, D.} (2022). 
Epigenetic cell memory: The gene's inner chromatin modification circuit. {\em PLOS Computational Biology}, Public Library of Science, vol. 18(4): 1 - 27.


\bibitem{ECC2022}
{\sc Bruno, S., Williams, R.J. and Del Vecchio, D.} (2022). 
Model reduction and stochastic analysis of the histone modification circuit. {\em 2022 European Control Conference (ECC)}, 264-271.

\bibitem{BrunoSontagPaper2023}
{\sc Bruno, S., Williams, R.J. and Del Vecchio, D.} (2023). 
Mathematical analysis of the limiting behaviors of a chromatin modification circuit. {\em Math. Control Signals Syst}, vol. 35: 399–432.



\bibitem{Monotonicitypaper}
{\sc Campos, F.A., Bruno, S., Fu, Y., Del Vecchio, D. and Williams, R.J.} (2023). Comparison theorems for stochastic chemical reaction networks. {\em Bull Math Biol}, vol. 85(39).



% \bibitem{Carey2013}
% {\sc Carey, N.} (2013).
% {\em The epigenetic revolution.} Columbia University Press.



% \bibitem{Cooper}
% {\sc Cooper, G.M.} (2000).
% {\em The Cell: A Molecular Approach}. Sunderland (MA): Sinauer Associates.




% \bibitem{BFS}
% {\sc Del Vecchio, D. and Murray, R.M.} (2014).
% {\em Biomolecular Feedback Systems}. Princeton University Press.


\bibitem{2007DoddCellPaper}
{\sc Dodd, I.B. and Micheelsen, M.A. and Sneppen, K. and Thon, G.} (2007). Theoretical Analysis of Epigenetic
Cell Memory by Nucleosome Modification. {\em Cell}, vol. 129.

{\color{black}
\bibitem{Feinberg}
{\sc Feinberg, M.} (1986)
Chemical Reaction Network Structure and the Stability of Complex Isothermal Reactors - I. the Deficiency Zero and Deficiency One Theorems {\em Chemical Engineering Science, vol. 42(10): 2229-2268}
}

% \bibitem{Groudine03}
% {\sc Felsenfeld, G. and Groudine, M.} (2003). Controlling the double helix. {\em Nature}, vol. 421.



%\bibitem{LinearAlgebra}
%{\sc Friedberg, S. H., Insel, A. J.; Spence, L. E.} (2014).
%{\em Linear Algebra, 4th edition.} Pearson Education.

%\bibitem{gardiner1994}
%{\sc Gardiner, C.W.} (1994).
%{\em Handbook of stochastic methods for physics, chemistry and the natural sciences.} Springer-Verlag.

\bibitem{Latouche}
{\sc Gaver, D. P., Jacobs, P. A. and Latouche, G.} (1984)
{\em Finite birth-and death models in randomly changing environments.} Advances in Applied Probability, vol. 16: 715-731.


%\bibitem{Gi:07}
%{\sc Gillespie, D. T.} (2006).
%{\em Stochastic Simulation of Chemical Kinetics.} Annual Review of Physical Chemistry, vol. 58(1): 35-55.

%\bibitem{Grassmann771}
%{\sc Grassmann, W.K.} (1977)
%{\em Transient solutions in markovian queueing systems.} Computers \& Operations Research, vol. 4(1): 47–53.

%\bibitem{Grassmann772}
%{\sc Grassmann, W.K.} (1977)
%{\em Transient solutions in Markovian queues: An algorithm for finding them and determining their waiting-time distributions.} European Journal of Operational Research, vol. 1(6): 396–402.

% \bibitem{Hanna}
% {\sc Hanna, J. and Saha, K. and Pando, B. and van Zon, J. and Lengner, C.J. and Creyghton, M.P. and van Oudenaarden, A. and  Jaenisch, R.} (2009). Direct cell reprogramming is a stochastic
% process amenable to acceleration. {\em Nature.} 462.


\bibitem{HaseltineRawlings}
{\sc Haseltine, E.L., and  Rawlings J.B.} (2002). 
Approximate simulation of coupled fast and slow reactions for stochastic chemical kinetics.{\em The Journal of Chemical Physics}, 117:6959--6969.



% \bibitem{HH}
% {\sc Hassin, R. and Haviv, M.} (1992). Mean passage times and nearly uncoupled Markov chains. {\em SIAM Journal on Discrete Mathematics.} 5 (3), 386-397.


% \bibitem{EpiReview2015}
% {\sc Heard, E. and Martienssen, R.A.} (2014).
% Transgenerational Epigenetic
% Inheritance: Myths and Mechanisms. {\em Cell}, vol. 157.

%\bibitem{HornJohnson}
%{\sc Horn, R.A, Johnson, C.R.} (2013).
%{\em Matrix Analysis, 2nd Edition} Cambridge University Press.


% \bibitem{Blakey}
% {\sc Huang, S. and Litt, M. and Blakey,  C. A.} (2015).
% {\em Epigenetic Gene Expression and Regulation.} Academic Press.

%\bibitem{IsaacsonMadsen}
%{\sc Issacson, D.L. and Madsen R.W.} (1985)
%{\em Markov Chains: Theory and Applications} R.E. Krieger Pub. Co. 


\bibitem{KangKhudaBukhshKoepplRempala}
{\sc Kang, H. W., KhudaBukhsh, W.R., Koeppl, H., and Rempala, G.A.} (2019). 
Quasi-steady-state approximations derived from the stochastic model of enzyme kinetics. {\em Bulletin of Mathematical Biology}, 81: 1303–1336.

{\color{black}
\bibitem{KelseyRoneyDougal2022}
{\sc Kelsey, V., Roney-Dougal, C.M.} (2022).
{\em Maximal cocliques in the generating graphs of the alternating and symmetric groups}. Combinatorial Theory, {\bf 2}(1).
}

{\color{black}
\bibitem{bib:JE2005}
{\sc Kleinberg, J. and Tardos, E.},
(2005).
{\em Algorithm Design (first edition)}.
Pearson Education, Inc.
}
    
% \bibitem{KaratzasShreve}
% {\sc Karatzas, I. and Shreve, S.} (1991).
% {\em Brownian Motion and Stochastic Calculus. } New York, Springer‐Verlag. Graduate Texts in Mathematics 113.

%\bibitem{FrankKelly}
%{\sc Kelly, F.P.} (1979) 
%{\em Reversibility and Stochastic Networks.} Cambridge University Press.

%\bibitem{KemenySnell}
%{\sc Kemeny, J.G. and Snell, J.L.} (1960)
%{\em Finite Markov Chain} Van Nostrand Comp. Int., New York

%\bibitem{Liggett}
%{\sc Liggett, T.M.} (2010) 
%{\em Continuous Time Markov Processes: An Introduction.} Graduate Studies in Mathematics, Volume 113 . American Mathematical Society.


\bibitem{SLopatatzidis2017}
{\sc Lopatatzidis, S., De Bock, J., de Cooman, G.} (2017)
{\em Computing lower and upper expected first-passage and return times in imprecise birth–death chains.} International Journal of Approximate Reasoning, vol. 80(2): 137–173.

%\bibitem{MasseyStochasticOrdering}
%{\sc Massey, W.A.} (1987)
%{\em Stochastic Orderings for Markov Processes on Partially Ordered Spaces.} Mathematics of Operations Research, vol. 12(2): 350–367.

% \bibitem{Meyer1975}
% {\sc Meyer, Jr. and Carl, D.} (1975).
% The role of the group generalized inverse in the theory of finite Markov chains. {\em SIAM Review}, vol. 17, no. 3, pp. 443–464.


\bibitem{RMilo2002}
{\sc Milo, R., Shen-Orr, S., Itzkovitz, S., Kashtan, N., Chklovskii, D., Alon, U.} (2002).
Network motifs: simple building blocks of complex networks.
{\em Science,} 298(5594):824-7.


%\bibitem{MullerStoyan}
%{\sc Muller, A., Stoyan, D.} (2002).
%{\em Comparison Methods for Stochastic Models and Risks.} Wiley Series in Probability and Statistics.

\bibitem{Norris}
{\sc Norris, J.R.} (1997).
{\em Markov Chains.} Cambridge University Press.


%\bibitem{PalaciosBruno2024}
%{\sc Palacios, S. and Bruno, S. and Weiss, R. and Salibi, E. and Kane, A. and Ilia, K. and Del Vecchio, D.} (2024).
%Analog epigenetic cell memory by graded DNA methylation. {\em bioRxiv}, DOI: 10.1101/2024.02.13.580200.


% \bibitem{Pan1999}
% {\sc Pan, Z.G. and Basar, T.} (1999).
% H infintity control of large-scale jump linear systems via averaging and aggregation. {\em International Journal of Control}, vol. 72, no. 10, pp. 866–881.



% \bibitem{Phillips1981}
% {\sc Phillips, R.G. and Kokotovic, P.V} (1981).
% A singular perturbation approach to modelling and control of Markov chains. {\em IEEE Trans. Automat. Control}, vol. 26, pp. 1087–1094.


%\bibitem{PenroseGeneralizedInverse}
%{\sc Penrose, R.} (1955)
%{\em A generalized inverse for matrices.} Mathematical Proceedings of the Cambridge Philosophical Society, 51(3), 406-413. doi:10.1017/S0305004100030401


% \bibitem{Ptashne2013}
% {\sc Ptashne, M.} (2013).
% Epigenetics: Core misconcept. {\em Proc. Natl. Acad. Sci}, vol. 110.


%\bibitem{Ross1996}
%{\sc Ross, S.M.} (1996).
%{\em Stochastic processes.} Wiley Series in Probability and Mathematical Statistics.


\bibitem{Srivastava}
{\sc Srivastava, R., You, L., Summers, J. and Yin, J.} (2002).
Stochastic vs. deterministic modeling of intracellular viral kinetics.
{\em Journal of Theoretical Biology,} 218(3):309-321.


\bibitem{vanMieghem2010}
{\sc van Mieghem, P.} (2010).
{\em Graph Spectra for Complex Networks}. Cambridge University Press.



\bibitem{WilsonRJ}
{\sc Wilson, R.J.} (1996).
{\em Introduction to Graph Theory.} 4th ed., Addison Wesley Longman.

% \bibitem{Ferrell}
% {\sc Xiong, W.  and Ferrell Jr, J.E.} (2003).
% A positive-feedback-based
% bistable ‘memory module’ that
% governs a cell fate decision. {\em Nature}, vol. 426.



% \bibitem{yinzhang2013}
% {\sc Yin, G.G. and Zhang, Q.} (2013).
% {\em Continuous-Time Markov Chains and Applications. A Two-Time-Scale Approach, Second Edition.} Springer, New York, NY.


% \bibitem{H3K4meMaintenance2015}
% {\sc Zhang, T. and  Cooper, S. and  Brockdorff, N.} (2015).
% The interplay of histone modifications -- writers
% that read. {\em EMBO Reports}, vol. 16.





%\bibitem{penrose}
%{\sc Penrose, R.} (1955)
%{\em On Best Approximate Solutions of Linear Matrix Equations} Proceedings of the Cambridge Philosophical Society.

%\bibitem{uniformization}
%{\sc Grassmann,W.K. } (1977)
%{\em Transient solutions in markovian queueing systems.} Computers & Operations Research, vol. 4, no. 1 pp. 47-53.

\end{thebibliography}

\clearpage
    
\begin{thebibliography}{9}

%\bibitem{Allis}
%{\sc Allis, C. D. and Caparros, M.-L. and Jenuwein, T. and  Reinberg, D.} (2015).
%{\em Epigenetics, 2nd Edition.} Cold Spring Harbor Laboratory Press.


%\bibitem{Kurtz:11}
%{\sc Anderson, D. and Kurtz, T.} (2011).
%{\em Continuous Time Markov Chain Models for Chemical Reaction Networks.} Design and Analysis of Biomolecular Circuits: Engineering Approaches to Systems and Synthetic Biology, pp. 3-42.

%\bibitem{AndersonKurtzBook}
%{\sc Anderson, D.F. and Kurtz, T.G.} (2015).
%{\em Stochastic Analysis of Biochemical Systems.} Springer International Publishing.

% \bibitem{AsmussenQueues2nd}
% {\sc Asmussen, S.} (2003).
% {\em Applied Probability and Queues, 2nd edition.} Stochastic Modelling and Applied Probability, Springer, New York, NY.

% \bibitem{AvrachenkovFilarHowlett2}
% {\sc Avrachenkov K.E., Filar, J.A. and Howlett, P.G.} (2013).
% {\em Analytic Perturbation Theory and Its Applications.} SIAM, Philadelphia.

%\bibitem{BlumenthalGetoor}
%{\sc Blumenthal, R.M., Getoor, R.K.} (1968).
%{\em Markov Processes and Potential Theory, 2nd edition.} Academic Press, New York and London.

\bibitem{bib:BWD20222}
{\sc Bruno, S., Williams, R.J. and Del Vecchio, D.} (2022). 
{\em Epigenetic cell memory: The gene's inner chromatin modification circuit.} PLOS Computational Biology, Public Library of Science, vol. 18(4): 1 - 27.

%\bibitem{ECC2022}
%{\sc Bruno, S., Williams, R.J., Del Vecchio, D.} (2022) 
%{\em Model reduction and stochastic analysis of the histone modification circuit.} 2022 European Control Conference (ECC), London, United Kingdom, 264-271.


% \bibitem{Meyer19752}
% {\sc Meyer, Jr. and Carl, D.} (1975).
% The role of the group generalized inverse in the theory of finite Markov chains. {\em SIAM Review}, vol. 17, no. 3, pp. 443–464.

 \bibitem{Monotonicitypaper2}
 {\sc Campos, F.A., Bruno, S., Fu, Y., Del Vecchio, D. and Williams, R.J.} (2023). Comparison theorems for stochastic chemical reaction networks. {\em Bull Math Biol}, vol. 85(39).

 \bibitem{LinearAlgebra}
 {\sc Friedberg, S. H., Insel, A. J. and Spence, L. E.} (2014).
 {\em Linear Algebra, 4th edition.} Pearson Education.

% \bibitem{gardiner1994}
% {\sc Gardiner, C.W.} (1994).
% {\em Handbook of stochastic methods for physics, chemistry and the natural sciences.} Springer-Verlag.

% \bibitem{Latouche}
% {\sc Gaver, D. P., Jacobs, P. A. and Latouche, G.} (1984).
% Finite birth-and death models in randomly changing environments. {\em Advances in Applied Probability}, Vol. 16, 715-731.


%\bibitem{Gi:07}
%{\sc Gillespie, D. T.} (2006).
%{\em Stochastic Simulation of Chemical Kinetics.} Annual Review of Physical Chemistry, vol. 58, no. 1, pp. 35-55.

% \bibitem{Grassmann771}
% {\sc Grassmann, W.K.} (1977).
% Transient solutions in Markovian queueing systems. {\em Computers \& Operations Research}, vol. 4, no. 1, pp. 47–53.

%\bibitem{Grassmann772}
%{\sc Grassmann, W.K.} (1977)
%{\em Transient solutions in Markovian queues: An algorithm for finding them and determining their waiting-time distributions.} European Journal of Operational Research, vol. 1, no. 6, pp. 396–402.

%\bibitem{HH}
%{\sc Hassin, R. and Haviv, M.} (1992). {\em Mean Passage Times and Nearly Uncoupled Markov Chains.} SIAM Journal on Discrete Mathematics 5 (3), 386-397


% \bibitem{HornJohnson}
% {\sc Horn, R.A and Johnson, C.R.} (2013).
% {\em Matrix Analysis, 2nd Edition} Cambridge University Press.

%\bibitem{Blakey}
%{\sc Huang, S. and Litt, M. and Blakey,  C. A.} (2015).
%{\em Epigenetic Gene Expression and Regulation.} Academic Press.

% \bibitem{IsaacsonMadsen}
% {\sc Issacson, D.L. and Madsen R.W.} (1985).
% {\em Markov Chains: Theory and Applications} R.E. Krieger Pub. Co. 

%\bibitem{KaratzasShreve}
%{\sc Karatzas, I. and Shreve, S.} (1991).
%{\em Brownian Motion and Stochastic Calculus. } New York, Springer‐Verlag. Graduate Texts in Mathematics 113.

% \bibitem{James1978}
% {\sc James, M.} (1978). The generalised inverse. {\em The Mathematical Gazette,}
% 62 (420), 109-114.

% \bibitem{FrankKelly}
% {\sc Kelly, F.P.} (1979). 
% {\em Reversibility and Stochastic Networks.} Cambridge University Press.

% \bibitem{KemenySnell}
% {\sc Kemeny, J.G. and Snell, J.L.} (1960).
% {\em Finite Markov Chains.} Van Nostrand Comp. Int., New York

% \bibitem{Liggett}
% {\sc Liggett, T.M.} (2010). 
% {\em Continuous Time Markov Processes: An Introduction.} Graduate Studies in Mathematics, Volume 113 . American Mathematical Society.

%\bibitem{MasseyStochasticOrdering}
%{\sc Massey, W.A.} (1987)
%{\em Stochastic Orderings for Markov Processes on Partially Ordered Spaces.} Mathematics of Operations Research, vol. 12, no. 2, 1987, pp. 350–367.

%\bibitem{MullerStoyan}
%{\sc Muller, A., Stoyan, D.} (2002).
%{\em Comparison Methods for Stochastic Models and Risks.} Wiley Series in Probability and Statistics.

\bibitem{Norris2}
{\sc Norris, J.R.} (1997). 
{\em Markov Chains.} Cambridge University Press.

% \bibitem{PenroseGeneralizedInverse}
% {\sc Penrose, R.} (1955).
% A generalized inverse for matrices. {\em Mathematical Proceedings of the Cambridge Philosophical Society}, 51(3), 406-413. 
% %doi:10.1017/S0305004100030401

%\bibitem{Ross1996}
%{\sc Ross, S.M.} (1996).
%{\em Stochastic processes.} Wiley Series in Probability and Mathematical Statistics.



%\bibitem{yinzhang2013}
%{\sc Yin, G.G. and Zhang, Q.} (2013)
%{\em Continuous-Time Markov Chains and Applications. A Two-Time-Scale Approach, Second Edition.} Springer, New York, NY.

\end{thebibliography}

\end{document}